\documentclass[10pt,reqno]{article}
\usepackage{amsmath,amsthm,amssymb,mathrsfs,stmaryrd,color,upgreek}
\usepackage{mathtools}
\usepackage[all]{xy}
\usepackage{url}

\usepackage[utf8]{inputenc}
\usepackage{tikz-cd}
\usetikzlibrary{fit}

\usepackage{geometry}
\usepackage{array}

\usepackage{enumitem}

\usepackage{relsize}
\usepackage[bbgreekl]{mathbbol}
\usepackage{amsfonts}
\DeclareSymbolFontAlphabet{\mathbb}{AMSb}
\DeclareSymbolFontAlphabet{\mathbbl}{bbold}
\newcommand{\Prism}{{\mathlarger{\mathbbl{\Delta}}}}

\usepackage[colorlinks=true,hyperindex, linkcolor=magenta, pagebackref, citecolor=cyan]{hyperref}
\usepackage{cleveref}

\newtheorem{theorem}{Theorem}[subsection]
\newtheorem{proposition}[theorem]{Proposition}
\newtheorem{lemma}[theorem]{Lemma}
\newtheorem{corollary}[theorem]{Corollary}

\theoremstyle{definition}
\newtheorem{definition}[theorem]{Definition}
\newtheorem{question}[theorem]{Question}
\newtheorem{remark}[theorem]{Remark}
\newtheorem{example}[theorem]{Example}
\newtheorem{notation}[theorem]{Notation}
\newtheorem{warning}[theorem]{Warning}
\newtheorem{construction}[theorem]{Construction}
\newtheorem{observation}[theorem]{Observation}
\newtheorem{assumption}[theorem]{Assumption}

\newcommand{\dR}{\mathrm{dR}}
\newcommand{\an}{\mathrm{an}}
\newcommand{\Fun}{\mathrm{Fun}}
\newcommand{\Sp}{\mathcal{S}}
\newcommand{\HT}{\mathrm{HT}}
\newcommand{\crys}{\mathrm{crys}}
\newcommand{\calF}{\mathcal{F}}
\newcommand{\NilPerfd}{\mathrm{NilPerfd}}
\newcommand{\daggerred}{{\mathrm{\dagger}\mathrm{-red}}}
\newcommand{\Ga}{\mathbf{G}_a}
\newcommand{\CAlg}{\mathrm{CAlg}}
\newcommand{\Spf}{\mathrm{Spf}}
\newcommand{\Spec}{\mathrm{Spec}}
\newcommand{\Shv}{\mathrm{Shv}}
\newcommand{\Perf}{\mathrm{Perf}}
\newcommand{\colim}{\mathop{\mathrm{colim}}}
\newcommand{\Weil}{\mathrm{Weil}}
\newcommand{\proet}{\mathrm{pro\acute{e}t}}
\newcommand{\et}{\mathrm{\acute{e}t}}
\newcommand{\sm}{\mathrm{sm}}

\newcommand{\nilraddagger}[1]{\mathrm{Nil}^{\dagger}_{#1}}
\newcommand{\Rig}{\mathrm{Rig}}
\newcommand{\fSch}{{f\mathrm{Sch}}}
\newcommand{\calO}{\mathcal{O}}
\newtagform{magenta}{\color{magenta}(}{)}
\usetagform{magenta}

\begin{document}

\title{Non-abelian $p$-adic Hodge theory and stacks}
\author{Bhargav Bhatt}
\date{}

\maketitle
\begin{abstract}
In the first half of these notes, we survey some aspects of the stack-theoretic perspective on $p$-adic Hodge theory, with an emphasis on the non-abelian aspects. In the second half, we report ongoing joint work on ``quantizing'' non-abelian $p$-adic Hodge theory. 
\end{abstract}

\setcounter{tocdepth}{2}
\tableofcontents

\newpage 
\section{Introduction}

\subsection{Goals}

These notes concern $p$-adic Hodge theory, i.e., the study of $p$-adic cohomological invariants (or fancier categorical ones) attached to $p$-adic formal schemes (or rigid spaces). In recent years, the subject has undergone extensive development and re-organization. In particular, the emergent stack-theoretic viewpoint on $p$-adic cohomology theories has proven especially valuable for studying coefficient systems (e.g., in non-abelian $p$-adic Hodge theory). Through a handful of examples, we illustrate this perspective. These notes are not intended to be comprehensive or self-contained --- several statements are recorded without proof --- and the reader is invited to consult the references for further exploration.

\subsection{What is done here?}

Given a $p$-adic formal scheme $X$, the broad goal of the stack-theoretic approach is to produce a stack $X^?$ whose algebraic geometry captures some essential feature ``$?$'' of the $p$-adic Hodge theory of $X$. For example, when $X$ is a smooth variety in characteristic $p$, one might want a stack $X^\crys$ whose category of vector bundles is the category of crystals in vector bundles on $X$. This procedure --- identifying a category of interest as quasi-coherent sheaves on a stack --- is often informally called {\em geometrization}. To construct such geometrizations, the following general mechanism  has proven fruitful:

\begin{construction}[Transmutation]
\label{cons:Transmutation}
Fix  categories $\mathcal{C}$ and $\mathcal{E}$ equipped with forgetful functors to the category of commutative rings (e.g., they could simply be slice categories). Let $\mathcal{R}:\mathcal{C} \to \mathcal{E}$ be a functor, which we informally call a ring stack on $\mathcal{C}$.
 Given $X \in \Fun(\mathcal{E},\mathrm{Sets})$, such as the functor represented by an algebraic variety, we can define $X^{\mathcal{R}} := X \circ \mathcal{R}:\mathcal{C} \to \mathrm{Sets}$; this is called the {\em transmutation of $X$ with respect to $\mathcal{R}$}. (In the body of the text, it will be important to replace rings with animated rings to accommodate some features specific to the $p$-adic world; see Remark~\ref{rmk:CharpvsChar0deRhamRingStack}.)
\end{construction}

Construction~\ref{cons:Transmutation} is itself purely formal and quite trivial. The art in applying it lies in concocting suitable ring stacks. The first (and still the model) example of transmutation is the following:

\begin{example}[Simpson's de Rham stack \cite{SimpsonHomotopy,SimpsonTeleman}]
\label{ex:SimpsondR}
Fix a field $k$ of characteristic $0$. Let $\mathcal{C}=\mathcal{E}=\mathrm{CAlg}_k$ be the category of $k$-algebras, and let $\mathcal{R}:\mathcal{C} \to \mathcal{E}$ be $R \mapsto R_{\text{red}} = R/\mathrm{Nil}(R)$. The functor $\mathcal{R}$ is almost representable: we have $\mathcal{R} = \Ga/\widehat{\Ga}$, where $\mathbf{G}_a$ represents the identity functor on $\CAlg_k$, and $\widehat{\Ga} \subset \Ga$ is the nilradical subfunctor. Given a $k$-variety $X$, regarded as a functor on $\mathrm{CAlg}_k$, its transmutation $X^{\mathcal{R}}:\mathrm{CAlg}_k \to \mathrm{Sets}$ is given by $R \mapsto X(R_{\text{red}})$. This functor is called the {\em de Rham stack} of $X$ and denoted $X^\dR$. Transmuting the projection $R \to R_{\text{red}}$ for  $R \in \CAlg_k$ gives a natural map $\eta_X:X \to X^\dR$. By the infinitesimal lifting property, the map $\eta_X$ is surjective for $X$ smooth. Moreover, via pullback along $\eta_X$, the category $\mathrm{QCoh}(X^\dR)$ identifies with the category of quasi-coherent sheaves on $X$ equipped with an integrable connection, with the functor $R\Gamma(X^\dR,-)$ corresponding to de Rham cohomology. Thus, the single stack $\mathcal{R} = \Ga/\widehat{\Ga}$, with its $k$-algebra structure, encodes the theory of algebraic $\mathcal{D}$-modules  via transmutation. 
\end{example}

 In these notes, we shall apply transmutation to define and study stacks that similarly geometrize both de Rham  and prismatic cohomology of $p$-adic formal schemes, and discuss applications. These cohomology theories are ``less transcendental'' than de Rham cohomology in characteristic $0$. Consequently, the stack-theoretic perspective is especially powerful for studying them, arguably more so than in characteristic $0$. To justify this, we offer the following in these notes (with more detailed introductions given in each section):
\begin{itemize}
    \item In \S \ref{sec:padicdR} and \S \ref{sec:Prismatization}, we explain how the stack-theoretic approach quickly recovers many known results on de Rham cohomology and non-abelian Hodge theory (in both characteristic $p$ and mixed characteristic); this material is  expository. Moreover, the characteristic $p$ story in \S \ref{sec:padicdR} does not require much background beyond basic algebraic geometry and some familiarity with stacks. 

\item In \S \ref{sec:Simpson}, we review joint work with Mingjia Zhang \cite{MingjiaSimpsonNotes,Bhatt2025Apsects} on the $p$-adic Simpson equivalence for rigid spaces. Our presentation emphasizes the link to the prismatic material in \S \ref{sec:Prismatization}: the equivalence is obtained via a spectral decomposition of the category of Hodge--Tate crystals made possible by the geometry of the Hodge--Tate stack (namely, its structure as a torsor for a commutative group).

\item In \S \ref{sec:QuantNAHT}, we use the stack-theoretic perspective to ``quantize'' non-abelian Hodge theory: over certain so-called $\varphi$-crystalline prisms (Definition~\ref{def:phicrystallineprism}), prismatic crystals are identified with $I$-connections on a lift, up to an explicit gerbe twist (Theorem~\ref{thm:QuantizationPhiCrys}). Similarly, for smooth partially proper rigid spaces over a perfectoid field $C/\mathbf{Q}_p$, we sketch a $B_\dR^+$-deformation of the $p$-adic Simpson equivalence (Theorem~\ref{thm:QuantNAHTrigid}). In both cases, the main idea again is to spectrally decompose the relevant categories by describing the stacks underlying them as torsors for commutative groups. This is a report on joint work in progress with Heuer, Petrov and Vologodsky \cite{BHPVQuanNAHT}.
 \end{itemize}

The unifying theme of these notes, which guided both the choice of material and the style of exposition, is the utility of the stack-theoretic perspective: highly non-trivial theorems about general varieties reduce to rather simple statements about a few specific ring stacks.

\subsection*{Acknowledgements}

I am most grateful to Michael Barz, Guido Bosco, H\'{e}l\`{e}ne Esnault, Michael Groechenig, Max Hauck, Ben Heuer, Jacob Lurie, Akhil Mathew, Arthur Ogus, Alexander Petrov, Peter Scholze, Kenta Suzuki, Gleb Terentiuk, Vadim Vologodsky, Daxin Xu and Mingjia Zhang for conversations related to the material presented here. Additional gratitude is due to Michael Barz, Ben Heuer, Alexander Petrov, Kenta Suzuki, Vadim Vologodsky, and Mingjia Zhang for valuable comments on earlier versions of these notes. 

These notes were prepared to accompany my {\em Takagi Lectures} in June 2026\footnote{The June 2026 version of this document was written without AI tools;  the exposition in the current version, especially typesetting and the bibliography, benefitted from LLM assistance.}. I would like to thank the organizers for the opportunity to contemplate this material. 

I am grateful to the Packard and Simons Foundations for supporting the research presented here.

\newpage
\section{$p$-adic de Rham cohomology}
\label{sec:padicdR}

In this section, we explain the stack-theoretic perspective on de Rham cohomology in the $p$-adic world, giving an analog of Example~\ref{ex:SimpsondR} in this context. The idea is  simple: replace the ind-group $\widehat{\Ga}$ with the group scheme $\Ga^\sharp$ (the divided power hull of the origin in $\Ga$). This replacement is analogous to the passage from the infinitesimal site to the crystalline site in Grothendieck's development of crystalline cohomology \cite{GrothendieckCrystals}. The stack-theoretic point of view nonetheless brings in some intuitions that are hidden in the classical treatment. In this section, we justify this assertion via the following examples:
\begin{itemize}
    \item The Cartier isomorphism and Cartier descent become consequences of a gerbe structure on the de Rham stack (Corollaries~\ref{cor:Cartier} and \ref{cor:CartierDescent}); the latter is itself a consequence of deformation theory applied to a square-zero extension structure on $\Ga^\dR := \Ga/\Ga^\sharp$ (Proposition~\ref{prop:SqZerodRring}).

    \item The Deligne--Illusie theorem \cite{DeligneIllusie} on decompositions of the de Rham complex is a consequence of an interesting $\mu_p$-action on the ring stack $\Ga^\dR$ (Corollary~\ref{cor:DeligneIllusie}).

    \item The existence of crystalline cohomology is a  formal consequence of  the fact that the ideal $(p) \subset \mathbf{Z}_p$ has divided powers, i.e., that $p \in \Ga^\sharp(\mathbf{Z}_p) \subset \Ga(\mathbf{Z}_p) = \mathbf{Z}_p$ (Corollary~\ref{cor:CrysMiracle}).

    \item Some fundamental statements in non-abelian Hodge theory in characteristic $p$ (due to \cite{BMRLocalization,OgusVologodsky}) are encoded, via a twisted version of Cartier duality, in the aforementioned gerbe structure on the de Rham stack (Corollary~\ref{cor:OVexplicit}).
\end{itemize}

Some references for the material in this section are \cite{BarzLogarithmic,BhattLuriePrismatization,DrinfeldStacky,DrinfeldPrismatization,OgusVologodsky,PetrovOV,TerentiukOV}.

\begin{notation}
In this section, we fix a perfect field $k$, let $A=W(k)$, and let $\varphi:A \to A$ be the Frobenius automorphism. Write $\CAlg^\an_{\Spf(A)}$ for the $\infty$-category of $p$-nilpotent animated $A$-algebras, and $\CAlg^\an_k$ for the $\infty$-category of animated $k$-algebras. Writing $\Sp$ for the $\infty$-category of spaces, our geometric objects over $A$ will be regarded as objects of $\Fun(\CAlg^\an_{\Spf(A)},\Sp)$, and similarly over $k$. For any $X \in \Fun(\CAlg^\an_{\Spf(A)},\Sp)$, one defines its quasi-coherent derived category $\mathcal{D}(X)$ as $ \lim \mathcal{D}(R)$, where the limit runs over the category of points of $X$, i.e., the category of pairs $(R \in \CAlg^\an_{\Spf(A)},\alpha \in X(R))$.  We shall also use the following convention\footnote{Another, perhaps more common in the literature, name for such objects is a ``(smooth) formal scheme $X$ over $\Spf(A)$''. However, a literal parsing of this name includes many more examples than we care to allow, e.g., $\mathrm{Spf}(B)$ where $B=A\llbracket t \rrbracket$ given the $(p,t)$-adic topology.}: a (smooth) scheme over $X/\Spf(A)$ refers to (smooth) schematic map $X \to \Spf(A)$.

\end{notation}

\subsection{The group scheme $\Ga^\sharp$}

The group scheme $\Ga^\sharp$ is an essential player not only for de Rham cohomology but for the entire prismatization story. In this subsection, we introduce this scheme (Definition~\ref{def:Gasharp}), discuss its representation theory (Proposition~\ref{prop:reptheoryGasharp}), and give a compact presentation coming from the Witt vectors (Proposition~\ref{prop:WittGasharp}).

\begin{definition}
\label{def:Gasharp}
The group scheme $\Ga^\sharp$ is the PD-hull of $0$-section in $\Ga$ over $\Spf(A)$. Explicitly, it is the functor on $\CAlg^\an_{\Spf(A)}$ corepresented by the divided power polynomial ring $A[ \{\gamma_n(x)\}_{n \geq 1}]$, where $\gamma_n(x)$ is a formal symbol standing for $\frac{x^n}{n!}$. Write $d:\Ga^\sharp \to \Ga$ for the natural map forgetting divided powers.

More generally, for a vector bundle $E$ over some $A$-stack $X$, we shall also write $E^\sharp \to E$ for the PD-hull of the $0$-section, so $E^\sharp = E \otimes_{\Ga} \Ga^\sharp$. 
\end{definition}

The representation theory of $\Ga^\sharp$ is simply the theory of nilpotent endomorphisms in linear algebra, as we explain next; this result underlies the appearance of operators (such as connections or $p$-curvatures or Higgs fields) in our later descriptions of quasi-coherent sheaves on certain stacks.

\begin{proposition}[A spectral decomposition via Cartier duality for $\mathrm{Rep}(\Ga^\sharp$)]
\label{prop:reptheoryGasharp}
The Cartier dual of $\Ga^\sharp$ is $\widehat{\Ga}$ (the formal completion of $\Ga$ at $0$). Moreover, the Cartier duality pairing defines a tautological line bundle on $B\Ga^\sharp \times \widehat{\Ga}$ that induces a (suitably normalized) Fourier-Mukai equivalence
\[ \mathcal{D}(B\Ga^\sharp) \simeq \mathcal{D}(\widehat{\Ga}),\]
carrying the standard symmetric monoidal structure on the left to convolution\footnote{As $\widehat{\Ga}$ is an ind-scheme, pushforward along a $\widehat{\Ga}$-torsor $f:X \to Y$ is not well-behaved, e.g., it does not commute with base change. However, $f^*$ has a left-adjoint whose formation commutes with base change; this left-adjoint is used in the definition of convolution.} on the right. 
\end{proposition}
\begin{proof}[Proof sketch]
The pairing $\Ga^\sharp \times \widehat{\Ga} \to \mathbf{G}_m$ is given by $(x,y) \mapsto \exp(xy) = \sum_{n \geq 0} \gamma_n(x) y^n$. For the rest, we refer to \cite[\S 12.1]{Bhatt2025Apsects}.
\end{proof}

\begin{example}[Cohomology of $\Ga^\sharp$-representations]
The equivalence in Proposition~\ref{prop:reptheoryGasharp} restricts, by passing to compact objects modulo $p^n$ and taking a limit, to an equivalence 
\[ \Perf(B\Ga^\sharp) \simeq \Perf_{A\text{-Perf}}(\widehat{\Ga}),\] 
where the RHS can be viewed as the $\infty$-category of $A[x]$-modules that are $A$-perfect and $(p,x)$-complete. Moreover, by symmetric monoidality, this equivalence can also be checked to be $t$-exact. In particular, the category $\mathrm{Vect}(B\Ga^\sharp)$ of representations of $\Ga^\sharp$ on finite projective $A$-modules is identified with the category $\mathrm{Coh}_{A\text{-Vect}}(\widehat{\Ga})$ of pairs $(M,T)$, where $M$ is a finite projective $A$-module and $T:M \to M$ is a ``quasi-nilpotent'' endomorphism, i.e., $T/p$ is nilpotent. Under this correspondence, the functor of forgetting the $\Ga^\sharp$-action corresponds to $(M,T) \mapsto M$, and the functor $R\Gamma(B\Ga^\sharp,-)$ of global sections corresponds to $(M,T) \mapsto \mathrm{fib}(M \xrightarrow{T} M)$. 
\end{example}

Next, we shall give a Witt vector model for $\Ga^\sharp$; this is essential not only for understanding $\Ga^\sharp$ well (e.g., its cohomology), but also to generalize from de Rham to prismatic cohomology later. Let us introduce some relevant notation first.

\begin{notation}[The Witt ring scheme]
\label{not:WittRing}
Write $W$ for the ring scheme of Witt vectors over $\Spf(A)$; it is corepresented by a polynomial ring $A[x_0,x_1,x_2,...]$ in countably many variables corresponding to the Witt components. There is a natural Frobenius endomorphism $F:W \to W$. Writing $F_* W$ for the $W$-module given by $W$ acting via $F$, we have $W$-linear maps $V:F_*W \to W$ and $F:W \to F_* W$ given by Verschiebung and Frobenius, satisfying $FV=p$. For a $p$-nilpotent ring $R$, the ring $W(R)$ is an inverse limit of locally nilpotent thickenings of $R$, so deformation theory promotes $W$ to an $A$-algebra scheme on $\CAlg^\an_{\Spf(A)}$. Moreover, by the same argument, the induced $A$-algebra structure on $F_* W$ agrees with the one given by $A \xrightarrow{\varphi} A \to W$. We write $F_*$ when emphasizing $W$-module structures and $\varphi_*$ when emphasizing $A$-algebra structures; by the preceding sentence, these describe the same object. 
\end{notation}

Recall that Witt vectors are intimately related to the theory of $\delta$-rings: the functor $W(-)$ can be regarded as a right adjoint to the forgetful functor from $\delta$-rings to rings \cite[Remark 2.7]{BhattScholzePrism}. Moreover, there is also a close connection between divided powers and the theory of $\delta$-rings: if $B$ is a $p$-torsionfree ring with a lift $\varphi:B \to B$ of the Frobenius and $x \in B$, then $x$ admits divided powers in $B$ if and only if $\varphi(x) \in pB$ \cite[Lemma 2.35]{BhattScholzePrism}; see also Lemma~\ref{lem:DeltaDividedPowers}. Combining these leads to the following Witt vector model for $\Ga^\sharp$:

\begin{proposition}[The Witt vector model for $\Ga^\sharp$, {\cite[Variant 3.4.12]{BhattLurieAbsPrismaticCoh}}, {\cite[Lemma 3.2.6]{DrinfeldPrismatization}}]
\label{prop:WittGasharp}
There is a natural identification 
\[ \Ga^\sharp \simeq W[F] := \ker(F:W \to W)\]
of group schemes over $\Spf(A)$. Write $\widetilde{d}:\Ga^\sharp \to W$ for the resulting map.
\end{proposition}

This model gives easy access to flat cohomology\footnote{In general, there are set-theory issues with making sense of flat cohomology: one must choose cutoff cardinals and it is, a priori, unclear if the result depends on the choice of such cardinal. However, for our purposes, Proposition~\ref{prop:WittGasharp} implies these issues are irrelevant.}  with coefficients in $\Ga^\sharp$. Indeed, Proposition~\ref{prop:WittGasharp} gives a short exact sequence
\[ 0 \to \Ga^\sharp \xrightarrow{\widetilde{d}} W \xrightarrow{F} W \to 0\]
of group schemes (see \cite[\S 3.4]{DrinfeldPrismatization} for the faithful flatness of $F$). This leads to:

\begin{proposition}[Cohomology with coefficients in $\Ga^\sharp$, {\cite[Lemma 5.7]{BhattLuriePrismatization}}]
\label{prop:CohGasharp}
For any $R \in \CAlg^\an_{\Spf(A)}$, we have a natural identification 
\[ R\Gamma(\Spec(R),\Ga^\sharp) \simeq \mathrm{fib}(W(R) \to W(R)).\]
In particular, this complex is concentrated in degrees $\leq 1$. 
\end{proposition}

It follows that the functor $R \mapsto R\Gamma(\Spec(R),\Ga^\sharp)$ on $\CAlg^{\an}_{\Spf(A)}$ is the left Kan extension of its restriction to the full (ordinary) subcategory spanned by polynomial algebras over some $p$-adic truncation $A/p^n$ of $A$. In particular, structural statements about this functor can be reduced to the case of classical rings.

\begin{proof}
It suffices to show that $R\Gamma_{fl}(\Spec(R),W) \simeq W(R)$, which follows by contemplating the $V$-adic filtration on $W$ to reduce to the analogous statement for $\Ga$ itself; the latter case follows from flat descent.
\end{proof}

\subsection{The de Rham stack}
\label{ss:GadR}

As we explained earlier, we will build de Rham stacks $p$-adically by replacing the group scheme $\widehat{\Ga}$ appearing in characteristic $0$ de Rham theory (Example~\ref{ex:SimpsondR}) by the group scheme $\Ga^\sharp$. To carry this out, we need a ring stack structure on $\mathrm{Cone}(\Ga^\sharp \to \Ga)$. This is a derived version of the following elementary fact: if $R$ is a torsionfree commutative ring, then the collection of elements $x \in R$ admitting divided powers forms an ideal $R^\sharp \subset R$. While one could construct this ring stack structure explicitly, the Witt vector model for $\Ga^\sharp$ (Proposition~\ref{prop:WittGasharp}) gives an elegant abstract approach, which we discuss next.

\begin{proposition}[The Witt vector model for $\Ga^\dR$, {\cite[Remark 7.9]{BhattLuriePrismatization}}, {\cite[Proposition 3.5.1]{DrinfeldPrismatization}}]
\label{prop:dRRingStackWittVectorModel}
There is a commutative diagram of $W$-module schemes on $\CAlg^\an_{\Spf(A)}$
\[ \begin{tikzcd}[column sep=6em, row sep=2.2em,
  execute at end picture={
    \foreach \c in {1,2,3}{
      \node[inner sep=2pt, fit=(\tikzcdmatrixname-1-\c)(\tikzcdmatrixname-2-\c)] (col\c) {};
      \draw[dotted] (col\c.north west) to[bend right=15] (col\c.south west);
      \draw[dotted] (col\c.north east) to[bend left=15] (col\c.south east);
    }}]
\Ga^\sharp \arrow[d, "d"] & \Ga^\sharp \oplus F_*W \arrow[l, "\text{pr}_1"'] \arrow[r, "\text{pr}_2"] \arrow[d, "{(\widetilde{d},V)}"] & F_*W \arrow[d, "p"] \\
\Ga & W \arrow[l, "\gamma"] \arrow[r, "F"] & F_*W
\end{tikzcd} \]
where $\gamma$ is projection to the $0$-th Witt component. Moreover,  the maps relating the middle column to the first and last columns are quasi-isomorphisms, where each column is regarded as a 2-term complex, as indicated by the dotted parentheses.
\end{proposition}
\begin{proof}
The existence of the commutative diagram is readily verified using the Witt vector identities.  As both $\gamma$ and $F$ are surjective, the quasi-isomorphy statement reduces to checking that $(\widetilde{d},V)$ induces identifications
\[ 0 \oplus F_* W \simeq \ker(\gamma) \quad \text{and} \quad \Ga^\sharp \oplus 0 \simeq \ker(F);\]
the first of these is standard, and the second is Proposition~\ref{prop:WittGasharp}.
\end{proof}

In particular, for any $R \in \CAlg^\an_{\Spf(A)}$, we obtain a natural identification of $A$-modules
\begin{equation}
\label{eq:RGammaConeWitt}
\tag{$\Ga^\dR\text{-Witt}$}
R\Gamma(\Spec(R),\mathrm{Cone}(\Ga^\sharp \to \Ga)) \simeq \mathrm{Cone}(\varphi_*W(R) \xrightarrow{p} \varphi_*W(R)),
\end{equation}
where, as in Notation~\ref{not:WittRing}, we write $\varphi_*$ to emphasize the $A$-module structure. To proceed via transmutation, we now simply need the following standard construction of animated rings, generalizing the Koszul complex:

\begin{construction}[Derived quotients, {\cite[Construction 2.1]{BhattLuriePrismatization}}, {\cite[Proposition 3.2.6]{KhanRydhVirtualCartier}}]
\label{cons:DerivedQuotients}
Given an animated ring $B$ and an element $f$ (incarnated as a map $\mathbf{Z}[x] \xrightarrow{x \mapsto f} B$), the cone $B/f := \mathrm{Cone}(B \xrightarrow{f} B)$ has a natural animated $B$-algebra structure via its realization as 
\[ B/f \simeq B \otimes_{\mathbf{Z}[x]} \mathbf{Z},\] where $x$ acts by $0$ on $\mathbf{Z}$ on the right. More generally, given a {\em generalized Cartier divisor} $d:I \to B$ (i.e., an invertible $B$-module $I$ and a $B$-module map $d:I \to B$), there is a natural animated $B$-algebra structure on $B/I = \mathrm{Cone}(d)$: classify $d$ via a map $\mathrm{Spec}(B) \to \mathbf{A}^1/\mathbf{G}_m$ and declare $\Spec(B/I)$ to be the pullback of $B\mathbf{G}_m \hookrightarrow \mathbf{A}^1/\mathbf{G}_m$. Taking powers, one defines $B/I^n$ as an animated $B$-algebra for any $n \geq 1$. 
\end{construction}

Using the above language, and in view of the formula \eqref{eq:RGammaConeWitt}, it is sensible to define:

\begin{definition}[The ring stack $\Ga^\dR$] 
\label{def:GadR}
The {\em de Rham ring stack} $\Ga^\dR$ is the functor
\[ \Ga^\dR:\CAlg^\an_{\Spf(A)} \to \CAlg^\an_{\Spf(A)}, \quad \quad R \mapsto \varphi_* W(R)/p. \]
\end{definition}

We shall informally write $\Ga^\dR = \Ga/\Ga^\sharp$. As the notation suggests, there is a natural map 
\[ \eta:\Ga \to \Ga^\dR\] of animated $A$-algebra stacks. For instance, in the Witt vector model, this corresponds to the map $\Ga = W/VF_*(W) \xrightarrow{F} F_* W/p$ induced by the Witt vector Frobenius. The map $\eta$ is surjective locally in the quasi-syntomic topology, and the kernel of the induced map on $\pi_0(-)$ is a locally nilpotent ideal.

\begin{remark}[The derivedness of $\Ga^\dR$]
\label{rmk:CharpvsChar0deRhamRingStack}
Unlike its characteristic $0$ counterpart $\Ga/\widehat{\Ga}$ (Example~\ref{ex:SimpsondR}), the ring stack $\Ga^\dR$ defined above is typically not set-valued on ordinary rings: the map $\Ga^\sharp \to \Ga$ (equivalently, $p:W(R) \to W(R)$) fails to be injective as soon as $R$ has nilpotents. This has two related consequences. First, even evaluating the transmutation $X^\dR$ of a smooth scheme $X/\Spf(A)$ on an ordinary $p$-nilpotent $A$-algebra $R$ involves the mildly derived ring $W(R)/p$, so derived algebraic geometry is unavoidable.\footnote{For most purposes, however, $1$-truncated animated rings (i.e., those with $\pi_i = 0$ for $i \geq 2$) suffice. These admit an explicit model: commutative differential graded algebras $d:A_1 \to A_0$, with the $A_0$-module $A_1$ regarded as an ``ideal'' of $A_0$. Drinfeld calls these ``quasi-ideals'' and studies them systematically in \cite{DrinfeldRingGroupoid}; this model makes computations with $1$-truncated animated rings straightforward, and will be used implicitly in the sequel.} Relatedly, $X^\dR$ is genuinely groupoid-valued, even on ordinary rings. In fact, passing to isomorphism classes has a drastic effect: for $X$ smooth over a perfect field $k$ of characteristic $p$, the sheafification of the corresponding set valued functor $R \mapsto \pi_0(X^\dR(R))$ is the Frobenius twist $X^{(1)}$ (Propositions~\ref{prop:SqZerodRring} and \ref{prop:deRhamGerbe}), and the resulting gerbe $X^\dR \to X^{(1)}$ has pullback and pushforward encoding Cartier descent and the Cartier isomorphism respectively (Corollaries~\ref{cor:CartierDescent} and \ref{cor:Cartier}).
\end{remark}

Using $\Ga^\dR$ and transmutation (Construction~\ref{cons:Transmutation}), we obtain:

\begin{definition}[The de Rham stack via transmutation]
\begin{enumerate}
\item For any $X \in \Fun(\CAlg^\an_{\Spf(A)},\Sp)$, define its {\em de Rham stack} $(X/A)^\dR \in \Fun(\CAlg^\an_{\Spf(A)},\Sp)$ via transmutation of $\Ga^\dR$, i.e., as the composition
\[ \CAlg^\an_{\Spf(A)} \xrightarrow{\Ga^\dR} \CAlg^\an_{\Spf(A)} \xrightarrow{X} \Sp,\]
i.e., for $R \in \CAlg^\an_{\Spf(A)}$, we have
\[ (X/A)^\dR(R) := X(\varphi_* W(R)/p) := \mathrm{Map}_{\Spf(A)}(\Spec(\varphi_* W(R)/p), X).\]

\item The map $\eta:\Ga \to \Ga^\dR$ induces, via transmutation, a map $\eta_X:X \to (X/A)^\dR$ called the {\em de Rham structure map}.  

\item For a map $f:X \to Y$ in $\Fun(\CAlg^\an_{\Spf(A)},\Sp)$, we call $(X/Y)^\dR := (X/A)^\dR \times_{(Y/A)^\dR} Y$ its {\em relative de Rham stack}, where the fibre product is formed using $\eta_Y:Y \to (Y/A)^\dR$, i.e., $(X/Y)^\dR$ is defined by the fibre square
\[ \xymatrix{
(X/Y)^\dR \ar[r] \ar[d]_-{f^\dR} \ar@{}[dr]|(.25){\lrcorner} & (X/A)^\dR \ar[d] \\
Y \ar[r]^-{\eta_Y} & (Y/A)^\dR.
} \]
The map $\eta_X$ refines to a map $\eta_{X/Y}:X \to (X/Y)^\dR$ by the universal property of the fibre product. 
\end{enumerate}
\end{definition}

The most important examples are:

\begin{example}
$(\Ga/A)^\dR = \Ga^\dR$ by definition. As the fibre of $\eta:\Ga \to \Ga^\dR$ is locally nilpotent, one computes that $\mathbf{G}_m^\dR = \mathbf{G}_m/\mathbf{G}_m^\sharp := \mathrm{Cone}(\mathbf{G}_m^\sharp \to \mathbf{G}_m)$, where $\mathbf{G}_m^\sharp$ is the PD-hull of the identity point in $\mathbf{G}_m$.
\end{example}

The following records some basic stability properties of this construction that help work with it:

\begin{proposition}[Properties of the de Rham stack construction]
\label{prop:dRStackBasicProp}
\begin{enumerate}

    \item If $X$ is a $k$-scheme, then $(X/A)^\dR$ is an fpqc sheaf.

    \item  The endofunctor $X \mapsto (X/A)^\dR$ of $\Fun(\CAlg^\an_{\Spf(A)},\Sp)$ commutes with limits and preserves representable quasi-syntomic covers. 

\item If $f:U \to X$ is representable and \'etale, then $(U/A)^\dR \to (X/A)^\dR$ is also representable and \'etale, and satisfies $U \simeq (U/A)^\dR \times_{(X/A)^\dR} X$ via the natural maps.

\item Fix an extension $K/k$ of perfect fields; let $B=W(K)$, given its natural $A$-algebra structure. For any $X \in \Fun(\CAlg^\an_{\Spf(A)},\Sp)$ with base change $X_B := X \times_{\Spf(A)} \Spf(B)$, regarded as an object of $\Fun(\CAlg^\an_{\Spf(B)},\Sp)$, there is a natural base change identification $(X_B/B)^\dR \simeq (X/A)^\dR \times_{\Spf(A)} \Spf(B)$.
\end{enumerate}
\end{proposition}

It follows that $(X/A)^\dR$ has good geometric properties, e.g., if $X/k$ is a smooth affine $k$-scheme, then $(X/A)^\dR$ is the quotient (as quasi-syntomic sheaves) of an affine smooth scheme over $\Spf(A)$ by an affine equivalence relation. The name ``de Rham stack'' is justified by the following:

\begin{proposition}[Comparison with de Rham cohomology]
\label{prop:dRcompStack}
Let $f:X \to Y$ be a smooth representable morphism, and write $f^\dR:(X/Y)^\dR \to Y$ for the relative de Rham stack. Then $\eta_{X/Y}:X \to (X/Y)^\dR$ is a quasi-syntomic cover. Moreover, there is a natural identification $Rf^\dR_* \calO_{(X/Y)^\dR} \simeq Rf_* \Omega^\bullet_{X/Y}$. 
\end{proposition}
\begin{proof}[Proof sketch]
We shall freely use the Rees module dictionary realizing filtrations on an object of $\mathcal{D}(Y)$ in terms of a lift along $Y = Y \times \mathbf{G}_m/\mathbf{G}_m \subset Y \times \mathbf{A}^1/\mathbf{G}_m$ (see \cite[\S 2.2]{BhattFGaugeNotes} and references therein). By scaling the differential $d:\Ga^\sharp \to \Ga$ by parameter $\lambda$ on $\mathbf{A}^1/\mathbf{G}_m$, one defines a new $\Ga$-algebra stack $\Ga^{\dR,+}$ over $\mathbf{A}^1/\mathbf{G}_m$ that degenerates from $\Ga^\dR$ over $\mathbf{G}_m/\mathbf{G}_m \subset \mathbf{A}^1/\mathbf{G}_m$ to $\Ga \oplus B\Ga^\sharp\{1\}$ at $B\mathbf{G}_m \subset \mathbf{A}^1/\mathbf{G}_m$, where $\{1\}$ records a twist coming from conormal bundle of $B\mathbf{G}_m \subset \mathbf{A}^1/\mathbf{G}_m$. By transmutation, this gives a stack\footnote{This stack may be reasonably called the ``Hodge filtered de Rham stack'': as indicated above, its cohomology encodes the Hodge filtration on de Rham cohomology.} $(X/Y)^{\dR,+} \to Y \times \mathbf{A}^1/\mathbf{G}_m$, realizing the relative de Rham stack $(X/Y)^\dR$ over $Y \times \mathbf{G}_m/\mathbf{G}_m$, and the relative Hodge stack $BT_{X/Y}^\sharp\{1\}/\mathbf{G}_m$ over $Y \times B\mathbf{G}_m$. Pushing forward the structure sheaf and using Proposition~\ref{prop:reptheoryGasharp}, one learns that $Rf^\dR_* \calO_{(X/Y)^\dR} \in \mathcal{D}(Y)$ admits a natural filtration with associated graded given by the graded object $\oplus_i Rf_* \Omega^i_{X/Y}[-i]$. Reducing to the case where $f$ is affine, it suffices to show that the extension classes associated to this filtration are given by the de Rham differential; this is calculated explicitly by reduction to the case $X=\mathbf{A}^1 \times Y$. 
\end{proof}

Going one step further, one can describe the entire derived category $\mathcal{D}( (X/Y)^\dR)$ in  terms of the sheaf $\mathcal{D}_{X/Y} \in \mathcal{D}(X)$ of relative crystalline differential operators. Indeed, the covering property in Proposition~\ref{prop:dRcompStack} and faithfully flat descent, combined with analysis of the descent datum using Proposition~\ref{prop:reptheoryGasharp}, show (see \cite[\S 3.4]{TerentiukOV}):

\begin{proposition}[Relation to $\mathcal{D}$-modules]
\label{prop:DmoduledRstackQCoh}
For a smooth representable map $f:X \to Y$,  pullback along $\eta_{X/Y}:X \to (X/Y)^\dR$ lifts to a fully faithful functor $\mathcal{D}( (X/Y)^\dR) \to \mathcal{D}(X,\mathcal{D}_{X/Y})$, with essential image exactly those objects whose mod $p$ reduction has locally nilpotent $p$-curvature\footnote{This condition on a $\mathcal{D}_{X/Y}$-module can be checked after base change on $Y$ and localization on $X$. Moreover, if $X$ admits an \'etale map $X \to \mathbf{A}^n_Y$, then the condition simply means that the vector fields $\{\frac{\partial}{\partial x_i}\}_{i=1,...,n}$, regarded as sections of $\mathcal{D}_{X/Y}$, act locally nilpotently modulo $p$. We refer to Remark~\ref{rmk:pCurvature} for a stack-theoretic perspective on  $p$-curvature.}.
\end{proposition}

\begin{remark}[Incorporating the $V$-adic topology]
\label{rmk:VadicTopdR}
Definition~\ref{def:GadR} of $\Ga^\dR$ is not quite ``correct'' from a principled perspective: we ignored the $V$-adic topology on $W(R)$ in our discussion. A better definition would take this topology into account and replace $\Spec(W(R)/p)$ with $\colim_n \Spec(W_n(R)/p)$ in the definition of de Rham stacks. Fortunately, the more pedestrian approach taken above also produces the correct result when transmuting schemes (our main purpose in these notes): if $X$ is a $k$-scheme and $R$ is a $p$-nilpotent $A$-algebra, then the natural map $\mathrm{Map}_k(\mathrm{Spec}(W(R)/p),X) \to \mathrm{Map}_k(\colim_n \mathrm{Spec}(W_n(R)/p),X)$ is an equivalence by Tannaka duality arguments (see \cite[Remark 4.1.5]{Bhatt2025Apsects} for this argument in a variant context).
\end{remark}

\begin{remark}[Sheared de Rham stacks]
\label{rmk:ShearedVariant}
The ``nilpotent $p$-curvature'' constraint in Proposition~\ref{prop:DmoduledRstackQCoh} raises an obvious question: is there a version of the stack-theoretic perspective on $\mathcal{D}$-modules that accommodates {\em all} $\mathcal{D}$-modules? The answer is ``yes'' and the subject of ongoing joint work with Kanaev, Mathew, Vologodsky and Zhang on ``sheared prismatization''. The basic idea is simple: replace the group scheme $\mathbf{G}_a^\sharp$ with the group ind-scheme $\widehat{\mathbf{G}_a^\sharp}$ given by the completed PD-hull of $0 \in \Ga$, since one has $\mathcal{D}(B\widehat{\Ga^\sharp}) \simeq \mathcal{D}(\Ga)$, extending Proposition~\ref{prop:reptheoryGasharp}. Slightly more precisely, one defines a square-zero extension ${}^sW \to W$ of the Witt ring scheme with the property that the map ${}^s W[F] \to W[F]$ identifies with $\widehat{\Ga^\sharp} \to \Ga^\sharp$; the ring stack ${}^s W/p$ then gives a theory of sheared de Rham stacks via transmutation. We refer to \cite{MathewVologodskyVideo, DrinfeldshearedW,BarzLogarithmic,BMVShearedWitt} for an introduction to these ideas.
\end{remark}

\subsection{de Rham stacks in characteristic $p$}

In this subsection, we work in characteristic $p$, i.e., over the perfect field $k$ itself. A fundamental feature of de Rham cohomology in characteristic $p$ is that the de Rham differential is linear over $p$-th powers, i.e., $d(f^p\omega) = f^p d(\omega)$ for a function $f$ and a form $\omega$. This property ensures that the de Rham complex $\Omega^\bullet_{X/k}$ of a smooth $k$-scheme $X$ is naturally quasi-coherent complex on the Frobenius twist $X^{(1)}$. Moreover, Cartier identified the cohomology sheaves of this complex: there is a natural isomorphism
\begin{equation}
\label{eq:CartierIsom}
\tag{Cart}
\mathcal{H}^*(\Omega^\bullet_{X/k}) = \Omega^*_{X^{(1)}/k}
\end{equation}
of graded coherent $\calO_{X^{(1)}}$-algebras. Our goal in this section is to study (and reprove) these statements via geometrization. The critical observation is the following:

\begin{proposition}[Birth of square-zero extensions]
\label{prop:SqZerodRring}
The $k$-algebra $\Ga^\dR$ is naturally and uniquely a square-zero extension of $F_* \Ga$ by $BF_*\Ga^\sharp$ (where $F_*$ tracks the $\Ga$-module, and thus $k$-module, structures). 
\end{proposition}
\begin{proof}[Proof sketch]
The functor $\Ga^\dR = \varphi_* W/p$ is left Kan extended from smooth $k$-algebras. Moreover, by flat descent, its restriction to smooth $k$-algebras is right Kan extended from semiperfect $k$-algebras (i.e., $k$-algebras with a surjective Frobenius). Thus, to show existence and uniqueness, it suffices to show the following: for a semiperfect $k$-algebra $R$, the animated $k$-algebra $\Ga^\dR(R)$ is uniquely a square-zero extension of $\varphi_* R$ by $\varphi_* \Ga^\sharp(R)[1]$. As $R$ has characteristic $p$, we have $FV=VF=p$ on $W(R)$. Thus, we can write $p:W(R) \to W(R)$ as the composition 
\[ W(R) \xrightarrow{F} F_* W(R) \xrightarrow{V} W(R).\]
The first map is surjective (by semiperfectness of $R$) with kernel $\Ga^\sharp(R)$ (Proposition~\ref{prop:WittGasharp}), while the second map is injective with cokernel $R$ via the $0$-th Witt component. It follows that $\pi_0(W(R)/p) = R$, $\pi_1(W(R)/p) =\Ga^\sharp(R)$, and the rest vanish. Applying $\varphi_*$, the claim follows from the general statement that any $1$-truncated animated $R$-algebra is uniquely the square-zero extension of its $\pi_0$ by $\pi_1$ (this follows from \cite[Corollary 25.3.6.3]{LurieSAG}).
\end{proof}

For any morphism $f:X \to Y$, write $X^{(1)/Y}$ for its Frobenius twist relative to $Y$ (and drop the $Y$ from the superscript when $Y=\mathrm{Spec}(k)$), defined by the fibre square
\[ \xymatrix{
X^{(1)/Y} \ar[r] \ar[d] \ar@{}[dr]|(.25){\lrcorner} & X \ar[d]^-f \\
Y \ar[r]^-{F_Y} & Y,
} \]
where $F_Y$ denotes the absolute Frobenius of $Y$. The features of de Rham cohomology mentioned above are mostly subsumed by:

\begin{proposition}[The gerbe structure on de Rham stack]
\label{prop:deRhamGerbe}
Given a smooth morphism $X \to Y$, there is a natural $BT_{X^{(1)/Y}/Y}^\sharp$-torsor $\pi_{X/Y}:(X/Y)^\dR \to X^{(1)/Y}$ over $Y$. Moreover, $R\pi_{X/Y,*} \calO_{(X/Y)^\dR}$ is naturally identified with the de Rham complex $\Omega^\bullet_{X/Y} \in \mathcal{D}(X^{(1)/Y})$.
\end{proposition}

In the sequel, we shall simply write $\pi$ for $\pi_{X/Y}$ unless there is risk of confusion.

\begin{proof}[Proof sketch]
The map $\pi_{X/Y}$ is obtained by transmutation from the map $\Ga^\dR \to F_* \Ga$ in Proposition~\ref{prop:SqZerodRring}; the torsor  structure on this map is obtained by deformation theory from the second half of the same proposition. The rest follows by running the proof of Proposition~\ref{prop:dRcompStack} directly over $X^{(1)/Y}$. 
\end{proof}

We shall informally refer to $\pi$ as {\em the de Rham gerbe}. The existence of the gerbe structure on $\pi_{X/Y}$ imposes strong constraints on the geometry of $(X/Y)^\dR$, as we shall see soon. But first, we discuss when the gerbe splits:

\begin{example}[Splitting the de Rham gerbe via Frobenius lifts]
\label{ex:Fliftsplitgerbe}
Say $X/k$ is a smooth $k$-scheme equipped with a lift $\widetilde{X}/A$ together with a $\delta$-structure. Then the gerbe $(X/k)^\dR \to X^{(1)}$ splits\footnote{In fact, with some further work, one can show the following:  the de Rham gerbe $\pi:(X/k)^\dR \to X^{(1)}$ splits if and only if $X$ admits a lift to $W_2(k)$ with a lift of Frobenius. In particular, if $X$ is a curve of genus $\geq 2$ or a supersingular elliptic curve, then the de Rham gerbe does {\em not} split.}: indeed, by transmutation, this follows from the fact that for any $k$-algebra $R$ and map $\mu:\mathrm{Spec}(R) \to X^{(1)}$, the induced map $\mathrm{Spec}(R) \to \widetilde{X}^{(1)} := \widetilde{X} \times_{\mathrm{Spf}(A),\varphi} \mathrm{Spf}(A)$ lifts uniquely to a $\delta$-map $\mathrm{Spf}(W(R)) \to \widetilde{X}^{(1)}$ by the description of $W(-)$ as an adjoint to the forgetful functor from $\delta$-rings; taking special fibres then gives a map $\mathrm{Spec}(W(R)/p) \to X^{(1)}$, and thus a point of $(X/k)^\dR(R)$ lifting $\mu$. 
\end{example}

The gerbe structure immediately leads to:

\begin{corollary}[Cartier isomorphism]
\label{cor:Cartier}
Fix a smooth map $f:X \to Y$. There is a natural identification
\[ \bigoplus_i R^i \pi_{*} \calO_{(X/Y)^\dR} \simeq \bigoplus_i \Omega^i_{X^{(1)/Y}/Y} \]
of graded coherent $\calO_{X^{(1)/Y}}$-algebras. 
\end{corollary}

In view of  Proposition~\ref{prop:deRhamGerbe}, this is  the Cartier isomorphism \eqref{eq:CartierIsom}.

\begin{proof}[Proof sketch]
The formal dg-algebra $\bigoplus_i \Omega^i_{X^{(1)/Y}/Y}[-i]$ is the pushforward of the structure sheaf along the trivial torsor $BT_{X^{(1)/Y}/Y}^\sharp \to X^{(1)/Y}$. The corollary now follows from the following general fact about gerbes: if $G$ is an affine commutative group scheme over a base $Z$ and $g:W \to Z$, $h:W' \to Z$ are two $BG$-torsors, then there is a canonical identification  $R^i g_* \calO_W \simeq R^i h_* \calO_{W'}$ of quasi-coherent sheaves, see \cite[Proof of Corollary 2.7.2]{BhattFGaugeNotes}. (We warn the reader that this is an ``abelian'' statement; the complexes $Rg_* \calO_W$ and $Rh_* \calO_{W'}$ need not be isomorphic.)
\end{proof}

Another, related, result of Cartier is also an immediate consequence of the gerbe structure:

\begin{corollary}[Cartier descent]
\label{cor:CartierDescent}
Fix a smooth map $f:X \to Y$. Then pullback along the de Rham gerbe $\pi:(X/Y)^\dR \to X^{(1)/Y}$ gives a fully faithful functor
\[ \pi^*:\mathrm{QCoh}(X^{(1)/Y}) \to \mathrm{QCoh}((X/Y)^\dR) \]
of abelian categories of quasi-coherent sheaves.
\end{corollary}

In terms of the $\mathcal{D}$-module description of $(X/Y)^\dR$ (Proposition~\ref{prop:DmoduledRstackQCoh}) and the notion of $p$-curvature (Remark~\ref{rmk:pCurvature}), the proof sketched below also shows that the essential image of $\pi^*$ consists of quasi-coherent $\mathcal{D}_{X/Y}$-modules with vanishing $p$-curvature.

\begin{proof}[Proof sketch]
As $\pi$ is a $BT_{X^{(1)/Y}/Y}^\sharp$-torsor (Proposition~\ref{prop:deRhamGerbe}) and the claim is  local on $X^{(1)/Y}$, we reduce to the following standard statement: for a flat affine group scheme $G$ over a ring $R$, pullback along $BG \to \Spec(R)$ is fully faithful on quasi-coherent sheaves.
\end{proof}

Next, observe that the Frobenius kills the group scheme $\Ga^\sharp$ in characteristic $p$: if an element $x$ of an $\mathbf{F}_p$-algebra has divided powers $\{\gamma_n(x)\}_{n \geq 1}$, then $p$-th power of each $\gamma_n(x)$ must be $0$. Hence, the Frobenius on $\Ga$ factors as 
\[ \Ga \xrightarrow{\eta} \Ga^\dR \xrightarrow{\pi} F_* \Ga. \]
By transmutation, we learn that for any $k$-scheme $X$, the relative Frobenius $F_{X/k}:X \to X^{(1)}$ factors as 
\[ X \xrightarrow{\eta_X} (X/k)^\dR \xrightarrow{\pi} X^{(1)}.\] 
In particular, the de Rham gerbe $\pi:(X/k)^\dR \to X^{(1)}$ splits after pullback along the relative Frobenius $F_{X/k}$. This leads to the following result of Petrov:

\begin{corollary}[$F$-split varieties have a split de Rham complex, \cite{PetrovQuasiFsplit}]
\label{cor:PetrovFSplit}
Let $X/k$ be a smooth $k$-scheme. Assume $X$ is $F$-split, i.e., $\calO_X \to F_* \calO_X$ is split in $\mathrm{Vect}(X)$. Then the de Rham complex splits, i.e., there is a natural (in the splitting) isomorphism
\[ \Omega^\bullet_{X/k} \simeq \bigoplus_i \Omega^i_{X^{(1)}/k}[-i] \]
in $\mathcal{D}(X^{(1)})$. 
\end{corollary}
\begin{proof}
Indeed, we have such an isomorphism after applying $F_{X/k}^*$ to both sides by the proof of Corollary~\ref{cor:Cartier}: the de Rham gerbe splits after pullback along $F_{X/k}$. Pushing forward along $F_{X/k}$ and applying the $F$-splitting (assumed to exist) then gives the desired statement.
\end{proof}

An earlier, and important, decomposition result for the de Rham complex was proven by Deligne--Illusie \cite{DeligneIllusie}. Their original proof required a careful bookkeeping of homotopies. Using the stack-theoretic point of view, Drinfeld recently explained a ``pure thought'' proof of a generalization (see \cite[\S 5.12.1]{DrinfeldPrismatization} and \cite[Corollary 2.7.14]{BhattFGaugeNotes}):

\begin{corollary}[Deligne--Illusie, \cite{DeligneIllusie}]
\label{cor:DeligneIllusie}
Say $X/k$ is a smooth $k$-scheme equipped with a lift $\widetilde{X}/W_2(k)$. If $\dim(X) < p$, then the de Rham complex splits, i.e., there is a natural isomorphism 
\[ \Omega^\bullet_{X/k} \simeq \bigoplus_i \Omega^i_{X^{(1)}/k}[-i] \]
in $\mathcal{D}(X^{(1)})$. 
\end{corollary}

An elaboration of the ideas sketched below leads to a natural operator $\Theta$, called the {\em Sen operator}, acting on $\Omega^\bullet_{X/k}$ with generalized eigenvalue $-i$ on $\mathcal{H}^i(\Omega^\bullet_{X/k}) \simeq \Omega^i_{X^{(1)}/k}$. Petrov used this structure in his proof \cite{PetrovdeRham} that the dimension assumption above cannot be weakened, answering a longstanding open problem.

\begin{proof}[Proof sketch]
We shall prove that the $k$-algebra square-zero extension 
\[ BF_* \Ga^\sharp \to \Ga^\dR \to F_* \Ga\]
from Proposition~\ref{prop:SqZerodRring}, when viewed as a $W_2(k)$-algebra square-zero extension, admits a $\mu_p$-action, where $\mu_p$ acts trivially on the ring $F_*\Ga$ and by weight $1$ on the ``ideal'' $BF_*\Ga^\sharp$. By transmutation, for any $X/k$ and $\widetilde{X}$ as in the corollary, there is an induced $\mu_p$-action on $(X/k)^\dR$ compatible with the trivial action on $X^{(1)}$ under $\pi:(X/k)^\dR \to X^{(1)}$ and the weight $1$ scaling action on the structure group $BT_{X^{(1)}/k}^\sharp$ of the torsor $\pi$. The resulting $\mu_p$-action on $R\pi_* \calO_{(X/k)^\dR} \simeq \Omega^\bullet_{X/k}$ then has weight $-i$ on $\mathcal{H}^i(\Omega^\bullet_{X/k}) \simeq \Omega^i_{X^{(1)}/k}$, which gives the desired splitting as the integers $\{0,1...,\dim(X)\}$ are pairwise distinct mod $p$ if $\dim(X)<p$. 

To construct the preceding action, we may assume $k=\mathbf{F}_p$ by universality, so $\varphi:A \to A$ is the identity. In this case, we shall use the Witt vector model
\[ \Ga^\dR = \mathrm{Cone}( W \xrightarrow{p} W).\]
Let $\mu_p \subset \mathbf{G}_m \stackrel{[\cdot]}{\hookrightarrow} W^*$ act  on the source $W$ by scaling, and trivially on the target $W$. By the Witt vector identities $p=FV=VF$, this induces a $\mu_p$-action on $\Ga^\dR$ as an animated $\mathbf{Z}$-algebra stack with the desired behaviour on $\pi_0(-)$ and $\pi_1(-)$; it is a pleasant exercise to check that this action is $\mathbf{Z}/p^2$-linear but not $\mathbf{Z}/p$-linear.
\end{proof}

\begin{remark}[The $p$-curvature and its horizontality]
\label{rmk:pCurvature}
Let $f:X \to Y$ be a smooth map over $k$. The $BT_{X^{(1)/Y}/Y}^\sharp$-action on $(X/Y)^\dR$ translates, via a vector bundle twisted version of Proposition~\ref{prop:reptheoryGasharp}, to the following structure on quasi-coherent sheaves: for any $E \in \mathcal{D}((X/Y)^\dR)$, there is a natural action of $\pi^* \mathrm{Sym}^*_{X^{(1)/Y}}(T_{X^{(1)/Y}/Y})$ on $E$ such that each local section of $T_{X^{(1)/Y}/Y}$ acts locally nilpotently. (Indeed, such an action is encoded, via Proposition~\ref{prop:reptheoryGasharp}, in the object 
\[ \mathrm{act}^*E \in \mathcal{D}( BT_{X^{(1)/Y}/Y}^\sharp \times_{X^{(1)/Y}} (X/Y)^\dR )\] 
obtained by pulling $E$ back along the action map $\mathrm{act}:BT_{X^{(1)/Y}/Y}^\sharp \times_{X^{(1)/Y}} (X/Y)^\dR \to (X/Y)^\dR$.) Under the equivalence from Proposition~\ref{prop:DmoduledRstackQCoh} with certain $\mathcal{D}_{X/Y}$-modules, this is simply the $p$-curvature action. In this dictionary, the fact that the preceding action takes place on $(X/Y)^\dR$ itself (rather than merely after pullback to $X$) encodes the horizontality of the $p$-curvature, i.e., that it is a map of $\mathcal{D}$-modules.
\end{remark}

\subsection{The crystalline stack}
\label{ss:CrysStack}

Grothendieck's crystalline cohomology was invented \cite{GrothendieckCrystals} to explain the following miracle: for a smooth variety $X/k$ with a lift $\widetilde{X}/\Spf(A)$, the de Rham cohomology $R\Gamma_\dR(\widetilde{X}/\Spf(A))$ of $\widetilde{X}/\Spf(A)$ is independent of the choice of lift $\widetilde{X}$ and, in fact, is  functorial in $X$. In particular, the Frobenius endomorphism  of $X$ acts on $R\Gamma_\dR(\widetilde{X}/\Spf(A))$, even though Frobenius typically does not lift to $\widetilde{X}$. From the stack-theoretic point of view, this miracle has a simple explanation, ultimately relying on the following:

\begin{proposition}
\label{prop:dRstackFpalgbera}
The animated $A$-algebra structure on the stack $\Ga^\dR$ on $\CAlg^\an_{\Spf(A)}$ refines to an animated $k$-algebra structure. 
\end{proposition}
\begin{proof}
As $\mathbf{F}_p$ is obtained from $\mathbf{Z}$ by freely setting $p=0$ in animated rings, the formula $\Ga^\dR = \varphi_* W/p$ shows that $\mathbf{Z}$-algebra structure on $\Ga^\dR$ refines to an $\mathbf{F}_p$-algebra structure, whence the $A$-algebra structure refines to one over $k=A \otimes_{\mathbf{Z}} \mathbf{F}_p$.
\end{proof}

Thanks to this, it is sensible to define a ``crystallization functor'' via transmutation of $\Ga^\dR$:

\begin{definition}
\label{def:CrysStack}
For $X \in \Fun(\CAlg^\an_k,\Sp)$, define its {\em crystallization} $(X/A)^\crys \in \Fun(\CAlg^\an_{\Spf(A)},\Sp)$ as  the composition
\[ \CAlg^\an_{\Spf(A)} \xrightarrow{\Ga^\dR} \CAlg^\an_k \xrightarrow{X} \Sp,\]
where the first arrow results from Proposition~\ref{prop:dRstackFpalgbera}. Thus, for $R \in \CAlg^\an_{\Spf(A)}$, we have
\[ (X/A)^\crys(R) := X(\varphi_* W(R)/p) := \mathrm{Map}_{\Spec(k)}(\Spec(\varphi_* W(R)/p), X). \]
\end{definition}

Unwinding definitions yields the following:

\begin{corollary}[The crystalline miracle]
\label{cor:CrysMiracle}
For any $X \in \Fun(\CAlg^\an_{\Spf(A)},\Sp)$ with special fibre (aka restriction) $X_k \in \Fun(\CAlg^\an_k,\Sp)$, we have a natural isomorphism $(X/A)^\dR \simeq (X_k/A)^\crys$. In particular, $(X/A)^\dR$ (and hence the de Rham cohomology of $X/A$) depends functorially on the special fibre $X_k$. 
\end{corollary}

\begin{remark}[Charts for the crystalline stack]
\label{rmk:ChartCrystallization}
Corollary~\ref{cor:CrysMiracle} and basic properties of the de Rham stack imply the following: given a smooth $k$-scheme $X$, any lift $\widetilde{X}/\Spf(A)$ of $X/k$ provides a quasi-syntomic cover $\widetilde{X} \to (X/A)^\crys$ of the crystalline stack of $X$.
\end{remark}

Going one step further, one can even identify crystals:

\begin{proposition}
For a smooth (or just quasi-syntomic) $k$-scheme $X$, there is a natural identification between $\mathcal{D}( (X/A)^\crys)$ and the $\infty$-category of crystals on Grothendieck's crystalline site of $X$ relative to $\Spf(A)$.
\end{proposition}
\begin{proof}[Proof sketch]
Let us construct the functor in one direction. The (big) crystalline site $(X/A)_\crys$ (see \cite[\href{https://stacks.math.columbia.edu/tag/07I5}{Tag 07I5}]{StacksProject}) consists of quadruples $(R,J,\gamma,\alpha)$, where $R$ is a $p$-nilpotent $A$-algebra, $J \subset R$ is an ideal, $\gamma$ is a divided power structure on $J$ compatible with the one on $(p)$, and $\alpha:\mathrm{Spec}(R/J) \to X$ is a $k$-map; the functor $(R,J,\gamma,\alpha) \mapsto R$ is the crystalline structure presheaf $\calO_\crys$.  Given $(R,J,\gamma,\alpha) \in (X/A)_\crys$, note that there is a natural animated $R$-algebra map
\[ R/J \to \Ga^\dR(R) \]
provided by the divided power structure $\gamma$ on $J$. Thus, $\alpha:\Spec(R/J) \to X$ determines a point of $(X/A)^\crys(R)$. Pulling back modules to these points for varying points yields a map
\[ \mathcal{D}( (X/A)^\crys) \to \lim_{(R,J,\gamma,\alpha) \in (X/A)_\crys} \mathcal{D}(R) \ \  =: \ \ \mathcal{D}_{\crys}( (X/A)_\crys,\calO_\crys),\]
which one then checks to be an equivalence. 
\end{proof}

\subsection{Non-abelian Hodge theory in characteristic $p$}
\label{ss:NAHTcharp}
For smooth projective complex varieties, building on deep analytic work of many people, Simpson established non-abelian Hodge theory by identifying suitable categories of Higgs bundles with suitable categories of vector bundles with flat connection, see \cite{SimpsonHiggs}. In \cite{OgusVologodsky}, building on a fundamental observation of \cite{BMRLocalization} (and independently Berthelot), Ogus and Vologodsky proved a characteristic $p$ counterpart to this statement. In this subsection, we explain a stack-theoretic perspective on their result; it amounts to the fact that the de Rham gerbe splits after a modification of the structure group (Corollary~\ref{cor:OVexplicit}).

Fix a smooth $k$-scheme $X$. By Proposition~\ref{prop:deRhamGerbe}, the stack $(X/k)^\dR$ is naturally a $BT_{X^{(1)}/k}^\sharp$-torsor over $X^{(1)}$. For such torsors, one has the following twisted version of Proposition~\ref{prop:reptheoryGasharp}:

\begin{proposition}[Twisted spectral decompositions via  Cartier duality]
\label{prop:AbstractCartier}
Fix a $k$-scheme $Y$, a vector bundle $E \in \mathrm{Vect}(Y)$ and a $BE^\sharp$-torsor $T \to Y$. The relative Picard stack $\mathcal{P}_T := \underline{\mathrm{Pic}}(T/Y)$ is naturally a multiplicative $\mathbf{G}_m$-gerbe over $\widehat{E^\vee}$ (the formal completion of the dual $E^\vee$ at the $0$-section). Moreover, the tautological line bundle on $T \times_Y \mathcal{P}_T$ induces a (suitably normalized) Fourier-Mukai equivalence
\[ \mathcal{D}(T) \simeq \mathcal{D}(\mathcal{P}_T)_{\deg=1} =: \mathcal{D}(\widehat{E^\vee})_{\mathcal{P}_T},\]
carrying the standard $\otimes$-product on the LHS to convolution on the RHS.
\end{proposition}
\begin{proof}[Proof sketch]
The proof can be found in \cite[Proposition 12.1.11]{Bhatt2025Apsects}. Let us simply explain here why $\mathcal{P}_T$ is a multiplicative $\mathbf{G}_m$-gerbe over $\widehat{E^\vee}$; all extensions that occur below are interpreted in the quasi-syntomic site of $Y$. The torsor $T$ determines, and is uniquely determined by, an extension
\[ BE^\sharp \to E_T \to \mathbf{Z}\]
endowed with a $BE^\sharp$-equivariant identification of $T$ with the fibre over $1$ of $E_T \to \mathbf{Z}$. Applying $\underline{\mathrm{RHom}}(-,B\mathbf{G}_m)$ and using the calculation in Proposition~\ref{prop:CohGasharp} then gives a fibre sequence
\[ B\mathbf{G}_m \to \widetilde{\mathcal{P}}_T \to \widehat{E^\vee}\]
in the derived category, so $\widetilde{\mathcal{P}}_T$ is a multiplicative $\mathbf{G}_m$-gerbe on $\widehat{E^\vee}$. There is also a tautological multiplicative line bundle on $E_T \times_Y \widetilde{\mathcal{P}}_T$; its restriction to $T \times_Y \widetilde{\mathcal{P}}_T$ then determines a map $\widetilde{\mathcal{P}}_T \to \mathcal{P}_T$, which one checks to be an isomorphism.
\end{proof}

Specializing the above to the de Rham gerbe then gives the following spectral description of the derived category of the de Rham stack (i.e., of $\mathcal{D}$-modules with nilpotent $p$-curvature, Proposition~\ref{prop:DmoduledRstackQCoh}):

\begin{corollary}[A twisted Simpson equivalence]
\label{cor:SimpsonNilpCharp}
There is a natural multiplicative $\mathbf{G}_m$-gerbe $\mathcal{P}_X$ on the formally completed cotangent bundle $\widehat{T^* X^{(1)}}$ and an equivalence
\[ \mathcal{D}((X/k)^\dR) \simeq \mathcal{D}(\widehat{T^* X^{(1)}})_{\mathcal{P}_X} \]
carrying the standard $\otimes$-product on the LHS to convolution on the RHS.
\end{corollary}

We shall refer to $\mathcal{P}_X$ as the {\em Simpson} gerbe. Passing to $\mathcal{D}$-modules via Proposition~\ref{prop:DmoduledRstackQCoh}, Corollary~\ref{cor:SimpsonNilpCharp} is a consequence of a more general statement\footnote{This more general statement can be recovered from the geometric perspective by passing to the ``sheared'' variant of Corollary~\ref{cor:SimpsonNilpCharp}, as in Remark~\ref{rmk:ShearedVariant}.} \cite{BMRLocalization}: the quasi-coherent sheaf $\mathcal{D}_{X/k}$ of crystalline differential operators on $X$ promotes naturally to an Azumaya algebra on $T^* X^{(1)}$, so there is an associated $\mathbf{G}_m$-gerbe $\mathcal{P}'_X$ on $T^* X^{(1)}$, whose restriction to $\widehat{T^*X^{(1)}}$ recovers $\mathcal{P}_X$, and an equivalence $\mathcal{D}(\mathcal{D}_{X/k}) \simeq \mathcal{D}(T^* X^{(1)})_{\mathcal{P}'_X}$.  Grosso modo, such statements are instances of non-abelian Hodge theory: up to twisting by the Simpson gerbe, they relate $\mathcal{D}_{X/k}$-modules to quasi-coherent sheaves on $T^* X^{(1)}$, i.e., Higgs modules on $X^{(1)}$. To obtain more precise statements, following \cite{OgusVologodsky}, one must trivialize the Simpson gerbe, at least on certain loci. We shall record two such statements. First, Example~\ref{ex:Fliftsplitgerbe} immediately shows:

\begin{corollary}[Splitting the Simpson gerbe via Frobenius lifts]
Assume $X/k$ comes equipped with a lift to $\Spf(A)$ together with a lift of Frobenius. Then the Simpson gerbe is trivial, and we have a natural equivalence
\[ \mathcal{D}((X/k)^\dR) \simeq \mathcal{D}(\widehat{T^* X^{(1)}}).\]
\end{corollary}

The preceding  recovers (special cases of) the local Cartier transform from \cite[Theorem 2.11]{OgusVologodsky} (see also \cite[Theorem 1.1]{TerentiukOV}). In practice, the Frobenius liftability is a  strong constraint, e.g., it never holds for a smooth projective curve of genus $\geq 2$. Another result, much more applicable in practice, is the following (see \cite[Theorem 2.8, Corollary 2.9]{OgusVologodsky} and \cite[Theorem 1.2]{TerentiukOV}):

\begin{corollary}[A non-abelian Hodge correspondence]
\label{cor:OVexplicit}
Assume $X$ is equipped with a lift to $W_2(k)$, and write $\Omega := T^*X^{(1)}$. The Simpson gerbe trivializes naturally (in the lift) after pullback along $\widehat{\Omega^\sharp} \to \widehat{\Omega}$. In particular, there is a natural (in the lift) equivalence:
\[
\left\{
\parbox{6.2cm}{\centering $\mathcal{D}_{X/k}$-modules $(M,\nabla)$ with $p$-curvature nilpotent of order $\leq p-1$}
\right\}
\xrightarrow{\sim}
\left\{
\parbox{6.2cm}{\centering Higgs modules $(N,\Theta)$ on $X^{(1)}/k$ with $\Theta$ nilpotent of order $\leq p-1$}
\right\}.
\]
\end{corollary}

In the above, $\widehat{\Omega^\sharp}$ is the completed PD-hull of the $0$-section in $\Omega = T^*X^{(1)}$; it maps to $\widehat{\Omega}$ as elements with divided powers are nilpotent in characteristic $p$.

\begin{proof}[Proof sketch]
The first part is Cartier dual to the statement that the de Rham gerbe $\pi:(X/k)^\dR \to X^{(1)}$, regarded as a $BT_{X^{(1)}/k}^\sharp$-torsor over $X^{(1)}$, splits naturally (in the lift) after pushout along the map $BT_{X^{(1)}/k}^\sharp \to BT_{X^{(1)}/k}$ forgetting divided powers. Write $(X/k)^{\dR,\gamma} \to X^{(1)}$ for this pushout. Unwinding definitions, the natural maps $(X/k)^\dR \to (X/k)^{\dR,\gamma} \to X^{(1)}$ are the transmutation of $X$ with respect to the maps 
\[ \Ga^\dR \to \Ga^{\dR,\gamma} \to F_* \Ga\]
of animated $k$-algebra stacks, where $\Ga^{\dR,\gamma}$ is obtained from pushing out the square-zero extension of Proposition~\ref{prop:SqZerodRring} along the natural map $BF_* \Ga^\sharp \to BF_* \Ga$. To prove this, consider the following map between short exact sequences of group schemes on $\mathbf{F}_p$-algebras:
\[ \xymatrix{
0 \ar[r] & \Ga^\sharp \ar[r] \ar[d]^-{\text{can}} & W \ar[r]^-{p} \ar[d]^-{\gamma_0} & W \ar[r] \ar[d]^-{\gamma_0} & \Ga \ar[r] \ar[d]^-{\mathrm{id}} & 0 \\
0 \ar[r] & \Ga \ar[r] & \Ga \ar[r]^-{p(=0)} & \Ga \ar[r]^-{\mathrm{id}} & \Ga \ar[r] & 0
} \]
The map between the middle horizontal arrows (given by multiplication by $p$) is a quasi-ideal presentation of a map of $1$-truncated animated rings that, by the above diagram, induces an identification on $\pi_0$ and the map $\text{can}:\Ga^\sharp \to \Ga$ on $\pi_1$. Consequently, we obtain an identification 
\[ \Ga^{\dR,\gamma} \simeq F_*\Ga/p := F_* \Ga \otimes_A k,\]
of animated $k$-algebra stacks, where the $k$-algebra structure comes from the second tensor factor. Under this identification, the projection map $\Ga^{\dR,\gamma} \to  \pi_0(\Ga^{\dR,\gamma}) = F_*\Ga$ is the evident one coming from the $k$-algebra structure on $F_* \Ga$. But the presentation above also gives a canonical map $F_* \Ga \to \Ga^{\dR,\gamma}$ splitting the projection to $\pi_0$; this latter map is not a $k$-algebra map (as $p$ equals $0$ for different reasons in the two animated rings) but it {\em is} a $W_2(k)$-algebra map by Lemma~\ref{lem:TwoAbarStructures} below. By transmutation, this gives the desired splitting claim.

The second part follows from the first part and Corollary~\ref{cor:SimpsonNilpCharp} by restricting the splitting further to the order $p$-infinitesimal neighbourhood of the $0$-section in $\Omega$, noting that the inclusion of this closed subscheme factors naturally over $\widehat{\Omega^\sharp} \to \widehat{\Omega}$.
\end{proof}

The following result in derived commutative algebra (which is perhaps a bit surprising from a classical perspective) was used above  and will be used again later. 

\begin{lemma}[Self-intersections of Cartier divisors]
\label{lem:TwoAbarStructures}
Let $A$ be an animated ring equipped with an invertible $A$-module $I$ and an $A$-module map $d:I \to A$. Write $\overline{A} := A/I$ for the corresponding quotient (Construction~\ref{cons:DerivedQuotients}). Let $R$ be an $\overline{A}$-algebra. Then the two natural $\overline{A}$-algebra structures on $R \otimes_A \overline{A}$, one coming from each tensor factor, induce the same $A/I^2$-algebra structure.
\end{lemma}
\begin{proof}
By base change along $\overline{A} \to R$, we may assume $R = \overline{A}$. Due to the definition of the derived quotients in Construction~\ref{cons:DerivedQuotients}, it is convenient to work geometrically. In this case, we claim that for any derived stack $X$ with Cartier divisors $D_1 \to X$ and $D_2 \to X$, the natural $X$-maps $D_1 \times_X D_2 \to D_1 \to D_1 + D_2$ and $D_1 \times_X D_2 \to D_2 \to D_1 + D_2$ are identified; note that this implies the lemma by specialization to the case $D_1=D_2$. But the universal example of this datum is $X=(\mathbf{A}^1/\mathbf{G}_m)^2$ with each $D_i$ being the pullback of $B\mathbf{G}_m \subset \mathbf{A}^1/\mathbf{G}_m$. As $D_1 \times_X D_2$ is classical here, the statement is obvious.
\end{proof}

\newpage
\section{Prismatization}
\label{sec:Prismatization}

In this section, we adapt the material of \S \ref{sec:padicdR} to mixed characteristic, i.e., to prismatic cohomology of $p$-adic formal schemes. We do not aim to be comprehensive, e.g., absolute prismatic cohomology is not discussed here. Rather, our goal is to convince the reader that the stack-theoretic perspective effortlessly translates the characteristic $p$ results of \S \ref{sec:padicdR}, especially on non-abelian Hodge theory, to mixed characteristic. The key change is to replace the ring stack $\Ga^\dR$ from \S \ref{sec:padicdR} with the prismatization ring stack $\Ga^\Prism$ (\S \ref{ss:RelPrism}). The constructions of \S \ref{sec:QuantNAHT} shall build on this material.

As in \S \ref{sec:padicdR}, there are several benefits of the stack-theoretic perspective. Let us highlight the following instances in this section:

\begin{itemize}
    \item The crystalline comparison theorem for prismatic cohomology \cite{BhattScholzePrism} (and thus the slightly surprising Frobenius descent property of crystalline cohomology over a crystalline prism) is encoded in an isomorphism $\varphi^* \Ga^\Prism \simeq \Ga^\dR$ of ring stacks over any crystalline prism (Example~\ref{ex:Prismatization} (2)).

    \item The Hodge--Tate comparison theorem for prismatic cohomology (the mixed characteristic version of the Cartier isomorphism) follows from a gerbe structure on the Hodge--Tate stack (Theorem~\ref{thm:HTstackformal} (1) and Remark~\ref{rmk:HTcomparison}), which in turn comes via deformation theory from a square-zero extension structure on $\Ga^\HT$, the restriction of $\Ga^\Prism$ to the Hodge--Tate locus (Proposition~\ref{prop:SqZeroHT}).

    \item Combining the aforementioned gerbe structure with Cartier duality leads to spectral decomposition for the category of Hodge--Tate crystals (Theorem~\ref{thm:HTstackformal} (2) and (3)); this is a prismatic variant of Corollary~\ref{cor:SimpsonNilpCharp}, and we regard it as the primordial version of the $p$-adic Simpson equivalence (e.g., in \S \ref{sec:Simpson}, it provides the starting point for the rigid analytic Simpson equivalence).

    \item The structure theory of the prismatization $(\mathbf{G}_m/A)^\Prism$ of the multiplicative group naturally leads to the prismatic logarithm, and thus to the theory of Chern classes in prismatic cohomology (Example~\ref{ex:Prismatization} (3); see also \cite[\S 2]{BhattLurieAbsPrismaticCoh} and \cite[\S 4]{BhattLuriePrismatization}).
\end{itemize}

With the exception of $\varphi$-crystalline prisms in \S \ref{ss:Prisms}, all the material in this section is expository. Some references are \cite{BhattLurieAbsPrismaticCoh,BhattLuriePrismatization,DrinfeldPrismatization,ALBH1}.

\subsection{Prisms}
\label{ss:Prisms}

In this subsection, we recall the notion of a prism. Unlike standard expositions, we shall directly work in the animated setting, as in \cite[\S 2]{BhattLuriePrismatization}:

\begin{definition}[Animated prisms]
An {\em animated prism} consists of a pair $(A,d:I \to A)$ (or sometimes just $(A,I)$), where $A$ is an animated $\delta$-ring and $d:I \to A$ is a generalized Cartier divisor such that $A$ is $(p,I)$-complete and the following holds: for each map $\alpha:A \to k$ with $k$ a perfect characteristic $p$ field, the ideal of $W(k)$ generated by the image of $I$ under the natural $\delta$-lift $\widetilde{\alpha}:A \to W(k)$ of $\alpha$ is $pW(k)$. 

A morphism  $(A,I) \to (B,J)$ of animated prisms is given by a map $A \to B$ of animated $\delta$-rings and an identification $(I \to A) \otimes_A B \simeq (J \to B)$ of generalized Cartier divisors. 
\end{definition}

The classical notion of a prism corresponds to an animated prism $(A,I)$ where $A/I$ is discrete; this implies $A$ is itself discrete (but the converse is not true). If $A/I$ additionally has bounded $p^\infty$-torsion, then we say $(A,I)$ is bounded.

\begin{remark}[Why animated prisms?]
There are several reasons we prefer to work in the animated setting in these notes.  At a conceptual level, the statements are much cleaner in the animated world. For example, pulling back along the Frobenius defines an endofunctor $(A,I) \mapsto (A,\varphi^*I)$ of the $\infty$-category of animated prisms (but it might not preserve classical prisms). Another example: the $\infty$-category of prisms over $(A,I)$ is identified with the $\infty$-category of $(p,I)$-complete animated $\delta$-$A$-algebras. There is also a practical reason for this choice, related to the notion of a $\varphi$-crystalline prism introduced below (Example~\ref{ex:PhiCrystallizeFontainePrism}): namely, the most important example (Example~\ref{ex:PhiCrystallizeFontainePrism}) of such a prism is not bounded. Absent boundedness, it is more natural to work directly in the animated context when doing formal geometry\footnote{For example, if $R$ is a $p$-complete ring whose $p$-power torsion is not bounded, then the classical formal scheme $\mathrm{Spf}(R)$ is poorly behaved, e.g., its quasi-coherent derived category does not coincide with the $p$-complete derived category of $R$, and the natural notion of flatness for $p$-complete $R$-algebras and affine formal schemes over $\mathrm{Spf}(R)$ also fails to coincide; both these issues disappear for the derived $p$-adic formal scheme attached to $R$.}.
\end{remark}

The following are some key examples of prisms:

\begin{example}
    \begin{enumerate}
        \item For any $p$-torsionfree $p$-complete $\delta$-ring $A$, the pair $(A,(p))$ is a prism. More generally, for any animated $p$-complete $\delta$-ring $A$, the pair $(A,p:A \to A)$ is a prism; such objects are called {\em (animated) crystalline prisms}.

        \item An animated prism $(A,I)$ is called {\em perfect} if the Frobenius on $A$ (or equivalently $A/p$) is an isomorphism. The category of perfect prisms is equivalent to the category of perfectoid rings via the functor $(A,I) \mapsto A/I$ (and in fact this is a reasonable definition of the latter); the inverse is given by Fontaine's construction $R \mapsto (A_{\inf}(R),\ker(\theta))$. In particular, all such animated prisms are classical and bounded.

        \item Fix a complete discretely valued extension $K/\mathbf{Q}_p$ with perfect residue field $k$ and a uniformizer $\pi \in \calO_K$. This datum determines {\em the Breuil--Kisin prism} $(W(k)\llbracket u \rrbracket, (E(u)))$; here the $\delta$-structure on $W(k)\llbracket u \rrbracket$ is uniquely determined by requiring $\varphi(u)=u^p$, and $E(u)$ is the monic generator of the kernel of $W(k)\llbracket u \rrbracket \xrightarrow{u \mapsto \pi} \calO_K$. 

        \item Consider the $\delta$-ring $A := \mathbf{Z}_p\llbracket q-1\rrbracket$, with $\delta$-structure determined by $\varphi(q)=q^p$. Then the pair $(A,([p]_q))$ is a prism, where $[p]_q = 1+q+...+q^{p-1}$ is the $q$-analog of $p$; this is called the {\em $q$-de Rham prism}.
    \end{enumerate}
\end{example}

In addition to these standard examples, for the material in \S \ref{sec:QuantNAHT}, we shall need another slightly more exotic class of prisms that will be studied in \cite{BHPVQuanNAHT}. To motivate these, recall the following relation between $\delta$-structures and divided powers (see \cite[\S 2.5]{BhattScholzePrism} for more):

\begin{lemma}
\label{lem:DeltaDividedPowers}
Let $A$ be a $p$-torsionfree $\delta$-ring and $x \in A$. Then $x$ admits divided powers in $A$ if and only if $\varphi(x) \in pA$. 
\end{lemma}

Thanks to the above, the following definition can be regarded as a derived analog of the class of prisms $(A,I)$ where $I$ admits divided powers in $A$.

\begin{definition}[$\varphi$-crystalline prisms]
\label{def:phicrystallineprism}
A {\em $\varphi$-crystalline structure} on an animated prism $(A,I)$ is an isomorphism $(A,\varphi^*(d):\varphi^*I \to A) \simeq (A,p:A \to A)$ of animated prisms. Such a structure is unique (when it exists) if $A/\varphi^*I$ or equivalently $A/p$ is discrete; in this case, we simply say $(A,I)$ is a {\em $\varphi$-crystalline prism}.
\end{definition}

Clearly, any crystalline prism has a $\varphi$-crystalline structure. Moreover, there is a universal mechanism for turning any animated prism into a $\varphi$-crystalline one. More precisely, the forgetful functor from the $\infty$-category of $\varphi$-crystalline animated prisms to the $\infty$-category of animated prisms has a left adjoint, given by $(A,d:I \to A) \mapsto (B,d_B:J\to B)$, where $B=A\{\frac{\varphi^*I}{p}\}^{\wedge}_{(p,I)}$ and $d_B$ is the base change of $d$. Other than crystalline prisms, the following is the main example of relevance in the sequel:

\begin{example}[$\varphi$-crystallizing the Fontaine prism]
\label{ex:PhiCrystallizeFontainePrism}
Fix a $p$-torsionfree perfectoid ring $R$, corresponding to a perfect prism $(A_{\inf}(R),I=\ker(\theta))$. Let $A_\crys(R) = D_{I}(A_{\inf}(R))$ be the $p$-completed divided power envelope of $I \subset A_{\inf}(R)$. Then $(A_\crys(R),IA_\crys(R))$ is a $\varphi$-crystalline prism. This prism is not crystalline; in fact, there is a natural map $A_\crys(R)/IA_\crys(R) \to R$ from the universal property of divided power envelopes\footnote{\label{fn:AcrysNbhdFontainePrism}Geometrically, one can regard $\mathrm{Spa}(A_\crys(R)[1/p],A_\crys(R)) \to \mathrm{Spa}(A_{\inf}(R))^\an$ as a rational open subset that contains the Hodge--Tate point and its positive Frobenius translates, but not its Frobenius pre-image. We were motivated to consider such prisms whilst contemplating quantization of rational non-abelian $p$-adic Hodge theory in a small neighbourhood of the Hodge--Tate point (carried out in \S \ref{ss:RigAnQuanNAHT}) from an integral perspective.}. Moreover, it is also not bounded: this reduces to the statement that the mod $(x)$ reduction of the $\delta$-ring $\mathbf{Z}_p\{x,\frac{\varphi(x)}{p}\} = D_{(x)}\mathbf{Z}_p\{x\}$ has unbounded $p$-power torsion.
\end{example}

\begin{remark}[The universal $\varphi$-crystalline prism]
The Cartier--Witt stack $\Spf(\mathbf{Z}_p)^\Prism$ from \cite{BhattLurieAbsPrismaticCoh,DrinfeldPrismatization} can be regarded as the non-existent final object in the opposite of the category of animated prisms. With the same abuse of notation, the universal $\varphi$-crystalline prism is identified with $\Ga^\sharp/\mathbf{G}_m^\sharp$, with the map to $\mathrm{Spf}(\mathbf{Z}_p)^\Prism$ provided informally by sending $\epsilon \in \Ga^\sharp \subset W$ to the animated prism defined by the distinguished element $p+\epsilon$. (We thank Jacob Lurie for this remark.)
\end{remark}

\subsection{Relative prismatization} 
\label{ss:RelPrism}

In this subsection, following \cite{DrinfeldPrismatization,BhattLuriePrismatization}, we introduce and study prismatization functors relative to a fixed prism via transmutation. 

\begin{notation}
\label{not:RelPrism}
Fix an animated prism $(A,I)$. Write $\overline{A} = A/I$ and endow $A$ and $\overline{A}$ with the $(p,I)$-adic topology. Write $\CAlg^\an_{\Spf(A)}$ for the $\infty$-category of $(p,I)$-nilpotent animated $A$-algebras, and similarly over $\Spf(\overline{A})$. As in \S \ref{sec:padicdR}, a (smooth) scheme over $\Spf(A)$ or $\Spf(\overline{A})$ is a (smooth) schematic map to it.
\end{notation}

The prism $(A,I)$ determines the following ring stack:

\begin{construction}[The ring stack $\mathbf{G}_a^\Prism$]
\label{cons:GaPrism}
Fix $R \in \CAlg^\an_{\Spf(A)}$. The structure map $A \to R$ refines uniquely to a $\delta$-map $A \to W(R)$. Using the $(p,I)$-nilpotence of $R$, one checks that $W(R)/I$ is $p$-nilpotent, and thus naturally an object of $\CAlg^\an_{\Spf(\overline{A})}$. Write $\Ga^\Prism$ for the resulting functor
\[\CAlg^\an_{\Spf(A)} \to \CAlg^\an_{\Spf(\overline{A})}, \quad R \mapsto W(R)/I.  \]
Informally, we refer to this as an animated $\overline{A}$-algebra stack on $\CAlg^\an_{\Spf(A)}$.
\end{construction}

Transmuting $\Ga^\Prism$ yields:

\begin{definition}[The relative prismatization]
For any $X \in \Fun(\CAlg^\an_{\Spf(\overline{A})},\Sp)$, write $(X/A)^\Prism$ for its transmutation with respect to $\Ga^\Prism$. Explicitly, $(X/A)^\Prism$ is the composition
\[ \CAlg^\an_{\Spf(A)} \xrightarrow{\Ga^\Prism} \CAlg^\an_{\Spf(\overline{A})} \xrightarrow{X} \Sp, \]
i.e., for $R \in \CAlg^\an_{\Spf(A)}$, we have
\[ (X/A)^\Prism(R) := X(W(R)/I) := \mathrm{Map}_{\Spf(\overline{A})}(\Spf(W(R)/I), X). \]
Frobenius on the Witt vectors defines a map $\varphi_{X/A}:(X/A)^\Prism \to \varphi_A^* (X/A)^\Prism$ over $\Spf(A)$ that we informally call the relative Frobenius; the corresponding $\varphi_A$-linear endomorphism of $(X/A)^\Prism$ yields a $\delta$-structure on $(X/A)^\Prism$ compatible with the given one on $A$. 
\end{definition}

\begin{remark}
In line with Remark~\ref{rmk:VadicTopdR}, it would have been more natural to replace $p$-adic formal scheme $\mathrm{Spf}(W(R)/I)$ appearing above with the non-$p$-adic formal scheme $\colim_n \mathrm{Spf}(W_n(R)/I)$ in the above discussion. Again, for the purposes of transmuting schemes, this difference is irrelevant by Tannaka duality.
\end{remark}

\begin{example}
\label{ex:Prismatization}
\begin{enumerate}
\item {\em The final object:} We have $(\Spf(\overline{A})/A)^\Prism = \Spf(A)$

\item {\em The crystalline case}: Say $k$ is a perfect field of characteristic $p$, and let $(W(k),(p))$ be the corresponding perfect prism. Then $\varphi^* \Ga^\Prism \simeq  \Ga^\dR$ as animated $k$-algebra stacks on $\CAlg^\an_{\Spf(A)}$.  Consequently, for any $k$-scheme $X$, we have $\varphi^* (X/W(k))^\Prism \simeq (X/W(k))^\crys$.

More generally, the above discussion applies over $(A,(p))$ for any $p$-complete animated $\delta$-ring $A$. The appearance of $\varphi$-pullback then shows that the value on $A$ of the crystalline cohomology of smooth  $A/p$-schemes canonically descends along Frobenius (compare \cite[Example 1.9 (1)]{BhattScholzePrism}).

\item {\em The multiplicative group:} Consider the multiplicative group $\mathbf{G}_{m,\Spf(\overline{A})}$ over $\Spf(\overline{A})$. By functoriality, its prismatization $(\mathbf{G}_{m,\Spf(\overline{A})}/A)^\Prism$ is naturally a commutative group stack over $\Spf(A)$. Moreover, using the natural $\delta$-structure on $\mathbf{G}_m$, the forthcoming Proposition~\ref{prop:CharPrismatization} gives a $\delta$-map $\mathbf{G}_{m,\Spf(A)} \to (\mathbf{G}_{m,\Spf(\overline{A})}/A)^\Prism
$ that can be shown to be a quasi-syntomic cover. Consequently, there is a fibre sequence
\[ G \to \mathbf{G}_{m,\Spf(A)} \to(\mathbf{G}_{m,\Spf(\overline{A})}/A)^\Prism \]
where $G$ is a commutative $\delta$-group scheme over $\Spf(A)$. This group scheme $G$ turns out to be quite interesting: it is the Cartier dual of a $1$-dimensional formal group over $\Spf(A)$ that interpolates between the multiplicative and additive formal groups, and its logarithm map agrees with the prismatic logarithm, first introduced in \cite[\S 2]{BhattLurieAbsPrismaticCoh} by a more explicit construction with the aim of defining Chern classes in prismatic cohomology; we refer to \cite{DrinfeldFormalGroup} and \cite[\S 4]{BhattLuriePrismatization} for more.
\end{enumerate}
\end{example}

Example~\ref{ex:Prismatization} (2) shows that the prismatization construction is a generalization of the crystallization construction and hence the de Rham stack construction too. It enjoys the same  formal properties as the latter (Proposition~\ref{prop:dRStackBasicProp}):

\begin{proposition}[Properties of relative prismatization]
\label{prop:RelPrismStackBasicProp}
\begin{enumerate}
    \item If $X/\Spf(\overline{A})$ is a $p$-adic formal scheme, then $(X/A)^\Prism$ is an fpqc sheaf.

    \item  The relative prismatization construction $(-/A)^\Prism:\Fun(\CAlg^\an_{\Spf(\overline{A})},\Sp) \to \Fun(\CAlg^\an_{\Spf(A)},\Sp)$ commutes with limits and preserves representable quasi-syntomic covers.

\item If $f:U \to X$ is a representable \'etale map in $\Fun(\CAlg^\an_{\Spf(\overline{A})},\Sp)$, then $(U/A)^\Prism \to (X/A)^\Prism$ is also representable and \'etale.

\item Fix a map $(A,I) \to (B,J)$ of prisms. For any $X \in \Fun(\CAlg^\an_{\Spf(\overline{A})},\Sp)$ with restriction $X_{\overline{B}} \in \Fun(\CAlg^\an_{\Spf(\overline{B})},\Sp)$, there is a natural base change identification $(X_{\overline{B}}/B)^\Prism \simeq (X/A)^\Prism \times_{\Spf(A)} \Spf(B)$ in $\Fun(\CAlg^\an_{\Spf(B)},\Sp)$.
\end{enumerate}
\end{proposition}

The above properties ensure that $(X/A)^\Prism$ has good geometric properties, e.g., if $X/\overline{A}$ is a smooth affine $p$-adic formal scheme, then $(X/A)^\Prism$ is the quotient (as quasi-syntomic sheaves) of an affine pro-smooth scheme over $\Spf(A)$ by an affine equivalence relation. The last property above, in particular, allows us to control the relative prismatization  well in terms of the relative de Rham stack: base changing along the ``crystallization'' map $A\to A\{\frac{I}{p}\}^{\wedge}_{(p)}$ is conservative and carries $(X/A)^\Prism$ to a crystalline stack. However, there is one respect in which the prismatization and crystallization differ: while lifts provide charts for the latter (Remark~\ref{rmk:ChartCrystallization}), one typically needs more structure to obtain charts for the prismatization.

\begin{proposition}[Charts from $\delta$-lifts]
\label{prop:CharPrismatization}
Say $X/\Spf(\overline{A})$ is a smooth scheme. Given a lift $\widetilde{X}/\Spf(A)$ and a $\delta$-structure on $\widetilde{X}$ compatible with the one on $A$, there is a natural $\delta$-map $\pi:\widetilde{X} \to (X/A)^\Prism$ which is a quasi-syntomic cover. 
\end{proposition}
\begin{proof}[Proof sketch]
We shall explain the construction of the map; the covering property needs an additional argument. Fix $R \in \CAlg^\an_{\Spf(A)}$ and a point $\alpha \in \widetilde{X}(R)$, regarded as a map $\alpha:\Spec(R) \to \widetilde{X}$. As $\widetilde{X}$ carries a $\delta$-structure compatible with the one on $A$, the map $\alpha$ refines uniquely to a $\delta$-map $\widetilde{\alpha}:\Spf(W(R)) \to \widetilde{X}$, where $W(R)$ is endowed with the $(p,I)$-adic topology. Base changing to $\Spf(\overline{A})$ then gives a map $\Spf(W(R)/I) \to X$, and thus a point $\pi(\alpha) \in (X/A)^\Prism(R)$. 
\end{proof}

Finally, let us mention that a primary  reason for constructing the prismatization was geometrization: realizing the category of prismatic crystals as the category of quasi-coherent sheaves on a stack. 

\begin{proposition}[{\cite[Proposition 8.15]{BhattLuriePrismatization}}]
    \label{prop:PrismaticCrystalsAsSheaves}
Say $X/\Spf(\overline{A})$ is a smooth scheme. Assume $\overline{A}$ has bounded $p^\infty$-torsion. Then there is a natural identification
\[ \mathcal{D}( (X/A)^\Prism) \simeq \mathcal{D}_\crys( (X/A)_\Prism, \mathcal{O}_\Prism) := \lim_{(B,J) \in (X/A)_\Prism} {\mathcal{D}}_{(p,J)\text{-comp}}(B),\]
where $(X/A)_\Prism$ is the relative prismatic site of \cite{BhattScholzePrism}.
\end{proposition}

\subsection{Non-abelian Hodge theory for formal schemes}
\label{ss:naPHTformal}

In this subsection, continuing with Notation~\ref{not:RelPrism}, we study the geometry and sheaf theory of the restriction of prismatization to the Hodge--Tate locus, i.e., to test objects over $\Spf(\overline{A})$. As in \S \ref{ss:NAHTcharp}, this leads to concrete statements about quasi-coherent sheaves that have a non-abelian Hodge theory flavour.

\begin{definition}[The Hodge--Tate stacks]
Let $\Ga^\HT$ be the endofunctor of $\CAlg^\an_{\Spf(\overline{A})}$ obtained by restricting the ring stack $\Ga^\Prism$ (Construction~\ref{cons:GaPrism}) along $\CAlg^\an_{\Spf(\overline{A})} \to \CAlg^\an_{\Spf(A)}$, i.e., $\Ga^\HT(R) = W(R)/I$. For any $X \in  \Fun(\CAlg^\an_{\Spf(\overline{A})},\Sp)$, the {\em Hodge--Tate stack of $X$ relative to $A$}, denoted $(X/A)^\HT$, is the transmutation of $X$ with respect to $\Ga^\HT$, i.e., the composition
\[ \CAlg^\an_{\Spf(\overline{A})} \xrightarrow{\Ga^\HT} \CAlg^\an_{\Spf(\overline{A})} \xrightarrow{X} \Sp,\]
so, for $R \in \CAlg^\an_{\Spf(\overline{A})}$, we have
\[ (X/A)^\HT(R) := X(W(R)/I) := \mathrm{Map}_{\Spf(\overline{A})}(\Spf(W(R)/I), X). \]
Thus, $(X/A)^\HT$ is also the restriction of $(X/A)^\Prism$ along $\Spf(\overline{A}) \to \Spf(A)$.
\end{definition}

Of course, the properties of the relative prismatization (Proposition~\ref{prop:RelPrismStackBasicProp}) pass to analogous properties of the Hodge--Tate stack. The principal new feature is that $(X/A)^\HT$ lives over $X$ and is a gerbe for an explicit group scheme;  this is the prismatic variant of the Frobenius linearity and gerbe structure on the de Rham stack in characteristic $p$ (Proposition~\ref{prop:deRhamGerbe}). Unsurprisingly, the starting point is a variant of Proposition~\ref{prop:SqZerodRring}:

\begin{proposition}[$\Ga^\HT$ as a square-zero extension]
\label{prop:SqZeroHT}
$\Ga^\HT$ is naturally and uniquely a square-zero extension of $\Ga$ by $B\Ga^\sharp\{1\}$ in animated $\overline{A}$-algebra stacks, where we write $M\{1\} = M \otimes_{\overline{A}} I/I^2$ for any $\overline{A}$-module $M$. 
\end{proposition}
\begin{proof}
By a left Kan extension argument, it suffices to show the statement after restricting to classical $p$-nilpotent $\overline{A}$-algebras. Moreover, by descent, it suffices to restrict attention to rings $R$ where $F:W(R) \to W(R)$ is surjective. For such an $R$, we shall show that the natural map
\[ (I \xrightarrow{d} A) \otimes_A W(R) = I \otimes_A W(R) \xrightarrow{d_R} W(R) \]
has cokernel $W(R) \to R$ (given by the $0$-th Witt component) and kernel $I \otimes_A \Ga^\sharp(R) \subset I \otimes_A W(R)$ (given by the Proposition~\ref{prop:WittGasharp}); this implies  the desired claim as any $1$-truncated animated ring $B$ is naturally and uniquely a square-zero extension of $\pi_0(B)$ by $\pi_1(B)$. This claim can be checked locally, so we may assume $I=(d)$ is principal. Now the image of $d$ under $A \to W(A) \to W(\overline{A})$ has the form $V(u)$ for a unit $u \in W(\overline{A})$, as $(A,I)$ is a prism. The identity $x \cdot V(u) = V(Fx \cdot u)$ and the surjectivity assumption on $F:W(R) \to W(R)$ then show that $\mathrm{coker}(d_R) = W(R)/VW(R) = R$ and $\ker(d_R) = \ker(F) = \Ga^\sharp(R)$, as wanted.
\end{proof}

Transmuting the above, we obtain:

\begin{theorem}
\label{thm:HTstackformal}
Fix a smooth scheme $X/\Spf(\overline{A})$. 
\begin{enumerate}
    \item {\em The Hodge--Tate gerbe:} There is a natural map $\pi:(X/A)^\HT \to X$ which is a $BT_X^\sharp\{1\}$-torsor; we call this the {\em Hodge--Tate gerbe}.
    \item {\em Liftings give reductions:} The pushout of the torsor in (1) along the natural map $BT_X^\sharp\{1\} \to BT_X\{1\}$ is the $BT_X\{1\}$-torsor of liftings of $X/\overline{A}$ to $A/I^2$.
 In particular, each lift of $X$ to $A/I^2$ provides a reduction of structure group of the Hodge--Tate gerbe in (1) along the natural map 
 \[ (T_X\{1\}/X)^\dR := T_X\{1\}/T_X^\sharp\{1\} \to BT_X^\sharp\{1\}\]
 of commutative group stacks over $X$. 

    \item  {\em A twisted Simpson equivalence:} There is a natural multiplicative $\mathbf{G}_m$-gerbe $\mathcal{P}_X$, the {\em Simpson gerbe}, on the formal completion $\widehat{T^* X}\{-1\}$ of $T^*X\{-1\}$ along its zero section, and an equivalence
    \begin{equation} 
    \label{eq:twistedSimpsonFS}
    \tag{SimpFS}
    \mathcal{D}( (X/A)^\HT) \simeq \mathcal{D}(\widehat{T^*X}\{-1\})_{\mathcal{P}_X},
    \end{equation}
    carrying the standard $\otimes$-product on the LHS to convolution on the RHS. 
\end{enumerate}
\end{theorem}

\begin{proof}
The proofs of these statements are entirely parallel to the characteristic $p$ case. More precisely, (1)  is deduced from Proposition~\ref{prop:SqZeroHT} by deformation theory as in Proposition~\ref{prop:deRhamGerbe}, (2) can be proven using the proof strategy of Corollary~\ref{cor:OVexplicit}, while (3) is proven as in Corollary~\ref{cor:SimpsonNilpCharp}, with the gerbe coming from dualizing the torsor as in Proposition~\ref{prop:AbstractCartier}.
\end{proof}

\begin{remark}[The Hodge--Tate comparison]
\label{rmk:HTcomparison}
As in Corollary~\ref{cor:Cartier}, part (1) of Theorem~\ref{thm:HTstackformal} implies the {\em Hodge--Tate comparison}: there is a natural identification 
\[ \bigoplus_i R^i \pi_* \calO_{(X/A)^\HT} \simeq \bigoplus_i \Omega^i_{X/\overline{A}}\{-i\}\] 
of graded coherent $\calO_X$-algebras. If $X$ comes equipped with a lift $\widetilde{X}/\Spf(A)$ carrying a compatible $\delta$-structure, then the mod $I$ reduction of the charts from Proposition~\ref{prop:CharPrismatization} split the Hodge--Tate gerbe, so one can upgrade the preceding cohomological statement to a derived one: the gerbe $\mathcal{P}_X$ is trivialized, and there is a natural identification
\[ R\pi_* \calO_{(X/A)^\HT} \simeq \bigoplus_i \Omega^i_{X/\overline{A}}\{-i\}\]
of commutative algebra objects in $\mathcal{D}(X)$. 
\end{remark}

\begin{remark}[The connection on the Simpson gerbe]
    \label{rmk:SimpsonGerbeConn}
Fix $X/\Spf(\overline{A})$ smooth with a lift $\widetilde{X}/\Spf(A/I^2)$. The Simpson gerbe $\mathcal{P}_X$ in Theorem~\ref{thm:HTstackformal} (3) is obtained via duality from the Hodge--Tate torsor in Theorem~\ref{thm:HTstackformal} (1). Now Theorem~\ref{thm:HTstackformal} (2) gives an additional property of the Hodge--Tate torsor in the presence of the lift $\widetilde{X}$. This statement dualizes to show that the $\mathbf{G}_m$-gerbe $\mathcal{P}_X$ admits a flat connection relative to $X$. In fact, one can show that $\mathcal{P}_X$ admits a flat connection relative to $\overline{A}$, i.e., it descends to the (sheared) de Rham stack of $\widehat{T^*X}\{-1\}$ relative to $\overline{A}$. (This result is due to Vadim Vologodsky.)
\end{remark}

\newpage

\section{The $p$-adic Simpson equivalence}
\label{sec:Simpson}

In this section, we explain an application of \S \ref{ss:naPHTformal} to rigid geometry over a perfectoid field; in the process, we  also introduce some notions that will be relevant to the forthcoming quantization discussion.

\begin{notation}
Let $C/\mathbf{Q}_p$ be a perfectoid field. Then $\calO_C$ is a rank $1$ perfectoid valuation ring with maximal ideal $\mathfrak{m}_C$, and we have the corresponding perfect prism $(A_{\inf}(\calO_C),\ker(\theta))$. 
\end{notation}

Our goal is to explain the $p$-adic Simpson theory (pioneered by Faltings \cite{FaltingsSimpson}) for a smooth rigid space $X/C$, relating local systems (topological side) on $X$ to Higgs bundles (Dolbeault side). The first step is to move from the topological side towards the coherent side by using the pro-\'etale site $X_{\proet}$ with its natural (completed) structure sheaf $\widehat{\mathcal{O}}_X$ (see \cite{ScholzePAdicHodge}): there is a natural exact $\otimes$-functor
\begin{equation}
    \label{eq:RH}
    \tag{RH}
\mathrm{Loc}(X;C) := \left\{
\parbox{2.6cm}{\centering $C$-local systems on $X$}
\right\}
\xrightarrow{L \mapsto L \otimes_C \widehat{\mathcal{O}}_X}
\left\{
\parbox{2.6cm}{\centering pro-\'etale vector bundles on $X$}
\right\} =: \mathrm{Vect}(X_{\proet},\widehat{\mathcal{O}}_X).
\end{equation}
By the primitive comparison theorem, this functor is fully faithful in the geometrically interesting situations\footnote{For example, if $C$ is algebraically closed and $X$ is additionally proper, the map is fully faithful, even in the derived sense, on $C$-local systems admitting an $\mathcal{O}_C$-lattice.}, so the question becomes: how to relate the target to Higgs bundles? We shall prove that the target is naturally equivalent to the category of twisted Higgs bundles, where the twist refers to a naturally constructed multiplicative $\mathbf{G}_m$-gerbe $\mathcal{P}_X \to T^*X\{-1\}$ --- the Simpson gerbe --- on the Tate-twisted cotangent bundle; we call this the {\em $p$-adic Simpson equivalence}. In fact, our methods naturally yield the derived enhancement: in Theorem~\ref{thm:padicSimpson}, we construct an equivalence
\[ \Perf(X_{\proet},\widehat{\calO}_X) \simeq \mathcal{D}_{X\text{-}\text{Perf}}(T^*X\{-1\})_{\mathcal{P}_X}.\]
Here $\mathcal{D}_{X\text{-}\text{Perf}}(T^*X\{-1\})$ formalizes the notion of a perfect complex with a Higgs structure, and the equivalence 
exchanges $\otimes$-product on the LHS with convolution on the RHS.

In constructing this equivalence, we  adopt Raynaud's point of view: rigid analytic geometry over $C$ is the generic fibre of formal geometry over $\Spf(\calO_C)$. More precisely, we obtain the above equivalence as, roughly, the generic fibre of the Simpson equivalence for formal schemes in Theorem~\ref{thm:HTstackformal} (3). Thus the stack-theoretic perspective is essential here twice over: the Simpson gerbe enters the very statement, and the formal scheme input is proven stack-theoretically.

\begin{remark}[Trivializing the gerbe]
    \label{rmk:SimpsonTrivGerbe}
When $C=\overline{C}$, $X/C$ is proper smooth, and one makes two auxiliary choices, the gerbe  $\mathcal{P}_X$ can be compatibly trivialized on the relevant subvarieties of $T^*X\{-1\}$, so one obtains an equivalence
\begin{equation}
    \label{eq:untwistedSimpson}
    \tag{Simp}
\begin{split}
\mathrm{Vect}(X_{\proet},\widehat{\mathcal{O}}_X) &\simeq \mathrm{Higgs}(X,\Omega^1_{X/C}\{-1\}) \\
&:= \left\{ (E,\Theta) \;\middle|\; \parbox{7cm}{\centering $E \in \mathrm{Vect}(X)$ and $\Theta:E \to E \otimes \Omega^1_{X/C}\{-1\}$ such that $\Theta$ is $\mathcal{O}_X$-linear and $\Theta \wedge \Theta = 0$} \right\},
\end{split}
\end{equation}
and similarly for perfect complexes. This was the form of the $p$-adic Simpson equivalence originally established by Faltings in dimension $1$ \cite{FaltingsSimpson} and Heuer \cite{HeuerSimpson} in general. We offer a simple deduction, due to Alexander Petrov, of the untwisted statement from the twisted one in Theorem~\ref{thm:padicSimpsonUntwisted}.
\end{remark}

\begin{remark}[The essential image of local systems]
In the situation of Remark~\ref{rmk:SimpsonTrivGerbe}, the major open question in the area (due to Faltings \cite[\S 5]{FaltingsSimpson}; see also \cite[Conjecture 6]{DeningerWernerII}) is the description of the essential image of $C$-local systems admitting a lattice under the composition of \eqref{eq:RH} and \eqref{eq:untwistedSimpson}. In particular does it coincide with semistable Higgs bundles with vanishing Chern classes when $C=\mathbf{C}_p$ and $X$ is projective? Both the restrictions are necessary, already in rank $1$, by \cite{HeuerRankOne}. We have nothing new to offer in this direction. 
\end{remark}

The work described in this section is joint work with Mingjia Zhang, see \cite{MingjiaSimpsonNotes} and \cite[\S 12]{Bhatt2025Apsects}. It is inspired by a number of previous works in this area, most notably \cite{FaltingsSimpson} and \cite{HeuerSimpson} as indicated above, but also \cite{ALBH1,AbbesGrosTopos,OgusVologodsky,TerentiukOV} (and \cite{MinWangPrismaticSimpson} over discretely valued fields). The only new mathematics in this section is the observation that the theory of the Simpson gerbe on smooth formal schemes already determines it on all smooth rigid spaces (Proposition~\ref{prop:SimpsonGerbeRigid}), and Petrov's simple proof of the splitting property in Remark~\ref{rmk:SimpsonTrivGerbe} (Theorem~\ref{thm:padicSimpsonUntwisted}).

\subsection{The small correspondence}

In this subsection, we brutally invert $p$ in the formal scheme Simpson equivalence to obtain a first approximation to the Simpson equivalence for rigid spaces. Thus, assume we can find a smooth scheme $\mathcal{X}/\Spf(\calO_C)$ with rigid generic fibre $X/C$. To interpret the result obtained by inverting $p$ in \eqref{eq:twistedSimpsonFS} in terms of the generic fibre, we need to understand each side separately. 

\begin{enumerate}
\item \label{item:SmallCrystalSide} {\em The crystal side:} 
Recall that the pro-\'etale site $X_{\proet}$ of $X$ has a basis consisting of maps $\mathrm{Spa}(R,R^+) \to X$ with source an affinoid perfectoid space. Any such map determines (and is determined by) a map $\mathrm{Spf}(R^+) \to (\mathcal{X}/A_{\inf}(\calO_C))^\HT$. Pulling back Hodge--Tate crystals along this construction and inverting $p$ therefore gives a functor
\[ \Perf( (\mathcal{X}/A_{\inf}(\calO_C))^\HT)[1/p] \to \Perf(X_{\proet},\widehat{\calO}_X). \]
The isogeny theorem for prismatic cohomology implies that this functor is fully faithful (see \cite[Proposition 12.3.3]{Bhatt2025Apsects}). Write $\Perf(X_{\proet},\widehat{\calO}_X)_{\sm} \subset \Perf(X_{\proet},\widehat{\calO}_X)$ for its essential image; objects in this category are called {\em small}.

\item {\em The Dolbeault side:} The integral model $\mathcal{X}$ of the rigid space $X$ provides a norm on the vector bundle $T^*X$ with the property that the rigid generic fibre  of the formal completion $\widehat{T^* \mathcal{X}}$ agrees with the open unit disc bundle $T^*X_{<1} \subset T^*X$. Consequently, inverting $p$ in the first half of Theorem~\ref{thm:HTstackformal} (3) (and twisting by $\calO_C\{-1\}$) gives a natural multiplicative gerbe on $T^*X\{-1\}_{< 1}$ banded by the ``unit tire''\footnote{We thank Brian Conrad for suggesting this terminology.} $\mathbf{G}_{m,\eta} = \{|x|=1\}$; let $\mathcal{P}_{\mathcal{X},\eta}$ be its pushout along $\mathbf{G}_{m,\eta} \subset \mathbf{G}_m$, regarded as a multiplicative $\mathbf{G}_m$-gerbe on the rigid space $T^*X\{-1\}_{<1}$.

\end{enumerate}

In terms of the above descriptions, after restricting to perfect complexes and inverting $p$, the equivalence \eqref{eq:twistedSimpsonFS} becomes
\begin{equation}
    \label{eq:twistedSimpsonSmall}
    \tag{$\mathrm{SimpRig}_{\sm}$}
\Perf(X_{\proet},\widehat{\calO}_X)_\sm \simeq \mathcal{D}_{X\text{-}\Perf}(T^*X\{-1\}_{<1})_{\mathcal{P}_{\mathcal{X},\eta}}
\end{equation}
exchanging $\otimes$ on the LHS with convolution on the RHS. This statement is often called the ``small'' $p$-adic Simpson equivalence for rigid spaces. In the next section, we explain how to extend past the small case.

\begin{remark}[Period sheaves]
Another approach to \eqref{eq:twistedSimpsonSmall} is via period sheaves (see \cite{LiuZhuRH, WangSimpsonRigid, ShengWangSemistable}). This approach is naturally explained by the stack-theoretic perspective: the relevant period sheaves arise from the splittings of the Hodge--Tate gerbe (or, rather, its pushout along $BT_{\mathcal{X}}^\sharp\{1\} \to BT_{\mathcal{X}}\{1\}$; this object is discussed in a broader context in Remark~\ref{rmk:HTlocusNaivePrism}), pulled back to $X_\proet$ as in item (\ref{item:SmallCrystalSide}) above; see \cite[\S 1.6, \S 7.5]{ALBH1}. We thank Ben Heuer for this remark.
\end{remark}

\subsection{Beyond small objects}

To even formulate a statement extending \eqref{eq:twistedSimpsonSmall} beyond small objects, we need to extend the Simpson gerbe in two ways: it should live on the entire cotangent bundle (not just the open unit disc bundle) and it should make sense for all smooth rigid spaces (not just those admitting a smooth formal model). Both obstacles can be overcome by, roughly, taking the generic fibre of the formal scheme Simpson gerbe theory:

\begin{proposition}[The Simpson gerbe for rigid spaces]
    \label{prop:SimpsonGerbeRigid}
To any smooth rigid space $X/C$, one can functorially attach a multiplicative $\mathbf{G}_m$-gerbe $\mathcal{P}_X \to T^*X\{-1\}$, compatibly with \'etale descent, and uniquely characterized by the following property:  when $X=\mathcal{X}_\eta$ for a smooth scheme $\mathcal{X}/\Spf(\calO_C)$, there is a natural identification $\mathcal{P}_X|_{T^*X\{-1\}_{<1}} \simeq \mathcal{P}_{\mathcal{X},\eta}$.  
\end{proposition}

Let us make the functoriality in the first sentence precise. Consider the category $\mathrm{Rig}_C^{\sm}$ of smooth rigid spaces over $C$, endowed with the \'etale topology. The cotangent module functor $\Omega^1_{-/C}$ is an \'etale sheaf on this site. The functoriality simply means that $X \mapsto \mathcal{P}_X$ lifts to map $\Omega^1_{-/C}\{-1\} \to B^2\mathbf{G}_m$ in $\mathcal{D}(\Rig_C^\sm)$.

\begin{proof}[Proof sketch]
Let $\mathrm{Aff}^{\mathrm{fp}}_{\Spf(\calO_C)}$ be the category of affine topologically finitely presented schemes over $\Spf(\calO_C)$, endowed with the \'etale topology; write $\mathrm{Aff}^{\sm}_{\Spf(\calO_C)} \subset \mathrm{Aff}^{\mathrm{fp}}_{\Spf(\calO_C)}$ for the full subcategory of smooth objects. Let $\widehat{L}_{-/\calO_C} \in \mathcal{D}(\mathrm{Aff}^{\mathrm{fp}}_{\Spf(\calO_C)})$ be the \'etale sheaf given by the $p$-completed algebraic cotangent complex (see \cite[Definition 7.1.23]{GabberRamero}); on $\mathrm{Aff}^{\sm}_{\Spf(\calO_C)}$, this is concentrated in degree $0$ and denoted $\Omega^1_{-/\calO_C}$.

The Simpson gerbe for formal schemes (Theorem~\ref{thm:HTstackformal}) defines a natural map\footnote{The factor $\mathfrak{m}_C$ appears because the gerbe in Theorem~\ref{thm:HTstackformal} lives on the formal completion $\widehat{T^*\mathcal{X}}\{-1\}$, so it can only be evaluated on topologically nilpotent $1$-forms; for $\mathcal{X}$ smooth over $\Spf(\calO_C)$, these are exactly the sections of $\mathfrak{m}_C \cdot \Omega^1_{\mathcal{X}/\calO_C}\{-1\}$. This tensor factor can be avoided using the theory of sheared Hodge--Tate stacks (see Remark~\ref{rmk:ShearedVariant} and \cite[Footnote 63]{Bhatt2025Apsects}), which gives such a gerbe directly on the full cotangent bundle.}
\[ \mathcal{P}^{\sm}_{\fSch}: \left(\mathfrak{m}_C \otimes_{\calO_C} \widehat{L}_{-/\calO_C}\{-1\}\right)|_{\mathrm{Aff}^{\sm}_{\Spf(\calO_C)}} \to B^2 \mathbf{G}_m \]
in $\mathcal{D}(\mathrm{Aff}^{\sm}_{\Spf(\calO_C)})$. Moreover, the existence of this map encodes both the functoriality and the \'etale descent properties of the Simpson gerbe for smooth formal schemes. We now use the following fact (whose proof is omitted here): for any $\Spf(S) \in \mathrm{Aff}^{\mathrm{fp}}_{\Spf(\calO_C)}$, the natural map
\[ \colim_{R \to S} \widehat{L}_{R/\calO_C} \to \widehat{L}_{S/\calO_C} \]
is an equivalence in $\mathcal{D}(\calO_C)$ (without $p$-completions!), where the colimit runs over all $\Spf(R) \in \mathrm{Aff}^{\sm}_{\Spf(\calO_C)}$ equipped with a map $R \to S$. In other words, $\widehat{L}_{-/\calO_C}$ is left Kan extended from its restriction to $\mathrm{Aff}^{\sm}_{\Spf(\calO_C)}$, and hence so is $\mathfrak{m}_C \otimes_{\calO_C} \widehat{L}_{-/\calO_C}\{-1\}$. Then $\mathcal{P}^{\sm}_{\fSch}$ extends uniquely to a map
\[ \mathcal{P}_{\fSch}: \mathfrak{m}_C \otimes_{\calO_C} \widehat{L}_{-/\calO_C}\{-1\} \to B^2 \mathbf{G}_m \]
in $\mathcal{D}(\mathrm{Aff}^{\mathrm{fp}}_{\Spf(\calO_C)})$, giving a theory of the Simpson gerbe for possibly singular formal schemes. 

Let $\Rig_C$ be the category of all rigid spaces over $C$, endowed with the \'etale topology. Taking generic fibres gives a morphism of topoi
\[ \nu: \Shv_\et(\Rig_C) \to \Shv_\et(\mathrm{Aff}^{\mathrm{fp}}_{\Spf(\calO_C)}),\]
characterized by $(\nu_* F)(\Spf(R)) = F(\Spf(R)_\eta)$. Write $L^{\mathrm{an}}_{-/C} \in \mathcal{D}(\Rig_C)$ for the analytic cotangent complex of rigid spaces (see \cite[Definition 7.2.31]{GabberRamero}), regarded as an \'etale sheaf. Thus $L^{\mathrm{an}}_{\Spf(R)_\eta/C} \simeq \widehat{L}_{R/\calO_C}[1/p]$ for $\Spf(R) \in \mathrm{Aff}^{\mathrm{fp}}_{\Spf(\calO_C)}$, and $L^{\mathrm{an}}_{U/C} \simeq \Omega^1_{U/C}$ for $U$ smooth. We obtain the following diagram in $\mathcal{D}(\Rig_C)$:
\[ \xymatrix{
\nu^{-1} \left(\mathfrak{m}_C \otimes_{\calO_C} \widehat{L}_{-/\calO_C}\{-1\}\right) \ar[rr]^-{\nu^{-1} \mathcal{P}_{\fSch}} \ar[d]_-{\text{can}_L} && \nu^{-1} B^2 \mathbf{G}_m \ar[d]^-{\text{can}_{\mathbf{G}_m}} \\
L^{\mathrm{an}}_{-/C}\{-1\} && B^2 \mathbf{G}_m,
} \]
where $\nu^{-1}$ denotes the derived pullback on sheaves of complexes along $\nu$, and the vertical maps are the natural ones induced by $\nu$. We claim:
\begin{itemize}
\item[$(\ast)$] The left vertical map $\text{can}_L$ is an isomorphism (so, in particular, $p$ acts invertibly on its source).
\end{itemize}

Granting $(\ast)$, there is a unique way to complete the above diagram to a commutative square by adding a map $\mathcal{P}_{\Rig}: L^{\mathrm{an}}_{-/C}\{-1\} \to B^2 \mathbf{G}_m$ in the bottom horizontal line. Restricting $\mathcal{P}_{\Rig}$ to smooth rigid spaces (where $L^{\mathrm{an}}_{-/C}\{-1\} \simeq \Omega^1_{-/C}\{-1\}$) gives the desired Simpson gerbe construction together with its properties from the proposition. For the uniqueness assertion, $(\ast)$ shows that $L^{\mathrm{an}}_{-/C}\{-1\}$ is the sheafified left Kan extension of its restriction to smooth rigid spaces; thus $\mathcal{P}_{\Rig}$ is determined by its restriction to smooth rigid spaces, and is characterized there by agreeing with $\mathcal{P}_{\fSch}$ on generic fibres of smooth formal schemes.

It remains to show $\text{can}_L$ is an isomorphism. The statement is clear after inverting $p$ by the definition of the target \cite[\S 7.2.29--7.2.31]{GabberRamero}. As $p$ acts invertibly on the target,  we must show that $p$ acts invertibly on the source, or equivalently that the source vanishes after mod $p$ reduction. For this, it is enough to show that $\nu^{-1} \widehat{L}_{-/\calO_C}/p \simeq 0$. This can be checked on stalks, which lets us ignore the sheafification implicit in computing $\nu^{-1}$. A conservative system of points of $\Shv_\et(\Rig_C)$ is given by geometric points $\overline{x} \to X$ of rigid spaces $X$, i.e., by the pro-objects $\{U_i = \mathrm{Spa}(S_i,S_i^+)\}_i$ formed by all affinoid \'etale neighbourhoods of $\overline{x}$ in $X$. As maps $U_i \to \Spf(R)_\eta$ correspond to maps $R \to S_i^+$, we must show that
\[ \colim L_{(R/p)/(\calO_C/p)} \simeq 0, \]
where the colimit runs over all $\Spf(R) \in \mathrm{Aff}^{\mathrm{fp}}_{\Spf(\calO_C)}$ with a map $R \to \{S_i^+\}$ to the displayed ind-object. But the formation of the cotangent complex commutes with colimits and, by writing each $S_i^+$ as a colimit of its open bounded subrings, we  have 
\[ \colim R/p = \colim_i S_i^+/p = k(\overline{x})^+/p,\] where $k(\overline{x})^+$ is the perfectoid valuation ring attached to the geometric point $\overline{x}$. The desired vanishing then follows from the vanishing of the $p$-completed cotangent complex for any map of perfectoid rings. 
\end{proof}

Having extended the Simpson gerbe to the entire cotangent bundle of the rigid space, there is  an obvious candidate for the non-small version of \eqref{eq:twistedSimpsonSmall}:

\begin{theorem}[Twisted $p$-adic Simpson equivalence]
\label{thm:padicSimpson}
Fix a smooth rigid space $X/C$. Then there is a $\Perf(X)$-linear equivalence
\begin{equation}
    \label{eq:twistedSimpson}
    \tag{SimpRig}
\Perf(X_{\proet},\widehat{\calO}_X) \simeq \mathcal{D}_{X\text{-}\text{Perf}}(T^*X\{-1\})_{\mathcal{P}_X}
\end{equation}
exchanging $\otimes$ on the LHS with convolution on the RHS. 
\end{theorem}
\begin{proof}[Idea of proof, following {\cite[Theorem 12.4.11]{Bhatt2025Apsects}}]
Consider the category $\text{Aff}\Rig_C^{\sm,ss} \subset \Rig_C^{\sm}$ consisting of smooth affinoid rigid spaces over $C$ whose canonical formal model is semistable. As every smooth rigid space over $C$ admits, \'etale locally, a semistable formal model over $\Spf(\calO_C)$ \cite[Corollary 3.3.2]{TemkinAlteredLocal}, the category $\text{Aff}\Rig_C^{\sm,ss}$ forms a basis for $\Rig_C^\sm$ in the \'etale topology. For a semistable scheme $\mathcal{U}/\Spf(\calO_C)$ with generic fibre $U = \mathcal{U}_\eta$, the natural logarithmic version of Theorem~\ref{thm:HTstackformal} (3) yields the following on taking generic fibres: 
\begin{equation}
    \label{eq:twistedSimpsonFSetalog}
    \tag{$\mathrm{SimpRig}^{\text{ss}}_{\sm}$}
\Perf(U_{\proet},\widehat{\calO}_U)_\sm \simeq \mathcal{D}_{U\text{-}\Perf}(T^*U\{-1\}_{<1})_{\mathcal{P}_{\mathcal{U},\eta}}
\end{equation}
where the notion of smallness on both sides are defined using the canonical formal model $\mathcal{U}/\Spf(\calO_C)$ (with its natural log structure). As $U$ varies, this gives an equivalence of presheaves on $\text{Aff}\Rig_C^{\sm,ss}$. We claim that the global sections of the corresponding equivalence of sheaves is the desired statement \eqref{eq:twistedSimpson}. This amounts to identifying the global sections of the sheafification of the LHS and RHS of \eqref{eq:twistedSimpsonFSetalog} with the LHS and RHS of \eqref{eq:twistedSimpson} separately. For the RHS, this amounts to the standard fact (already used implicitly in Proposition~\ref{prop:SimpsonGerbeRigid}) that pulling back along rig-\'etale covers of formal schemes can make differential forms $p$-adically small. For the LHS, this amounts to showing that object of $\Perf(U_\proet,\widehat{\calO}_U)$ comes, \'etale locally on $U$, from a Hodge--Tate crystal on $\mathcal{U}$; this is deduced from the isogeny theorem for prismatic cohomology by a deformation theory argument due to Faltings, see \cite[Corollary 12.3.10]{Bhatt2025Apsects}.
\end{proof}

\begin{remark}[The canonical Higgs field]
\label{rmk:CanonicalHiggs}
Theorem~\ref{thm:padicSimpson} implies that any $E \in \Perf(X_\proet,\widehat{\calO}_X)$ comes equipped with a $\mathrm{Sym}^*(T_X\{1\})$-module structure (i.e., a Higgs field in the category $\Perf(X_\proet,\widehat{\calO}_X)$). Given its construction, this ``canonical Higgs field'' can be regarded as a rigid analog of the $p$-curvature for $\mathcal{D}$-modules in characteristic $p$ (Remark~\ref{rmk:pCurvature}). 

The canonical Higgs field was first constructed by Rodr\'iguez Camargo \cite[Theorem 1.0.3]{CamargoGeometricSen}, building on work of Pan \cite[\S 3.1]{PanLocallyAnalytic} (see also \cite[Theorem 4.8]{HeuerSimpson}). Their construction relies on a ``decompletion'' argument using local toric charts and explicit calculations with the associated towers of rigid spaces, which makes the Higgs field quite concrete. The approach above is more abstract but needs no local co-ordinates or explicit computations: the Higgs field on Hodge--Tate crystals (and thus in the small correspondence) comes from twisted Cartier duality, while the passage to arbitrary pro-\'etale perfect complexes uses the isogeny theorem for prismatic cohomology\footnote{The original proof of the isogeny theorem in \cite{BhattScholzePrism} {\em did} involve explicit local computations. However, the Nygaard filtered prismatization \cite{BhattFGaugeNotes} gives a conceptual construction of the Nygaard filtration on prismatic cohomology and its basic properties, from which the isogeny theorem follows immediately.} (recasting Faltings' original ideas \cite{FaltingsSimpson} in more modern language). This method also applies to fairly general coefficients, see \cite[\S 12]{Bhatt2025Apsects}.
\end{remark}

\begin{remark}[The analytic Hodge--Tate stack perspective]
\label{rmk:AnalyticHTstack}
Grosso modo, the proof of \eqref{eq:twistedSimpson} sketched above had two steps: first prove the corresponding statement \eqref{eq:twistedSimpsonFS} for formal schemes using the torsor structure on the Hodge--Tate stack and  Cartier duality, and then pass to generic fibres. In recent work, Ansch\"utz, Le Bras, Rodr\'iguez Camargo, and Scholze \cite{ALBRCSRatAnalPrism} have developed a theory of the {\em analytic} Hodge--Tate stack (and more generally rational {\em analytic} prismatization) directly for rigid spaces; this enables them to directly work with rigid spaces to prove \eqref{eq:twistedSimpson} via an enhancement of Cartier duality \cite{CamargoCartierDual} made possible by condensed mathematics. We shall return to this topic when quantizing \eqref{eq:twistedSimpson} in \S \ref{ss:RigAnQuanNAHT}.
\end{remark}

\subsection{Trivialize the gerbe on all spectral varieties}

To obtain more explicit statements from Theorem~\ref{thm:padicSimpson}, following the characteristic $p$ discussion (reviewed in \S \ref{ss:NAHTcharp}), one must find loci in $T^*X\{-1\}$ where the Simpson gerbe trivializes. The cleanest instance is following:

\begin{proposition}
\label{prop:LiftsTrivializeGerbe}
Fix a smooth rigid space $X/C$, and write $B_2(C) := B_\dR^+(C)/\ker(\theta)^2$. There is an equivalence of groupoids
\[
\left\{ \parbox{3.6cm}{\centering Liftings of $X$ to $B_2(C)$} \right\}
\simeq
\left\{ \parbox{6cm}{\centering Trivializations of the multiplicative $\mathbf{G}_m$-gerbe $\mathcal{P}_X$  on $\widehat{T^*X}\{-1\}$} \right\}.
\]
\end{proposition}

The preceding leads to an explicit nilpotent Simpson correspondence, see \cite[Corollary 12.5.5]{Bhatt2025Apsects}, and also \cite{SoloveraNeryUnipotent}. Beyond the nilpotent setting, the original work of Faltings \cite{FaltingsSimpson} and Heuer \cite{HeuerSimpson} can be recast as the following remarkable statement:

\begin{theorem}[Faltings, Heuer]
    \label{thm:padicSimpsonUntwisted}
Assume $C = \overline{C}$. Fix a smooth rigid space $X/C$, as well as the following:
\begin{enumerate}
    \item A lift of $X$ to $B_2(C)$.
    \item A continuous\footnote{Continuity is automatic if we demand that $\widetilde{\exp}$ agree with the usual exponential on $p^n \calO_C$ for some $n \gg 0$.} group homomorphism $\widetilde{\exp}:C \to 1+\mathfrak{m}_C$ splitting the logarithm $\log:1+\mathfrak{m}_C \to C$.
\end{enumerate}
For any proper $Z \subset T^*X\{-1\}$, the gerbe $\mathcal{P}_X|_Z$ is naturally (in the above choices) trivial. In particular, when $X$ itself is proper, we obtain a natural (in the above choices) equivalence
\[ \Perf(X_{\proet},\widehat{\mathcal{O}}_X) \simeq \mathcal{D}_{X\text{-}\Perf}(T^*X\{-1\})\]
exchanging $\otimes$ on the LHS with convolution on the RHS.
\end{theorem}

Heuer's work \cite{HeuerSimpson} essentially contains this statement. However, it predates the Simpson gerbe, so does not prove this statement as written. One proof appears in \cite[\S 12]{Bhatt2025Apsects}. We give instead a simpler argument due to Petrov, which replaces $p$-adic Hodge theory for singular rigid spaces used in {\em loc.\ cit.}\ with the classical primitive comparison theorem. This argument also avoids the choices of base points needed in \cite{HeuerSimpson}. (We thank Petrov for allowing us to include it.)

\begin{proof}[Proof (Petrov)]
The second part follows from the first by trivializing $\mathcal{P}_X$ on the union of all spectral varieties (i.e., on analytic closed subspaces of $T^*X\{-1\}$ finite over $X$) and observing that any object of $\mathcal{D}_{X\text{-}\Perf}(T^*X\{-1\})_{\mathcal{P}_X}$ is supported on such a $Z$. For the first, we shall use a different construction of the Simpson gerbe that we recall first. 

Let $T$ be any smooth rigid space equipped with a lift to $B_2(C)$, and let $\nu:T_\proet \to T_\et$ be the natural map. Then there is a natural (in the lift) algebra isomorphism
\[ R\nu_* \widehat{\calO}_T \simeq \bigoplus_i \Omega^i_{T/C}\{-i\}[-i],\]
see, e.g., \cite[Proposition 7.2.5]{GuoHTsingular}. In particular, there is a natural map 
\[ c:\Omega^1_{T/C}\{-1\} \to R\nu_*\widehat{\calO}_T[1]\] 
in $\mathcal{D}(T_\et)$. On the other hand, there is a natural short exact sequence
\[ 0 \to \mu_{p^\infty} \to 1+\mathfrak{m}_C\widehat{\calO}^+_T \xrightarrow{\log} \widehat{\calO}_T \to 0\]
on $T_\proet$. Writing $\partial_T$ for the boundary map and pushing forward to $T_\et$ then allows us to define a composition
\[ \Omega^1_{T/C}\{-1\} \xrightarrow{c} R\nu_*\widehat{\calO}_T[1] \xrightarrow{\partial_T} \mu_{p^\infty}[2] \xrightarrow{\text{can}} \mathbf{G}_m[2],\]
where the third term uses the identification $\mu_{p^\infty} \simeq R\nu_* \mu_{p^\infty}$ on $T_\et$, and the last map comes from the inclusion $\mu_{p^\infty} \to \mathbf{G}_m$. Taking global sections gives
\[ R\Gamma(T_\et,\Omega^1_{T/C}\{-1\}) \xrightarrow{R\Gamma(c)} R\Gamma(T_\proet,\widehat{\calO}_T)[1] \xrightarrow{R\Gamma(\partial_T)} R\Gamma(T_\et,\mu_{p^\infty})[2] \xrightarrow{\text{can}} R\Gamma(T_\et,\mathbf{G}_m)[2].\]
We shall use without proof the following compatibility\footnote{This alternate description is a characteristic $0$ version of the description of the Brauer class of the Azumaya algebra of differential operators in characteristic $p$ via the Milne sequence and the tautological form on the cotangent bundle, see \cite[Proposition 4.4]{OgusVologodsky} or \cite[Remark 2.2.4 (4)]{BMRLocalization}.}:

\begin{itemize}
\item[$(\ast)$] Assume $T=T^*X\{-1\}$, equipped with the $B_2(C)$-lift coming from the chosen lift of $X$. Then evaluating the above map on the tautological section $\omega_\tau$ of $\Omega^1_{T/C}\{-1\}$ (i.e., precomposing with $\mathbf{Z} \xrightarrow{\omega_\tau} R\Gamma(T_\et,\Omega^1_{T/C}\{-1\})$) gives the point of $R\Gamma(T_\et,\mathbf{G}_m[2])$ representing the Simpson gerbe on $T$.
\end{itemize}

Thanks to the previous description and functoriality, in order to prove the theorem, it suffices to show that for any proper $Z \subset T=T^*X\{-1\}$ the map 
\[ R\Gamma(\partial_Z):R\Gamma(Z_\proet,\widehat{\calO}_Z) \to R\Gamma(Z_\et,\mu_{p^\infty})[1] \]
is trivialized naturally in $Z$. But $Z$ is proper, so the primitive comparison theorem \cite[Theorem 3.17]{ScholzeCDM} identifies this with the map 
\[ R\Gamma(Z_\proet,C) \to R\Gamma(Z_\proet,\mu_{p^\infty})[1] \]
coming from the boundary map for the exponential sequence
\[ 0 \to \mu_{p^\infty} \to 1+\mathfrak{m}_C \xrightarrow{\log} C \to 0,\]
thought of as a sequence of abelian sheaves on $Z_\proet$. The choice of the exponential then splits this sequence, and thus trivializes $R\Gamma(\partial_Z)$ naturally in $Z$. 
\end{proof}

\newpage
\section{Quantized non-abelian $p$-adic Hodge theory}
\label{sec:QuantNAHT}

\centerline{\fbox{The work described in this section is joint work in progress \cite{BHPVQuanNAHT}.}}

\vspace{0.1in}

The Simpson equivalences of \S \ref{ss:NAHTcharp} and \S \ref{sec:Simpson} have the same shape: a category of ``de Rham'' or ``topological'' coefficients is identified, up to an explicit gerbe twist, with Higgs modules. In this section, we explain why these equivalences admit natural deformations.  On the Higgs side, the deformation replaces Higgs modules with $\lambda$-connections for a suitable $\lambda$, so we refer to such deformed equivalences as ``quantizations.''\footnote{Passing from Higgs fields to $\lambda$-connections deforms the commutative algebra of functions on the cotangent bundle into a non-commutative algebra of differential operators; the corresponding passage in physics, with $\lambda$ playing the role of Planck's constant $\hbar$, is called {\em quantization}.} On the other side, the deformations come from the prismatic or $B_\dR^+$-directions.

As explained below, the stack-theoretic perspective is critical to this work: the quantized Simpson equivalences (Theorems~\ref{thm:QuantizationPhiCrys} and \ref{thm:QuantNAHTrigid}) arise via Cartier duality from an identity of torsors for a commutative group stack over a common base (Corollary~\ref{cor:TorsorIdentity}), which in turn is the transmutation of an identity
\[ [\Ga^\Prism] = [\Ga^{I,\crys}] + [\Ga^\calF] \]
in the group of square-zero extensions of a fixed ring stack by a fixed group stack (Proposition~\ref{prop:TorsorAdditionIdentity}).

\subsection{Introduction}
\label{ss:QuantIntro}
In this extended introductory section, we present the main objects being considered, the primary quantization questions motivating us, and the general shape of our answers. 

\subsubsection{$\lambda$-connections}
\label{sss:lambdaconn}
We begin by recalling Deligne's notion of a $\lambda$-connection \cite[\S 4]{SimpsonLambda}, its role in non-abelian Hodge theory over $\mathbf{C}$, and a stack-theoretic perspective on it that will be important in the sequel.

\begin{notation}
\label{not:lambdaconn}
Fix a commutative ring $A$ equipped with an invertible ideal $I$. Write $\overline{A}=A/I$, and $M\{1\} = M \otimes_{\overline{A}} I/I^2$ for the displayed twist of an $\overline{A}$-module $M$ (and similarly $M\{n\}$ for any $n \in \mathbf{Z}$). Assume $A/I$ has characteristic $0$.  Finally, fix a smooth scheme $\mathfrak{X}/\Spec(A)$, and set $X=\mathfrak{X}_{\overline{A}}$. (The definitions below readily transplant to smooth schemes over $\Spf(A)$ when $A$ is given the $I$-adic topology.)
\end{notation}

The principal definition is the following:

\begin{definition}[$I$-connections]
An {\em $I$-connection on $\mathfrak{X}/A$} consists of a vector bundle $E \in \mathrm{Vect}(\mathfrak{X})$ and a map $\nabla:E \to E \otimes \frac{1}{I}\Omega^1_{\mathfrak{X}/A}$ that obeys the Leibniz rule and satisfies $\nabla^2=0$\footnote{When $I=(\lambda)$, by trivializing the fractional ideal $\frac{1}{I}$ via $1/\lambda \in \frac{1}{I}$,  the map $\nabla$ becomes a map $\nabla:E \to E \otimes \Omega^1_{\mathfrak{X}/A}$ satisfying the identity $\nabla(fm) = f \nabla(m) + \lambda m \otimes df$ and $\nabla^2=0$; this is the standard presentation of this notion.}. (In particular, we shall only work with {\em flat} connections.)
\end{definition}

\begin{remark}[The $\mathcal{D}$-module perspective]
Say $I = (\lambda)$. Let $\mathcal{D}_{\mathfrak{X}/A,\lambda} \subset \mathcal{D}_{\mathfrak{X}/A}$ be the subalgebra of the ring of differential operators generated by $\mathcal{O}_{\mathfrak{X}}$ and $\lambda T_{\mathfrak{X}/A}$. Given local co-ordinates $x_1,...,x_n$ on $\mathfrak{X}$, it is generated by the operators $\{x_i\}_{i=1,...,n}$ and $\{\lambda \frac{\partial}{\partial x_i}\}_{i=1,...,n}$, subject to the obvious commutation relations (see \cite[p.~87]{SimpsonModuliI}). An $I$-connection on $\mathfrak{X}/A$ is then the same as a vector bundle on $\mathfrak{X}$ whose $\mathcal{O}_{\mathfrak{X}}$-module structure has been extended to a $\mathcal{D}_{\mathfrak{X}/A,\lambda}$-module structure.
\end{remark}

Thus, if $I=A$, an $I$-connection is just a flat connection. On the other hand, reduction modulo $I$ defines a functor from $I$-connections on $\mathfrak{X}/A$ towards Higgs bundles on $X/\overline{A}$ (where the Higgs fields are valued in the Tate-twist $\Omega^1_{X/\overline{A}}\{-1\}$). Thus, we can regard the notion of $I$-connections on $\mathfrak{X}/A$ as an $I$-adic deformation of the notion of Higgs bundles on the ``special fibre'' $X/\overline{A}$ towards flat connections on the ``generic fibre'' $\mathfrak{X}_{\Spec(A[1/I])}$. The motivating example is the following:

\begin{example}[The $\lambda$-connection attached to a filtered flat connection]
Let $X$ be a smooth $\mathbf{C}$-scheme. Let $(E,\nabla,F^\bullet E)$ be a vector bundle with flat connection $(E,\nabla)$ endowed with a  Griffiths transversal filtration $F^\bullet E$ by subbundles; Gauss--Manin connections provide natural examples of such triples. Let $A=\mathbf{C}[\lambda]$, let $I=(\lambda)$, and let $\mathfrak{X} = X \times \Spec(A)$. Then the Rees module
\[ \mathcal{R}(E,\nabla,F^\bullet E) := \bigoplus_i F^i E \cdot \lambda^{-i} \subset E \otimes_{\mathbf{C}} \mathbf{C}[\lambda^{\pm 1}], \]
endowed with the operator induced by the connection, is a $\lambda$-connection on $\mathfrak{X}/A$. It interpolates between the flat connection $(E,\nabla)$ (at $\lambda=1$) and the Higgs bundle $(\mathrm{gr}^\bullet_F E,\Theta := \overline{\nabla})$ (at $\lambda=0$).
\end{example}

Given the last example, it is unsurprising that $\lambda$-connections appear quite fundamentally in non-abelian Hodge theory over $\mathbf{C}$. Let us summarize some of the main results in this direction (working in the analytic category, as that is most natural to consider here):

\begin{theorem}[Corlette {\cite{CorletteFlat}}, Simpson {\cite{SimpsonHiggs}, \cite[\S 4]{SimpsonLambda}}]
Let $X$ be a compact K\"ahler manifold, and let $\mathbf{A}^1_\lambda$ denote the complex $\lambda$-line, regarded as a complex manifold. There is a commutative diagram of symmetric monoidal functors
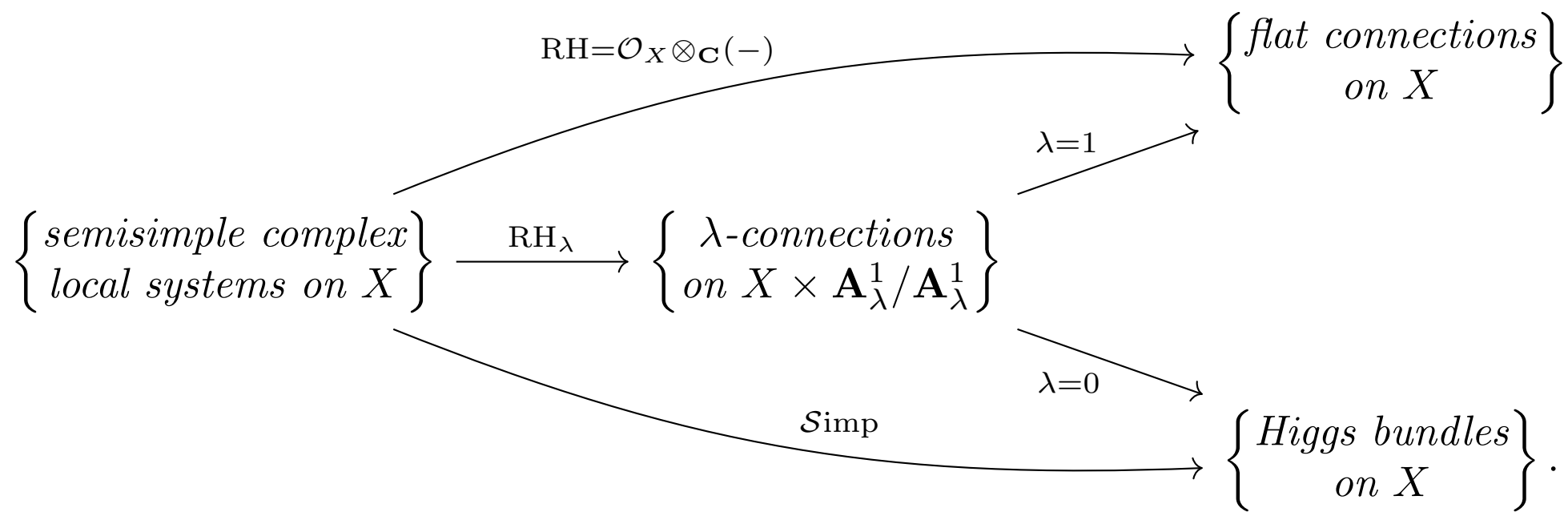
\[
\begin{tikzcd}[column sep=4em, row sep=1.5em]
& & \left\{ \begin{tabular}{@{}c@{}} flat connections \\ on $X$ \end{tabular} \right\} \\
\left\{ \begin{tabular}{@{}c@{}} semisimple complex \\ local systems on $X$ \end{tabular} \right\} \arrow[r, "\mathrm{RH}_\lambda"]
\arrow[urr, bend left=12, "\mathrm{RH} = \mathcal{O}_X \otimes_{\mathbf{C}} (-)"]
\arrow[drr, bend right=12, "\mathcal{S}\mathrm{imp}"]
& \left\{ \begin{tabular}{@{}c@{}} $\lambda$-connections \\ on $X \times \mathbf{A}^1_\lambda/\mathbf{A}^1_\lambda$ \end{tabular} \right\} \arrow[ur, "\lambda=1"] \arrow[dr, "\lambda=0"']
& \\
& & \left\{ \begin{tabular}{@{}c@{}} Higgs bundles \\ on $X$ \end{tabular} \right\}.
\end{tikzcd}
\]
The functor $\mathrm{RH}$ is the restriction of the classical Riemann--Hilbert equivalence, and is thus an equivalence onto semisimple flat connections, while the functor $\mathcal{S}\mathrm{imp}$ is an equivalence onto polystable Higgs bundles with vanishing Chern classes. Moreover, all functors are compatible with cohomology in a natural sense. 
\end{theorem}

The ultimate goal of the material presented here is to give an analog of the functor $\mathrm{RH}_\lambda$, after $\lambda$-adic completion, in the setting of $p$-adic analytic geometry; this will be discussed further in \S \ref{ss:BdRQuantization}.  In the rest of this section, we explain a stack-theoretic perspective on $I$-connections that seems potent near $\{I=0\}$. 

\begin{assumption}
\label{ass:IComplete}
For the rest of \S \ref{sss:lambdaconn}, assume $A$ is $I$-adically complete. Endow $A$ with the $I$-adic topology to define $\Spf(A)$, and repurpose notation: $\mathfrak{X}$ now denotes a smooth scheme over $\Spf(A)$. Write  $X=\mathfrak{X}_{\overline{A}}$ as before, and write $\widetilde{X} = \mathfrak{X}_{A/I^2}$ for the mod $I^2$ lift.
\end{assumption}

We shall proceed by transmutation. Since $I \subset \widehat{\Ga}(\Spf(A))$, we have a map
\[ \tau:\Ga^{I,\crys} := \Ga/I\widehat{\Ga} \to \Ga^\Weil := \Ga/I \]
of animated $A/I^2$-algebra stacks on $\CAlg^\an_{\Spf(A)}$, which is a square-zero extension as $I^2\Ga \subset I\widehat{\Ga}$. Transmuting gives a map
\[ \pi_{I,\crys}:(\widetilde{X}/A)^{I,\crys} \to (X/A)^\Weil,\]
whose target is the (derived) Weil restriction of $X$ along $\Spf(\overline{A}) \hookrightarrow \Spf(A)$ (Definition~\ref{def:WeilStack}). The source, which we call the {\em $I$-crystallization} of $\widetilde{X}$, controls $I$-connections, much as the de Rham stack controls flat connections:

\begin{proposition}[$I$-connections depend only on the mod $I^2$ lift]
\label{prop:IConnModI2}
The category of $I$-connections on $\mathfrak{X}/A$ is naturally identified with $\mathrm{Vect}( (\widetilde{X}/A)^{I,\crys})$. In particular, this category depends only on $\widetilde{X}$.
\end{proposition}

The last assertion above is the $I$-connections variant of the standard fact that the category of flat connections on $\mathfrak{X}/A$ only depends on $X$. The square-zero property of $\tau$ mentioned above transmutes to the following structure, whose prismatic variant (Corollary~\ref{cor:GerbeIdR}) will be critical to the sequel.

\begin{proposition}[$I$-connections as a torsor]
\label{prop:IConnTorsor}
The map $\pi_{I,\crys}$ is naturally a torsor for a commutative group stack (which is, up to some twists, the de Rham stack of $T_X$ relative to $X$).
\end{proposition}

Over $\{I=0\}$, the map $\pi_{I,\crys}$ becomes (up to twists) $B\widehat{T}_X \to BT_X$, whose pullback is the inclusion of nilpotent Higgs bundles into all Higgs bundles.

\subsubsection{Quantization over a ($\varphi$-)crystalline prism}
\label{ss:QuantizePhiCrys}

In this subsection, we discuss the quantization problem for the twisted Simpson equivalence over certain prisms. As the story is already interesting in the crystalline case, let us first formulate the question there:

\begin{question}
    \label{q:quantizecharp}
For a smooth variety $X$ over a perfect field $k$ of characteristic $p$, Corollary~\ref{cor:SimpsonNilpCharp} yields an identification
\[
\begin{tikzcd}[column sep=6em]
\left\{
\parbox{5cm}{\centering nilpotent $\mathcal{D}$-modules on $X/k$}
\right\}
\arrow[r, dashed, leftrightarrow, "\text{twist by } \mathcal{P}_X"]
&
\left\{
\parbox{6cm}{\centering nilpotent Higgs modules on $X^{(1)}/k$}
\right\}.
\end{tikzcd}
\]
Here $\mathcal{P}_X$ denotes the Simpson gerbe on $\widehat{T^*X^{(1)}}$. It is relatively clear how to twist the RHS with a gerbe on $\widehat{T^*X^{(1)}}$; to twist the LHS, one uses $p$-curvature (Remark~\ref{rmk:pCurvature}). Given a lift $\widetilde{X}_\infty/\Spf(W(k))$ of $X/k$, both sides above admit natural $W(k)$-deformations (that in fact depend only on the mod $p^2$ lift $\widetilde{X}$) given by
\[
\begin{tikzcd}[row sep=1.5em]
\left\{ \parbox{3.2cm}{\centering crystals on $X/W(k)$} \right\}
\arrow[d, phantom, sloped, "\simeq"]
\\
\left\{ \parbox{3.6cm}{\centering nilpotent $\mathcal{D}$-modules on $\widetilde{X}_\infty/W(k)$} \right\}
\end{tikzcd}
\quad \text{and} \quad
\left\{
\parbox{4.5cm}{\centering quasi-nilpotent $p$-connections on $\widetilde{X}_\infty/W(k)$}
\right\}
\]
respectively. This raises a natural question: is there a twisted $W(k)$-linear equivalence between these categories? We call such a statement a \emph{quantization}. The primary obstacle is the twist: the Simpson gerbe lives on $\widehat{T^*X^{(1)}}$, and it is unclear what it means to lift it to $W(k)$ or twist the categories appearing above by such a lift (e.g., there is no $p$-curvature for $\mathcal{D}$-modules on $\widetilde{X}_\infty/W(k)$). 
\end{question}

We shall formulate the answer in the more general setting of $\varphi$-crystalline prisms: the extra generality introduces no new complications, and the resulting ideas have direct bearing on our later work in \S \ref{ss:BdRQuantization} on quantization in rigid geometry. Thus, let $(A,I)$ be a $\varphi$-crystalline animated prism, i.e.,  $\varphi^* (A,I)$ is crystalline (Definition~\ref{def:phicrystallineprism}). Fix a smooth scheme $X/\Spf(\overline{A})$. Theorem~\ref{thm:HTstackformal} (3) gives an  equivalence
\[
\begin{tikzcd}[column sep=6em]
\left\{
\parbox{4.5cm}{\centering Hodge--Tate crystals on $X/A$}
\right\}
\arrow[r, dashed, leftrightarrow, "\text{twist by } \mathcal{P}_X"]
&
\left\{
\parbox{3.6cm}{\centering quasi-nilpotent Higgs modules on $X/\overline{A}$}
\right\}
\end{tikzcd}
\]
where the Higgs field takes values in $\Omega := \Omega^1_{X/\overline{A}}\{-1\}$ and $\mathcal{P}_X \to \widehat{\Omega}$ is the Simpson gerbe. To answer Question~\ref{q:quantizecharp}, we must deform the Simpson gerbe and the preceding twisted equivalence in the presence of lifting data. The ideal $I$ admits divided powers by $\varphi$-crystallinity, so there is a sensible notion of crystallizations $(Y/A)^\crys$ for geometric objects $Y/\Spf(\overline{A})$ (see \S \ref{ss:Crystallization}). Using this notion, we prove:

\begin{proposition}[The deformed Simpson gerbe]
    \label{prop:DeformedSimpsonGerbe}
Any lift $\widetilde{X}$ of $X$ to $A/I^2$ determines a crystal structure on the  $\mathbf{G}_m$-gerbe $\mathcal{P}_X \to \widehat{\Omega}$, i.e., a descent to a  $\mathbf{G}_m$-gerbe $\mathcal{P}_{\widetilde{X}}$ along the structure map $\widehat{\Omega} \to (\widehat{\Omega}/A)^\crys$ (lifting the Hodge--Tate version of this descent seen earlier in Remark~\ref{rmk:SimpsonGerbeConn}).
\end{proposition}

Using the above, our task is then to construct an equivalence of the form 
\[
\begin{tikzcd}[column sep=6em]
\left\{
\parbox{4.1cm}{\centering\normalfont prismatic crystals on $X/A$}
\right\}
\arrow[r, dashed, leftrightarrow, "\textnormal{twist by } \mathcal{P}_{\widetilde{X}}"]
&
\left\{
\parbox{5.1cm}{\centering\normalfont quasi-nilpotent $I$-connections on (any $A$-lift of) $\widetilde{X}/A$}
\right\}
\end{tikzcd}
\]
in the presence of a lift $\widetilde{X}$. Let us view the LHS as the derived category $\mathcal{D}((X/A)^\Prism)$ of the prismatization $(X/A)^\Prism$ (Proposition~\ref{prop:PrismaticCrystalsAsSheaves}). As indicated in \S \ref{sss:lambdaconn} and made precise in \S \ref{ss:IConnStack}, there is a stack $(\widetilde{X}/A)^{I,\crys}$ such that the RHS above is $\mathcal{D}( (\widetilde{X}/A)^{I,\crys})$. To relate these categories, we shall relate the stacks $(X/A)^\Prism$ and $(\widetilde{X}/A)^{I,\crys}$ themselves: we show that both these stacks are torsors over the same base for the same group (Corollary~\ref{cor:GerbePrismatic} and Corollary~\ref{cor:GerbeIdR}). Identifying their difference in terms of  (the dual to)  the gerbe above (Proposition~\ref{prop:TorsorAdditionIdentity}) then gives the desired statement on passage to derived categories:

\begin{theorem}[Quantization over $\varphi$-crystalline prisms]
    \label{thm:QuantizationPhiCrys}
Fix a lift $\widetilde{X}$ of $X$ to $A/I^2$ with induced gerbe $\mathcal{P}_{\widetilde{X}}$ (Proposition~\ref{prop:DeformedSimpsonGerbe}).
Then there is a natural geometrically defined action of $\mathcal{D}( (\widehat{\Omega}/A)^\crys)$ on both 
$\mathcal{D}((X/A)^\Prism)$ and $\mathcal{D}( (\widetilde{X}/A)^{I,\crys})$, as well as a natural identification of ${\mathcal{D}( (\widehat{\Omega}/A)^\crys)}$-modules
\[ \mathcal{D}( (X/A)^\Prism) \simeq
\mathcal{D}( (\widetilde{X}/A)^{I,\crys})_{\mathcal{P}_{\widetilde{X}}},\]
where the subscript $\mathcal{P}_{\widetilde{X}}$ indicates twisting by the invertible $\mathcal{D}( (\widehat{\Omega}/A)^\crys)$-module attached to $\mathcal{P}_{\widetilde{X}}$.
\end{theorem}

In other words, the notions of prismatic crystals and $I$-connections coincide up to an explicit gerbe twist. On the Hodge--Tate locus, this recovers the identification from Theorem~\ref{thm:HTstackformal} (3), with the additional feature that the Simpson gerbe $\mathcal{P}_X$ has now been endowed with a crystal structure.

\begin{remark}
The results discussed in this section involve certain nilpotence constraints; these can be avoided once the machinery of {\em sheared} prismatization is available, see Remark~\ref{rmk:ShearedVariant} and references therein.
\end{remark}

\begin{remark}
One may wonder about the necessity of the $\varphi$-crystallinity assumption on $(A,I)$ for the material discussed above: after all, the ``classical limit'' corresponding to Theorem~\ref{thm:HTstackformal} needs no such assumption. Based on \cite[Remark 7.8]{BMS1}, we expect some assumption akin to $\varphi$-crystallinity of $(A,I)$ is necessary to have a close relationship between prismatic crystals and $I$-connections.
\end{remark}

When $(A,I)$ is a crystalline prism, $I$-connections are simply $p$-connections, and prismatic crystals are  crystals in the sense of Grothendieck's crystalline site, up to a Frobenius twist. In this case, Theorem~\ref{thm:QuantizationPhiCrys} is, roughly,  a $p$-adic lift of the Ogus--Vologodsky correspondence \cite{OgusVologodsky}. Such lifts have been studied extensively. Without attempting to be exhaustive, let us mention \cite{ShihoLocalOV} for a local version in the presence of a Frobenius lift,  the  paper \cite{DaxinLiftCartier} (building on \cite{OyamaPDHiggs}) for a global version without Frobenius lifts but restricted to PD-nilpotent objects (see also Remark~\ref{rmk:padicCartierTransform}), and  the recent works \cite{OgusCrysPrism, WangPrismaticCrysCharp, YuPrismaticCrysCharp} for a prismatic perspective on the results just mentioned. 

\subsubsection{$B_\dR^+$-quantization over a perfectoid field}
\label{ss:BdRQuantization}

Let us now move to the rigid setting, where the quantization question is the following:

\begin{question}
\label{q:BdRQuantization}
For a smooth rigid space $X$ over a perfectoid field $C/\mathbf{Q}_p$, the twisted $p$-adic Simpson equivalence (Theorem~\ref{thm:padicSimpson}) yields an identification
\[
\begin{tikzcd}[column sep=6em]
\left\{
\parbox{4.7cm}{\centering pro-\'etale vector bundles on $X$}
\right\}
\arrow[r, dashed, leftrightarrow, "\text{twist by } \mathcal{P}_X"]
&
\left\{
\parbox{3.6cm}{\centering Higgs bundles on $X/C$}
\right\}.
\end{tikzcd}
\]
Here the Higgs field takes values in $\Omega := \Omega^1_{X/C}\{-1\}$, and $\mathcal{P}_X$ denotes the Simpson gerbe on $\Omega$. The field $C$ has a one-parameter thickening $B_\dR^+$ constructed by Fontaine, with $I := \ker(\theta:B_\dR^+ \to C)$ being principal; write $\xi \in I$ for a generator (whose choice will not matter).  Given a lift $\widetilde{X}_\infty/B_\dR^+$ of $X/C$, both sides above admit natural deformations (that in fact depend only on the mod $I^2$ lift $\widetilde{X}$) to $B_\dR^+$ given by
\[
\left\{
\parbox{3.2cm}{\centering pro-\'etale $\mathbb{B}_\dR^+$-vector bundles on $X$}
\right\}
\quad \text{and} \quad
\left\{
\parbox{3.9cm}{\centering $\xi$-connections on $\widetilde{X}_\infty/B_\dR^+$}
\right\}
\]
respectively. Naturally, one wonders: is there a twisted $B_\dR^+$-linear equivalence between these categories deforming Theorem~\ref{thm:padicSimpson}? The primary obstacle is again the twist: the Simpson gerbe lives on $\Omega$, and it is unclear what it means to deform it to $B_\dR^+$ or twist the categories appearing above by such a deformation.
\end{question}

We answer this question in \S\ref{ss:RigAnQuanNAHT} by essentially transplanting\footnote{It may be possible to bypass analytic prismatization as follows: first specialize the discussion in \S \ref{ss:QuantizePhiCrys} to the $\varphi$-crystalline prism $(A_\crys(\calO_C), IA_\crys(\calO_C))$ from Example~\ref{ex:PhiCrystallizeFontainePrism}, and then ``\'etale localize to the generic fibre'', similar to our proof of the Simpson equivalence for rigid spaces in \S \ref{sec:Simpson}. Indeed, as indicated in footnote~\ref{fn:AcrysNbhdFontainePrism}, this was the motivation to consider this prism in the first place.} the ideas discussed in \S \ref{ss:QuantizePhiCrys} into the analytic prismatization recently constructed by Ansch\"utz--Le Bras--Rodr\'iguez Camargo--Scholze \cite{ALBRCSRatAnalPrism} (see also \cite{ABLBCSdRFF, LeBrasHabThesis, HauckRatAnalSyn}). As the resulting $B_\dR^+$-deformation  (Theorem~\ref{thm:QuantNAHTrigid}) of $p$-adic Simpson equivalence is almost a verbatim copy of Theorem~\ref{thm:QuantizationPhiCrys}, we do not repeat it here. Instead, let us record the result of combining this equivalence with gerbe trivializations provided by Theorem~\ref{thm:padicSimpsonUntwisted}:

\begin{theorem}[The  untwisted quantized $p$-adic Simpson equivalence]
    \label{thm:untwistedDeformedSimpson}
Let $C$ be an algebraically closed perfectoid field. Let $X/C$ be smooth proper. Fix a $B_2(C)$-lift $\widetilde{X}$ of $X$ as well as an exponential function $\widetilde{\exp}:C \to 1+\mathfrak{m}_C$ as in Theorem~\ref{thm:padicSimpsonUntwisted}. Then there is $B_\dR^+$-linear  equivalence
\[ \Perf(X_\proet,\mathbb{B}_\dR^+) \simeq
\left\{
\parbox{6cm}{\centering $\xi$-connections on perfect complexes on (any $B_\dR^+$-lift of) $\widetilde{X}/B_\dR^+$}
\right\}
=: \Perf((\widetilde{X}/B_\dR^+)^{\xi,\dR}) \]
compatible with $\otimes$-products.
\end{theorem}

The questions considered in this section have a rich history, which we do not attempt to survey here. We mention only Fontaine's work relating almost de Rham representations of the Galois group of a discretely valued field to $t$-connections \cite{FontaineAlmostdR}, its higher-dimensional analogue due to Liu--Zhu \cite{LiuZhuRH}, the recent works \cite{LMNQWSmallRH, LMNQNilpotentRH, YuBdRRH}, and the works described in the next remark.

\begin{remark}[History of our work]
The results discussed above partly overlap with \cite[\S 4]{ALBRCSRatAnalPrism}. A first version of our work was completed and circulated informally in spring 2026, before the {\em Takagi Lectures} in June \cite{BhattTakagiVideos}. At that time \cite{ALBRCSRatAnalPrism} was not yet available; instead, we learned some key features of analytic prismatization from \cite{ABLBCSdRFF, LeBrasHabThesis} (and later used \cite{HauckRatAnalSyn} to confirm our understanding). The paper \cite{ALBRCSRatAnalPrism} appeared while these notes were being finalized for submission; we refer the reader there for its treatment of related questions.
\end{remark}

\subsubsection{Outline of the constructions}
\label{sss:OutlineQuant}

Our construction of the aforementioned quantizations relies on introducing and relating several geometric objects. In this subsection, as a blueprint to the rest of this section, we outline the main objects introduced  and the overall narrative linking them.

\begin{notation}
Fix a $\varphi$-crystalline prism $(A,I)$ throughout this section, and write $\overline{A}=A/I$. For any $Z \in \Fun(\CAlg^\an_{\Spf(A)},\Sp)$, write $\overline{Z}$ for its Hodge--Tate reduction, i.e., restriction to $\Fun(\CAlg^\an_{\Spf(\overline{A})},\Sp)$. Fix a smooth scheme $X/\Spf(\overline{A})$.
\end{notation}

  In the sequel, using transmutation, we shall introduce the following constructions:

\begin{enumerate}
    \item {\em Crystallizations:} The $\varphi$-crystalline structure on $(A,I)$ ensures that $I$ has divided powers in $A$. Consequently, there is a crystalline stack $(X/A)^\crys$, lifting the relative de Rham stack $(X/\Spf(\overline{A}))^\dR$ over the Hodge--Tate locus either once a lift of $X$ to $A/I^2$ is fixed or if one restricts attention to $\overline{A}$-algebras where $\gamma_n(I)$ vanishes for all $n \geq 1$ (Construction~\ref{cons:PhiCrysTransmute}, Lemma~\ref{lem:TwoAbarStructures}, and Warning~\ref{warn:TwoAbarOnGadR}).

    \item {\em Weil restrictions:} The stack $(X/A)^\Weil$ over $\Spf(A)$ is the (derived) Weil restriction of scalars of $X$ along $\Spf(\overline{A}) \hookrightarrow \Spf(A)$ (Definition~\ref{def:WeilStack}). The Hodge--Tate reduction $\overline{(X/A)^\Weil}$ of this stack is simply the $BT_X\{1\}$-torsor of liftings of $X$ to $A/I^2$. Since the ideal $I$ has divided powers, one can view $(X/A)^\Weil$ as a divided power thickening of a stack over $X$, so there is a natural map $\nu_X:(X/A)^\Weil \to (X/A)^\crys$ (Construction~\ref{cons:QuantizeRingStackMaps} (5)).

    \item {\em A group over $(X/A)^\Weil$:} The twisted tangent bundle $T_X\{1\} \to X$ is a relative group scheme. By functoriality, it defines a relative group stack $\mathcal{T}^\crys := (T_X\{1\}/A)^\crys \to (X/A)^\crys$. Write $\mathcal{T} = \mathcal{T}^\crys \times_{(X/A)^\crys} (X/A)^\Weil$ for its pullback along $\nu_X$ (Definition~\ref{def:dRWeil}). The fibres of $\mathcal{T}$ over geometric points of $(X/A)^\Weil$ are de Rham stacks of the tangent space at the corresponding point of $X$. 

    \item {\em A torsor structure on $(X/A)^\Prism$:} There is a natural $\mathcal{T}$-torsor $\pi_\Prism:(X/A)^\Prism \to (X/A)^\Weil$. Its Hodge--Tate reduction is the natural map from $(X/A)^\HT$, regarded as a $BT_X^\sharp\{1\}$-torsor over $X$, to its pushout along the natural map $BT_X^\sharp\{1\} \to BT_X\{1\}$ of group stacks (Proposition~\ref{prop:SqZeroGaPrism} and Corollary~\ref{cor:GerbePrismatic}).
\end{enumerate}

To discuss $I$-connections, as indicated in \S \ref{sss:lambdaconn}, we need to fix a lift $\widetilde{X}/\Spf(A/I^2)$ of $X/\Spf(\overline{A})$. Such a lift induces an identification $\overline{ (X/A)^\Weil} \simeq BT_X\{1\}$ by (2). We shall next geometrize the theory of $I$-connections:

\begin{enumerate}[resume]

    \item {\em The $I$-crystallization:} There is a stack $(\widetilde{X}/A)^{I,\crys}$ over $\Spf(A)$ whose quasi-coherent sheaf theory governs quasi-nilpotent $I$-connections on any lift of $\widetilde{X}$ to $\Spf(A)$ (Definition~\ref{def:IConnStack}). As the notation emphasizes, the further $A$-lift is not part of the data, so the category of ``quasi-nilpotent $I$-connections on an $A$-lift'' only depends on the lift modulo $I^2$.

    \item {\em A torsor structure on $(\widetilde{X}/A)^{I,\crys}$:}  There is a natural $\mathcal{T}$-torsor $\pi_{I,\crys}:(\widetilde{X}/A)^{I,\crys} \to (X/A)^\Weil$ (Corollary~\ref{cor:GerbeIdR}).

\end{enumerate}

Comparing (4) and (6) shows that $(X/A)^\Prism$ and $(\widetilde{X}/A)^{I,\crys}$ are both torsors for the same group over the same base, so the notion of prismatic crystal on $X$ relative to $A$ is a twist of the notion of a quasi-nilpotent $I$-connection on an $A$-lift of $X$. To make this useful, the twist must be made more explicit. For this, we shall show:

\begin{enumerate}[resume]
    \item {\em The Faltings stack:} The $\mathcal{T}$-torsor $\pi_\calF:(\widetilde{X}/A)^\calF \to (X/A)^\Weil$ that measures the difference between $\pi_\Prism$ and $\pi_{I,\crys}$ has an intrinsic definition via transmutation, using (roughly, but critically) the ring stack $W/I^2$ (Definition~\ref{def:FaltingsStack}, Corollary~\ref{cor:TorsorIdentity}). Moreover, it is naturally pulled back from a $\mathcal{T}^\crys$-torsor $(\widetilde{X}/A)^{\calF'} \to  (X/A)^\crys$ along the map $\nu_X$ (Corollary~\ref{cor:GerbeFaltings}).
\end{enumerate}

Thus, we obtain a geometric description of the twist relating $\pi_\Prism$ with $\pi_{I,\crys}$. To convert these geometric statements to ones involving sheaves (Theorem~\ref{thm:QuantizationPhiCrys}), one may pass to Cartier duals as in Theorem~\ref{thm:HTstackformal} (3). As these notes focus on the geometric underlying statements, we defer this conversion to \cite{BHPVQuanNAHT}.

\begin{enumerate}[resume]
\item  In \S \ref{ss:RigAnQuanNAHT}, we sketch a rigid-analytic adaptation of all the preceding constructions, leading to a $B_\dR^+$-deformation (Theorem~\ref{thm:QuantNAHTrigid}) of the twisted $p$-adic Simpson equivalence proven in Theorem~\ref{thm:padicSimpson}.
\end{enumerate}

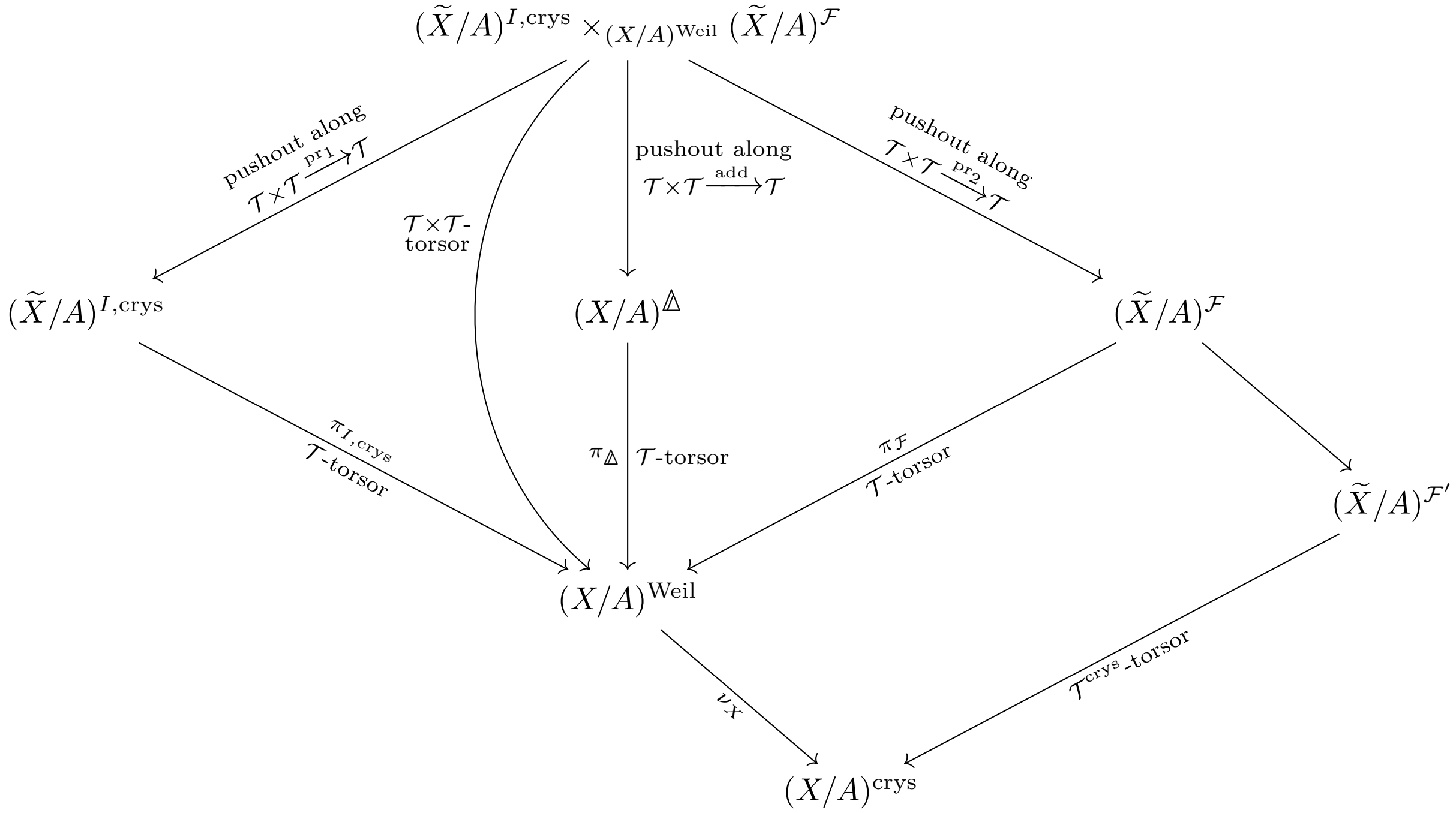
\begin{figure}[h!]
\centering
\begin{tikzcd}[column sep={8.5em,between origins}, row sep={9em,between origins}]
  & & (\widetilde{X}/A)^{I,\crys} \times_{(X/A)^{\Weil}} (\widetilde{X}/A)^{\calF}
      \arrow[dll, sloped, "\substack{\text{pushout along} \\ \mathcal{T}\times\mathcal{T} \xrightarrow{\mathrm{pr}_1} \mathcal{T}}"{pos=0.6}]
      \arrow[d, "\substack{\text{pushout along} \\ \mathcal{T}\times\mathcal{T} \xrightarrow{\mathrm{add}} \mathcal{T}}"]
      \arrow[drr, sloped, "\substack{\text{pushout along} \\ \mathcal{T}\times\mathcal{T} \xrightarrow{\mathrm{pr}_2} \mathcal{T}}"{pos=0.6}]
      \arrow[ddd, bend right=50, "\substack{\mathcal{T}\times\mathcal{T}\text{-} \\ \text{torsor}}"'{pos=0.4}]
  &[-1.5em] &[1.5em] &[-1.5em] \\
  (\widetilde{X}/A)^{I,\crys}
      \arrow[ddrr, sloped, "\pi_{I,\crys}", "\mathcal{T}\text{-torsor}"']
  & & (X/A)^{\Prism}
      \arrow[dd, "\pi_{\Prism}"', "\mathcal{T}\text{-torsor}"]
  & & (\widetilde{X}/A)^{\calF}
      \arrow[ddll, sloped, "\pi_{\calF}", "\mathcal{T}\text{-torsor}"']
      \arrow[dr]
  & \\[-3em]
  & & & & & (\widetilde{X}/A)^{\calF'}
      \arrow[ddll, sloped, "\mathcal{T}^{\crys}\text{-torsor}"']
  \\[-6em]
  & & (X/A)^{\Weil}
      \arrow[dr, sloped, "\nu_X"']
  & & & \\[-3em]
  & & & (X/A)^{\crys} & &
\end{tikzcd}
\caption{A summary of our geometric constructions (where $\mathcal{T} \sim T_X\{1\}$)}
\label{fig:GeometricSummary}
\end{figure}

The constructions introduced above are summarized by the commutative diagram in Figure~\ref{fig:GeometricSummary}, where the two parallelograms are pullback squares.

\subsection{Crystallization}
\label{ss:Crystallization}

The $\varphi$-crystalline structure on $(A,I)$ ensures that $I \subset A$ has divided powers, so there is a good notion of crystalline cohomology relative to $A$ for smooth schemes over $\Spf(\overline{A})$. We formulate this stack-theoretically, extending the construction in \S \ref{ss:CrysStack}:

\begin{construction}[The crystallization]
\label{cons:PhiCrysTransmute}
The $A$-algebra structure on $\Ga^\dR = F_*W/p$ (Proposition~\ref{prop:dRRingStackWittVectorModel}) promotes naturally to an $\overline{A}$-algebra structure: indeed, the $\varphi$-crystalline structure on $(A,I)$ gives an isomorphism 
\[ F_* W/p \simeq (F_*W)/I\]
of animated $A$-algebra stacks on $\CAlg^\an_{\Spf(A)}$; transporting the natural $\overline{A}$-algebra structure on the target gives one $\Ga^\dR = F_* W/p$. Consequently, for any $X \in \Fun(\CAlg^\an_{\Spf(\overline{A})},\Sp)$, we can define its {\em $\varphi$-crystallization $(X/A)^{\crys}$} via transmutation as the composition
\[ \CAlg^\an_{\Spf(A)} \xrightarrow{\Ga^{\dR}} \CAlg^\an_{\Spf(\overline{A})} \xrightarrow{X} \Sp. \]
If we wish to emphasize the $\overline{A}$-algebra structure on $\Ga^\dR$, we shall call it $\Ga^\crys$.
\end{construction}

\begin{example}
When $(A,(p))$ is a crystalline prism with its natural $\varphi$-crystalline structure, Construction~\ref{cons:PhiCrysTransmute} agrees with Definition~\ref{def:CrysStack}, so the notation is not abusive.
\end{example}

\begin{observation}[Two $\overline{A}$-algebra structures on $\Ga^\dR/I$]
\label{obs:TwoAbarStruct}
The animated $A$-algebra stack $\Ga^\dR/I := \Ga^\dR \otimes_A \overline{A}$ admits two natural $\overline{A}$-algebra structures: 
\begin{itemize}
    \item The {\em internal one}, coming from the first tensor factor. 
    \item The {\em external one}, coming from the second tensor factor.
\end{itemize}
    The natural map $\mu:\Ga^\dR/I \to \Ga^\dR$ is an $\overline{A}$-algebra map for either structure, while the tautological quotient map $\tau:\Ga^\dR \to \Ga^\dR/I$ is an $\overline{A}$-algebra map with respect to the internal $\overline{A}$-algebra structure on the target. {\em In the sequel, unless otherwise specified, we shall always endow $\Ga^\dR/I$ with the external $\overline{A}$-algebra structure.}
\end{observation}

While the two $\overline{A}$-structures on $\Ga^\dR/I$ are distinct, they induce the same $A/I^2$-algebra structure by Lemma~\ref{lem:TwoAbarStructures}.

\begin{warning}
\label{warn:TwoAbarOnGadR}
On the Hodge--Tate locus, there are {\em two} natural $\overline{A}$-algebra structures on $\Ga^\dR$ itself: one from Construction~\ref{cons:PhiCrysTransmute}, and the other given by the map $\Ga \to \Ga^\dR$ and the evident $\overline{A}$-algebra structure on $\Ga$ on the Hodge--Tate locus. These two structures are typically different (but induce the same $A/I^2$-structure by Lemma~\ref{lem:TwoAbarStructures}). To understand the difference, assume $(A,I)$ is a classical prism. Each $\overline{A}$-structure corresponds to an $A$-module null-homotopy of the composition $I \to A \to \Ga^\dR(\overline{A})$, or equivalently to an $A$-module lift $I \to \Ga^\sharp(\overline{A})$ of $I \to A \to \overline{A}$. The lift in Construction~\ref{cons:PhiCrysTransmute} corresponds to the composition $I \to \Ga^\sharp(A) \to \Ga^\sharp(\overline{A})$, where the first map comes from the divided power structure on $I$ provided by the $\varphi$-crystalline structure and the second map is the obvious one. The lift corresponding to the map $\overline{A} \to \Ga \to \Ga^\dR$ results from noticing that $I \to A \to \overline{A}$ is $0$, so the $0$ map $I \to \Ga^\sharp(\overline{A})$ provides the desired lift. For the two to coincide, we need $\gamma_n(I) \subset I$ for all $n \geq 1$; this holds true for crystalline prisms (as $\mathrm{val}_p(\gamma_n(p)) \geq 1$ for all $n \geq 1$) but not for the prism $(A_\crys(R),IA_\crys(R))$ from Example~\ref{ex:PhiCrystallizeFontainePrism}; in the latter case, the difference disappears if one restricts to test algebras over the quotient $R$ of $A_\crys(R)/IA_\crys(R)$. 
\end{warning}

\subsection{Weil restriction}

Classically, given a finite flat map $f:X \to Y$ of schemes and an $X$-scheme $Z$, its Weil restriction $\mathrm{Res}^X_Y(Z) \to Y$ is the $Y$-scheme representing the functor $f_* Z$, see \cite[\href{https://stacks.math.columbia.edu/tag/05YF}{Tag 05YF}]{StacksProject}. Once one passes to the derived setting and also allows stacks, this construction works without flatness constraints on $f$ and can be applied to closed immersion $\Spf(\overline{A}) \hookrightarrow \Spf(A)$, see \cite[Theorem C]{HKRWeilRestriction}. Let us give a transmutation based construction of it:

\begin{definition}[The Weil stack]
\label{def:WeilStack}
Let $\Ga^{\Weil}$ be the animated $\overline{A}$-algebra stack on $\CAlg^\an_{\Spf(A)}$ determined by $\Ga^{\Weil}(R) = R \otimes_A \overline{A}$, with its obvious $\overline{A}$-algebra structure. For any $X \in \Fun(\CAlg^\an_{\Spf(\overline{A})},\Sp)$, write $(X/A)^{\Weil}$ for its transmutation with respect to $\Ga^{\Weil}$, i.e., $(X/A)^{\Weil}$ is the composition
\[ \CAlg^\an_{\Spf(A)} \xrightarrow{\Ga^{\Weil}} \CAlg^\an_{\Spf(\overline{A})} \xrightarrow{X} \Sp. \]
We shall call this the {\em Weil stack} of $X$ relative to $A$.
\end{definition}

Unwinding definitions, for any $X \in \Fun(\CAlg^\an_{\Spf(\overline{A})},\Sp)$, there is a natural projection map $\overline{(X/A)^{\Weil}} \to X$. In fact, in the realization of $(X/A)^{\Weil}$ as the Weil restriction of $X$ along $\Spf(\overline{A}) \to \Spf(A)$ (which is right adjoint to Hodge--Tate reduction), the preceding map is the counit map of the adjunction.

\begin{remark}[The Hodge--Tate reduction of $(X/A)^{\Weil}$]
\label{rmk:HTlocusNaivePrism}
On the Hodge--Tate locus, the $A$-algebra structure on $\Ga^{\Weil}$ admits two $\overline{A}$-refinements: the external one (that exists before Hodge--Tate reduction), and the internal one (coming from the natural map $\overline{A} \to \Ga \to \Ga^{\Weil}$). These induce the same $A/I^2$-structure by Lemma~\ref{lem:TwoAbarStructures}. The natural map $\Ga^{\Weil} = \Ga/I \to \Ga$ realizes the source as a square-zero extension of $\Ga$ by $B\Ga\{1\}$ in $\overline{A}$-algebras (for either $\overline{A}$-algebra structure on the source), while the natural map $\Ga \to \Ga^{\Weil} = \Ga/I$ is $\overline{A}$-linear for the internal $\overline{A}$-algebra structure on the target, and hence $A/I^2$-linear for either $\overline{A}$-algebra on the target. By transmutation, given $X \in \Fun(\CAlg^\an_{\Spf(\overline{A})},\Sp)$, the choice of a lift $\widetilde{X}$ to $A/I^2$ yields a section to the projection map $\overline{(X/A)^{\Weil}} \to X$. In fact, for $X/\Spf(\overline{A})$ smooth, this construction identifies the projection $\overline{(X/A)^{\Weil}} \to X$ with the $BT_X\{1\}$-torsor of liftings of $X$ to $A/I^2$.
\end{remark}

\begin{remark}[The Weil stack as the naive prismatization]
The stack $(X/A)^\Weil$ can also be regarded as the ``naive prismatization'' of $X$. More precisely, if $X/\Spf(\overline{A})$ is smooth, and one defines $(X/A)_{\Prism,\text{naive}}$ exactly as in the definition of prismatic site \cite{BhattScholzePrism} but omitting all $\delta$-structures, then one can show that $R\Gamma( (X/A)^\Weil, \calO) \simeq R\Gamma( (X/A)_{\Prism,\text{naive}},\calO)$. The equal characteristic $0$ counterpart of this cohomology theory is well-behaved: it essentially recovers filtered de Rham cohomology by the unpublished Bonn masters thesis of Nils Wassmuth  (see also \cite[Theorem 1.6]{AntieauFiltrationsII}). On the other hand, in our $p$-adic world, the object $R\Gamma( (X/A)^\Weil, \calO)$ is large and somewhat pathological (e.g., via Remark~\ref{rmk:HTlocusNaivePrism}); nevertheless, the stack $(X/A)^\Weil$ itself plays an important role. 
\end{remark}

\subsection{A bestiary of ring stack maps}

Our subsequent constructions shall rely  heavily on transmutation of various maps of ring stacks. Instead of defining the maps as we need them, we collect them  for ease of reference and comparison:

\begin{construction}[Compendium]
\label{cons:QuantizeRingStackMaps}
We shall need the following maps of animated $\overline{A}$-algebra stacks:
\begin{enumerate}
    \item The $0$-th Witt component map $W \to \Ga$ reduces modulo $I$ to define a map 
    \[ \pi_{\Prism}:\Ga^\Prism \to \Ga^{\Weil}\] of animated $\overline{A}$-algebra stacks.

    \item The Witt vector Frobenius $F: W \to F_*W$ reduces modulo $I$ to a map
    \[ F:\Ga^\Prism = W/I \to (F_*W)/I \simeq \Ga^\dR \]
    of animated $\overline{A}$-algebra stacks, where the last identification uses the $\varphi$-crystalline structure (Construction~\ref{cons:PhiCrysTransmute}).

    \item There is a natural animated $\overline{A}$-algebra map 
    \[ \mu:\Ga^\dR/I = \Ga^\dR \otimes_A \overline{A} \to \Ga^\dR\] 
    coming from the natural $\overline{A}$-algebra structure on $\Ga^\dR$ (Observation~\ref{obs:TwoAbarStruct}); this is an $\overline{A}$-algebra map for both the external and internal $\overline{A}$-algebra structures on the source. (Recall that our default $\overline{A}$-structure on $\Ga^\dR/I$ is the external one.)

    \item Reducing the tautological map $\Ga \to \Ga^\dR$ modulo $I$ gives a natural animated $\overline{A}$-algebra map 
    \[ \alpha:\Ga^{\Weil} =\Ga/I \to \Ga^\dR/I.\] 

    \item Combining the previous two maps gives the map
    \[ \nu = \mu \circ \alpha:\Ga^{\Weil} \to \Ga^\dR\]
    which is an $\overline{A}$-algebra map for the natural $\overline{A}$-algebra structures on both sides. Using this, we can regard any $\Ga^\dR$-module as a $\Ga^{\Weil}$-module. 

    \item The tautological quotient map 
    \[ \tau:\Ga^\dR \to \Ga^\dR/I\]
    is an animated $\overline{A}$-algebra map with respect to the internal $\overline{A}$-algebra structure on the target, and hence an $A/I^2$-algebra map (see Observation~\ref{obs:TwoAbarStruct}). Moreover:

    \begin{lemma}
    \label{lem:tauSqZero}
The map $\tau$ promotes to a square-zero extension of $\Ga^\dR/I$ by $\Ga^\dR\{1\}$ (regarded as a module over $\Ga^\dR/I$ via $\mu$) in $A/I^2$-algebras. 
    \end{lemma}
    \begin{proof}
    In the Witt vector model $\Ga^\dR = F_*W \otimes_A \overline{A}$, the map $\tau$ is given by the map
    \[ F_*W \otimes_A \overline{A} \to \left(F_* W \otimes_A \overline{A}\right) \otimes_A \overline{A}  \]
    sending the LHS to the first tensor factor of the RHS (with $\overline{A}$-structures on both sides coming from the outermost copies of $\overline{A}$). In other words, it is obtained by applying $F_* W \otimes_A (-)$ to the map $i_1:\overline{A} \to \overline{A} \otimes_A \overline{A}$ given by the first factor inclusion. Endowing $\overline{A} \otimes_A \overline{A}$ with $\overline{A}$-structure via the second factor, the map $i_1$ is naturally a square-zero extension of its target by $\overline{A}\{1\}$ in $A/I^2$-algebras (see \cite[Theorem 7.4.1.26]{LurieHA}), so the claim follows by base change.
    \end{proof}

    \item Reducing the tautological map $F_*W \to F_*W/p \simeq (F_*W)/I$ modulo $I$, and regarding $F_*W/p$ as $\Ga^\dR$, gives an animated $\overline{A}$-algebra map 
    \[\overline{\theta}:\Ga^\dR \simeq F_*W/I \to \Ga^\dR/I\] with respect to the external $\overline{A}$-algebra structure on the target (Observation~\ref{obs:TwoAbarStruct}).
    \begin{lemma} 
    \label{lem:ThetaSqZero}
    The map $\overline{\theta}$ promotes to a square-zero extension of $\Ga^\dR/I$ by $\Ga^\dR\{1\}$ (regarded as a module over $\Ga^\dR/I$ via $\mu$) in $\overline{A}$-algebras. 
    \end{lemma}
    \begin{proof}
    In the Witt vector model, $\overline{\theta}$ is simply the map
    \[ \left(F_* W \to F_* W \otimes_A \overline{A}\right) \otimes_A \overline{A},\]
    with the $\overline{A}$-algebra structures coming from the outer copy of $\overline{A}$. Thus, $\overline{\theta}$ is obtained by applying $F_*W \otimes_A (-)$ to the second factor inclusion $i_2:\overline{A} \to \overline{A} \otimes_A \overline{A}$. The latter is naturally a square-zero extension of the target by $\overline{A}\{1\}$ in $\overline{A}$-algebras (for the $\overline{A}$-algebra structure coming from the second factor), so the claim follows by base change. 
    \end{proof}
\end{enumerate}
The maps above may be summarized in the commutative diagram
\begin{equation}
\tag{Rings}
\label{eq:Bestiary}
\begin{tikzcd}[column sep=4em, row sep=2.5em]
W/I = \Ga^\Prism \arrow[r, "F"] \arrow[d, "\pi_\Prism"'] & (F_*W)/I \simeq F_*W/p \simeq \Ga/\Ga^\sharp = \Ga^\dR \arrow[d, shift right=0.7ex, "\overline{\theta}"'] \arrow[d, shift left=0.7ex, "\tau"] \\
\Ga/I = \Ga^{\Weil} \arrow[r, "\alpha"] \arrow[dr, "\nu"'] & \Ga^\dR/I \arrow[d, "\mu"] \\
& \Ga^\dR.
\end{tikzcd}
\end{equation}
Components of this diagram will be elaborated upon in the subsequent sections. But let us indicate right away that the difference between $\overline{\theta}$ and $\tau$ is ultimately responsible for the twist appearing in both the $p$-adic Simpson correspondence and its quantization (see Corollary~\ref{cor:TorsorIdentity}).

\end{construction}

\subsection{The prismatic stack as a torsor}
\label{ss:PrismaticGerbe}
In this subsection, we explain why a shadow of the gerbe structure on the Hodge--Tate stack survives quantization; this is a prismatic variant of Proposition~\ref{prop:SqZeroHT} and was the  starting point of the work discussed in this section:

\begin{proposition}[A square-zero ideal in $\Ga^\Prism$]
\label{prop:SqZeroGaPrism}
There is a natural pullback diagram of animated $\overline{A}$-algebra stacks on $\CAlg^\an_{\Spf(A)}$ (extracted from \eqref{eq:Bestiary}):
\begin{equation}
\tag{$\Ga^\Prism$-$\square$}
\label{eq:WmodIsqzero} \xymatrix{ \Ga^\Prism := W/I \ar[r]^-F \ar[d]^-{\pi_\Prism} & (F_*W)/I \simeq \Ga^\dR \ar[d]^{\overline{\theta}} \\ 
\Ga^{\Weil} := \Ga/I \ar[r]^{\alpha} & \Ga^\dR/I. }
\end{equation}
Moreover, the right vertical map, and hence the left vertical map, is a square-zero extension of its target by $\Ga^\dR\{1\}$ in animated $\overline{A}$-algebras. Here $\Ga^\dR\{1\}$ is regarded as a $\Ga^\dR/I$-module via $\mu$ (Lemma~\ref{lem:ThetaSqZero}), and hence as a $\Ga^{\Weil}$-module via $\nu = \mu \circ \alpha$ (see \eqref{eq:Bestiary}).
\end{proposition}
\begin{proof}
The natural identification $F_* W/p \simeq \Ga^\dR$ can be rewritten as a pullback square
\[ \xymatrix{ W \ar[r] \ar[d] & F_*W \ar[d] \\ \Ga \ar[r] & \Ga^\dR \simeq F_*W/p }\]
of animated $A$-algebra stacks. Reducing modulo $I$ then gives the first part. The second part follows from Lemma~\ref{lem:ThetaSqZero} and the first part.
\end{proof}

\begin{remark}
\label{rmk:PrismFibNoPhiCrys}
The fibre square \eqref{eq:WmodIsqzero}, minus the isomorphism of its top right term with $\Ga^\dR$, exists without the $\varphi$-crystalline assumption on the animated prism $(A,I)$, as the proof of Proposition~\ref{prop:SqZeroGaPrism} makes clear. The $\varphi$-crystalline assumption is needed only for the square-zero property: it identifies $(F_*W)/I$ with $\Ga^\dR$ (Construction~\ref{cons:PhiCrysTransmute}), and this makes $\overline{\theta}$ a square-zero extension (Lemma~\ref{lem:ThetaSqZero}).
\end{remark}

To extract geometric consequences, it is convenient to make the following definition (whose notation is inspired by the last sentence of Construction~\ref{cons:PhiCrysTransmute}):

\begin{definition}[Deforming the relative de Rham stack to the Weil stack]
\label{def:dRWeil}
For any $X \in \Fun(\CAlg^\an_{\Spf(\overline{A})},\Sp)$ and a map $E \to X$, define $(E/X)^\crys_\Weil$ as a pullback
\[ \xymatrix{ (E/X)^\crys_\Weil \ar[r] \ar[d] & (E/A)^\crys \ar[d] \\
(X/A)^{\Weil} \ar[r]^{\nu_X} & (X/A)^\crys, }\]
where the bottom horizontal map is the transmutation of $X$ along $\nu:\Ga^{\Weil} \to \Ga^\dR$. 

For $X/\Spf(\overline{A})$ smooth, we write $\mathcal{T}^\crys := (T_X\{1\}/A)^\crys \to (X/A)^\crys$ and $\mathcal{T} := (T_X\{1\}/X)^\crys_\Weil \to (X/A)^\Weil$, so $\mathcal{T}$ is the pullback of $\mathcal{T}^\crys$ along $\nu_X$ (as in  \S \ref{sss:OutlineQuant}).
\end{definition}

%We shall apply this construction when $E$ is a vector bundle over $X$ (e.g., the tangent or cotangent bundle). 
Note that if $E=\Ga^d \times X$ for some integer $d \geq 0$, then $(E/X)^\crys_\Weil \simeq (\Ga^\dR)^d \times (X/A)^{\Weil}$.  In terms of this object, Proposition~\ref{prop:SqZeroGaPrism} transmutes to the following:

\begin{corollary}[The torsor structure on the prismatization]
\label{cor:GerbePrismatic}
For a smooth scheme $X/\Spf(\overline{A})$, there is a natural map $\pi_\Prism:(X/A)^\Prism \to (X/A)^{\Weil}$, induced via transmutation from $\pi_\Prism:\Ga^\Prism \to \Ga^{\Weil}$. Moreover, this map is naturally a $\mathcal{T}$-torsor. 
\end{corollary}

The above can be regarded as a deformation towards $\Spf(A)$ of the reduction of structure group of the Hodge--Tate stack provided by an $A/I^2$-lift of $X$:

\begin{remark}[Relation to the Hodge--Tate gerbe]
In the setup of Corollary~\ref{cor:GerbePrismatic}, assume $X$ comes equipped with a lift $\widetilde{X}$ to $A/I^2$. This choice defines a map $s:X \to \overline{(X/A)^{\Weil}}$ splitting the projection (see Remark~\ref{rmk:HTlocusNaivePrism}); the composition $X \xrightarrow{s} (X/A)^{\Weil} \xrightarrow{\nu_X} (X/A)^\crys$ is the tautological map. The pullback of $\pi_\Prism:(X/A)^\Prism \to (X/A)^{\Weil}$  along $s$ therefore gives a $(T_X\{1\}/X)^\dR$-torsor over $X$. It is a pleasant exercise to show that this coincides with the $(T_X\{1\}/X)^\dR$-reduction of the Hodge--Tate stack $(X/A)^\HT \to X$ provided by the choice of $\widetilde{X}$ (Theorem~\ref{thm:HTstackformal} (2)).
\end{remark}

\subsection{The $I$-crystallization}
\label{ss:IConnStack}

In \S \ref{sss:lambdaconn}, we reviewed the notion of an $I$-connection (an $I$-adic deformation of the notion of a Higgs bundle), and gave a stack-theoretic perspective on this notion  in characteristic $0$ (Proposition~\ref{prop:IConnModI2}). In this subsection, we explain the analogous construction in the prismatic setting. The key player is the following:

\begin{definition}[The $I$-crystallization]
\label{def:IConnStack}
\begin{enumerate}
\item Let $\Ga^{I,\crys}$, the {\em $I$-crystalline ring stack}, be the animated $A/I^2$-algebra stack on the category $\CAlg^\an_{\Spf(A)}$ defined by the following pullback diagram of animated $A/I^2$-algebra stacks (extracted from \eqref{eq:Bestiary}):
\begin{equation}
\tag{$\Ga^{I,\crys}$-$\square$}
\label{eq:IconnRingStack}
     \xymatrix{ \Ga^{I,\crys} \ar[r] \ar[d]^-{\pi_{I,\crys}} & \Ga^\dR \ar[d]^-\tau \\ 
\Ga^\Weil = \Ga/I \ar[r]^-\alpha & \Ga^\dR/I }
\end{equation}
where $\alpha$ and $\tau$ come from Construction~\ref{cons:QuantizeRingStackMaps}, and their $A/I^2$-algebra structure comes from Observation~\ref{obs:TwoAbarStruct}.

\item For any $Y \in \Fun(\CAlg^\an_{\Spf(A/I^2)},\Sp)$, write $(Y/A)^{I,\crys}$ for its transmutation with respect to $\Ga^{I,\crys}$, i.e., 
\[ \CAlg^\an_{\Spf(A)} \xrightarrow{\Ga^{I,\crys}} \CAlg^\an_{\Spf(A/I^2)} \xrightarrow{Y} \Sp. \]
We call $(Y/A)^{I,\crys}$ the {\em $I$-crystallization} of $Y$; it controls $I$-connections on lifts of $Y$ (Remark~\ref{rmk:IConnLiftDepend}). If $Y_{\overline{A}}$ denotes the base change of $Y$ to $\Spf(\overline{A})$, the map $\pi_{I,\crys}$ in \eqref{eq:IconnRingStack} transmutes to a map $\pi_{I,\crys}:(Y/A)^{I,\crys} \to (Y_{\overline{A}}/A)^{\Weil}$. 
\end{enumerate}
\end{definition}

More informally, the ring stack $\Ga^{I,\crys}$ has a quotient description
\[ \Ga^{I,\crys} = \Ga/I\Ga^\sharp,\] 
with fibre square \eqref{eq:IconnRingStack} expressing Milnor glueing for the product ideal $I\Ga^\sharp$. Notice that the shape of fibre square \eqref{eq:IconnRingStack} is analogous to that of \eqref{eq:WmodIsqzero}. Moreover, by Lemma~\ref{lem:tauSqZero}, there are even more similarities:

\begin{proposition}[A square-zero ideal in $\Ga^{I,\crys}$]
\label{prop:SqZeroGaIdR}
The right vertical map, and hence the left vertical map, of \eqref{eq:IconnRingStack} is naturally a square-zero extension of the target by $\Ga^\dR\{1\}$ (regarded as a $\Ga^{\Weil}$-module via $\nu$, as in Proposition~\ref{prop:SqZeroGaPrism}) in $A/I^2$-algebras.
\end{proposition}

As in Remark~\ref{rmk:PrismFibNoPhiCrys}, the fibre square \eqref{eq:IconnRingStack} exists (in $A$-algebra stacks) without the $\varphi$-crystalline assumption; the latter is only used to realize the vertical maps as square-zero extensions.

\begin{remark}[Relating to $I$-connections]
\label{rmk:IConnLiftDepend}
Fix a smooth scheme $X/\Spf(\overline{A})$ with a lift $\widetilde{X}_\infty/\Spf(A)$, and write $\widetilde{X} := \widetilde{X}_\infty \otimes_A A/I^2$. There is a natural map $\widetilde{X}_\infty \to (\widetilde{X}/A)^{I,\crys}$ coming from the $\Ga$-algebra structure on $\Ga^{I,\crys}$. This map is faithfully flat, and making descent along this map explicit leads to the following: the category $\mathrm{Vect}( (\widetilde{X}/A)^{I,\crys})$ identifies with the category of quasi-nilpotent flat $I$-connections on $\widetilde{X}_\infty/A$, proving that this notion only depends on $\widetilde{X}$ and not $\widetilde{X}_\infty$ itself; see \S \ref{sss:lambdaconn} for a parallel discussion in the  characteristic $0$ setting.
\end{remark}

By deformation theory, Proposition~\ref{prop:SqZeroGaIdR} yields:

\begin{corollary}[The torsor structure on the $I$-crystallization]
\label{cor:GerbeIdR}
For a smooth scheme $X/\Spf(\overline{A})$ with a lift $\widetilde{X}/\Spf(A/I^2)$, the natural map $\pi_{I,\crys}:(\widetilde{X}/A)^{I,\crys} \to (X/A)^{\Weil}$ is naturally a $\mathcal{T}$-torsor. 
\end{corollary}

\begin{remark}[The Hodge--Tate locus of the $I$-crystallization]
Fix the notation of Corollary~\ref{cor:GerbeIdR}. By the choice of $\widetilde{X}$, the Hodge--Tate reduction $\overline{ (X/A)^\Weil}$ is identified with $BT_X\{1\}$ as in Remark~\ref{rmk:HTlocusNaivePrism}. Using the obvious $\Ga$-algebra structure on $\Ga^{I,\crys}$, it is a pleasant exercise to identify $\overline{ (\widetilde{X}/A)^{I,\crys} }$ with $BT_X^\sharp\{1\}$, and the map $\overline{\pi_{I,\crys}}$ with the evident map $BT_X^\sharp\{1\} \to BT_X\{1\}$. 
\end{remark}

\begin{remark}[The $p$-adic Cartier transform]
\label{rmk:padicCartierTransform}
For $X/\Spf(\overline{A})$ smooth, the existence of the torsors $\pi_\Prism$ and $\pi_{I,\crys}$ over the common base $(X/A)^\Weil$ yields a $\varphi$-crystalline generalization of the $p$-adic Cartier transform \cite{DaxinLiftCartier}, identifying the categories of prismatic crystals and quasi-nilpotent $I$-connections equipped with a system of PD-nilpotent divided powers. Indeed, the latter categories can both be identified with $\mathcal{D}( (X/A)^\Weil)$ via pullback along $\pi_\Prism^*$ or $\pi_{I,\crys}^*$ respectively.  
\end{remark}

\subsection{The main theorem: an identification of torsors}
\label{ss:MainThmTorsors}

 \begin{notation} 
 In this subsection, fix a smooth scheme $X/\Spf(\overline{A})$ with a lift $\widetilde{X}/\Spf(A/I^2)$. 
 \end{notation}

 In \S \ref{ss:PrismaticGerbe} and \S \ref{ss:IConnStack}, we have used transmutation to construct two maps
\[ \pi_\Prism:(X/A)^\Prism \to (X/A)^{\Weil} \quad \text{and} \quad \pi_{I,\crys}:(\widetilde{X}/A)^{I,\crys} \to (X/A)^\Weil, \]
both of which are $\mathcal{T}$-torsors over $(X/A)^{\Weil}$ (Corollaries~\ref{cor:GerbePrismatic} and \ref{cor:GerbeIdR}). In this subsection, we describe their difference via transmutation of the following:

\begin{definition}[The Faltings stack]
\label{def:FaltingsStack}
\begin{enumerate}
\item Let $\Ga^\calF$ be the animated $A/I^2$-algebra stack on $\CAlg^\an_{\Spf(A)}$ defined by the following pullback diagram:
\begin{equation}
    \tag{$\Ga^\calF$-fib}
    \label{eq:FaltingsRingStack}
    \xymatrix{ \Ga^\calF \ar[r] \ar[d]^-{\pi_\calF} & F_*W/I^2 \ar[d]^-{\text{can}} \\ 
    \Ga^{\Weil} = \Ga/I \ar[r]^-\nu & \Ga^\dR \simeq (F_*W)/I, }
\end{equation}
where $\nu$ comes from Construction~\ref{cons:QuantizeRingStackMaps} and the right vertical map is the canonical one.

\item For any $Y \in \Fun(\CAlg^\an_{\Spf(A/I^2)},\Sp)$, write $(Y/A)^\calF$ for its transmutation with respect to $\Ga^\calF$, i.e., 
\[ \CAlg^\an_{\Spf(A)} \xrightarrow{\Ga^\calF} \CAlg^\an_{\Spf(A/I^2)} \xrightarrow{Y} \Sp. \]
If $Y_{\overline{A}}$ denotes the base change of $Y$ to $\Spf(\overline{A})$, the map $\pi_\calF$ in \eqref{eq:FaltingsRingStack} transmutes to a map $\pi_\calF:(Y/A)^\calF \to (Y_{\overline{A}}/A)^{\Weil}$. 
\end{enumerate}
\end{definition}

\begin{remark}[Why ``Faltings stack''?]
The name ``Faltings stack'' references the Faltings extension in $p$-adic Hodge theory \cite[Corollary 6.14]{ScholzePAdicHodge}: in the rigid setting of \S \ref{ss:RigAnQuanNAHT}, the counterpart of the right vertical arrow in \eqref{eq:FaltingsRingStack} is the square-zero extension that gives rise to the Faltings extension.
\end{remark}

The following is immediate:

\begin{proposition}[A square-zero ideal in $\Ga^{\calF}$]
\label{prop:SqZeroGaF}
The right vertical map, and hence the left vertical map, of \eqref{eq:FaltingsRingStack} is naturally a square-zero extension of the target by $\Ga^\dR\{1\}$ (regarded as a $\Ga^{\Weil}$-module via $\nu$, as in Proposition~\ref{prop:SqZeroGaPrism}) in animated $A/I^2$-algebra stacks.
\end{proposition}

By deformation theory, Proposition~\ref{prop:SqZeroGaF} and the pullback square \eqref{eq:FaltingsRingStack} yield:

\begin{corollary}[The torsor structure on the Faltings stack and its crystalline descent]
\label{cor:GerbeFaltings}
The natural map $\pi_\calF:(\widetilde{X}/A)^\calF \to (X/A)^{\Weil}$ is naturally a $\mathcal{T}$-torsor. Moreover, this $\mathcal{T}$-torsor is naturally the pullback of a $\mathcal{T}^\crys$-torsor $\pi_{\calF'}:(\widetilde{X}/A)^{\calF'} \to (X/A)^\crys$ along $\nu_X:(X/A)^{\Weil} \to (X/A)^\crys$.
\end{corollary}

The Faltings ring stack $\Ga^\calF$ measures the failure of the prismatic ring stack $\Ga^\Prism$ and the $I$-crystalline ring stack $\Ga^{I,\crys}$ to agree as square-zero extensions of $\Ga^\Weil$:

\begin{proposition}
    \label{prop:TorsorAdditionIdentity}
We have an identity
\[ [\Ga^\Prism] = [\Ga^{I,\crys}] + [\Ga^\calF]\]
in the $E_\infty$-group of square-zero extensions of $\Ga^{\Weil}$ by $\Ga^\dR\{1\}$ in $A/I^2$-algebras. 
\end{proposition}
\begin{proof}[Proof idea]
The map $\nu:\Ga^\Weil \to \Ga^\dR$ appearing in \eqref{eq:FaltingsRingStack} factors as $\Ga^\Weil = \Ga/I \xrightarrow{\alpha} \Ga^\dR/I \xrightarrow{\mu} \Ga^\dR$ (see Construction~\ref{cons:QuantizeRingStackMaps}). Contemplating the squares \eqref{eq:WmodIsqzero} and \eqref{eq:IconnRingStack}, the square \eqref{eq:FaltingsRingStack} then shows that all three square-zero extensions of $\Ga^\Weil$ under consideration are pulled back from square-zero extensions of $\Ga^\dR/I$ under $\alpha:\Ga^\Weil \to \Ga^\dR/I$. Moreover, the only non-formal input needed to write down all three is the $A$-deformation $F_*W$ of $\Ga^\dR=F_*W/I$. Thus, one reduces to a general assertion about an animated ring $B$ equipped with an invertible ideal $J \to B$; this reduces to the universal case, in which case one can argue using natural quasi-ideal models. 
\end{proof}

Transmuting the above gives the main result of this section:

\begin{corollary}
\label{cor:TorsorIdentity}
There is a natural identification
    \[ [ (X/A)^\Prism ] = [ (\widetilde{X}/A)^{I,\crys}] + [ (\widetilde{X}/A)^\calF ] \]
    in the $E_\infty$-group of $\mathcal{T}$-torsors over $(X/A)^{\Weil}$. 
\end{corollary}

The sheaf-theoretic consequence mentioned in Theorem~\ref{thm:QuantizationPhiCrys} is obtained by passing to Cartier duals, as was done in the Hodge--Tate story in Theorem~\ref{thm:HTstackformal}; we do not elaborate further on this here.

\subsection{Passage to rigid analytic spaces}
\label{ss:RigAnQuanNAHT}

In this subsection, we answer the $B_\dR^+$-quantization question (Question~\ref{q:BdRQuantization}) using the analytic prismatization of \cite{ALBRCSRatAnalPrism} (which is built on analytic geometry \cite{ABLBCSdRFF} and hence on condensed mathematics \cite{ScholzeCondensed}). In fact, our constructions from \S\ref{ss:QuantIntro}--\ref{ss:MainThmTorsors} can be transplanted into this framework quite effortlessly to prove Theorem~\ref{thm:QuantNAHTrigid}, so we shall be quite brief. We begin by introducing the basic geometric categories.

\begin{construction}[Nilperfectoids]
There is an $\infty$-category $\NilPerfd_{\mathbf{Q}_p}$ of (animated) nilperfectoid $\mathbf{Q}_p$-algebras $R$ that is a subcategory of condensed $\mathbf{Q}_p$-algebras that minimally contains perfectoid $\mathbf{Q}_p$-algebras as well as their thickenings by $\dagger$-nilpotent ideals. More precisely, any $R \in \NilPerfd_{\mathbf{Q}_p}$ sits in a canonical fibre sequence
\begin{equation}
\tag{$\dagger$-red}
\label{eq:DaggerNilRad}
    \nilraddagger{R} \to R \to R^\daggerred,
\end{equation} 
where $R^\daggerred$ is perfectoid, and $\nilraddagger{R}$ is the ideal of $\dagger$-nilpotent elements. This category provides the test objects under consideration, so all geometric objects are viewed as functors\footnote{This move is analogous to the following more familiar one from algebraic geometry: the category of schemes over $\mathbf{F}_p$ embeds fully faithfully inside  presheaves on the category of $\mathbf{F}_p$-algebras which are semiperfect (i.e., have a surjective Frobenius): this follows from faithfully flat descent as any $\mathbf{F}_p$-algebra admits a faithfully flat cover which is semiperfect. An even closer analog is provided by viewing derived $p$-adic formal schemes as functors on semiperfectoids, i.e., $p$-nilpotent animated rings $R$ with $\pi_0(R/p)$ semiperfect.} on $\NilPerfd_{\mathbf{Q}_p}$.
\end{construction}

Two examples of relevance to the sequel are:

\begin{example}[de Rham stacks]
The sequence in \eqref{eq:DaggerNilRad} is functorial in $R$; the corresponding fibre sequence in $\Fun(\NilPerfd_{\mathbf{Q}_p},\Sp)$ is denoted
\[ \Ga^\dagger \to \Ga \to \Ga^\dR.\]
Note that $\Ga^\dR$ is an endofunctor of $\NilPerfd_{\mathbf{Q}_p}$. By transmutation, this yields Rodr\'iguez Camargo's theory of analytic de Rham stacks $X \mapsto X^\dR$ for $X \in \Fun(\NilPerfd_{\mathbf{Q}_p},\Sp)$; write $(X/Y)^\dR$ for the relative analytic de Rham stack of a morphism $X \to Y$ in $\Fun(\NilPerfd_{\mathbf{Q}_p},\Sp)$. For more on this notion, see \cite{camargo2024analytic}.
\end{example}

\begin{example}[Rigid spaces]
The category of partially proper rigid spaces\footnote{
As far as the author understands, the restriction to partially proper rigid spaces is a reflection of the fact that a nilperfectoid $R$ does not come equipped with a preferred subring $R^+ \subset R$ of ``functions bounded by $1$''. While somewhat restrictive (e.g., the closed unit disc is excluded), this category still includes many rigid spaces of interest, e.g., any proper rigid space, or the analytification of any algebraic variety.} over $\mathbf{Q}_p$ embeds fully faithfully (via Yoneda) into $\Fun(\NilPerfd_{\mathbf{Q}_p},\Sp)$. In particular, unlike Scholze's theory of diamonds \cite{ScholzeDiamonds} where the test objects are perfectoid, this theory detects non-topological information, such as nilpotents in the structure sheaf.
\end{example}

Specializing to the problem at hand, we shall work in the following geometric context.

\begin{notation}
Fix a perfectoid field $C/\mathbf{Q}_p$,  let $\mathcal{Y}_{(0,\infty)} := \mathrm{Spa}(A_{\inf}(\calO_C))^\an$ be the associated curve, and let $\mathrm{Spa}(C) \hookrightarrow \mathcal{Y}_{(0,\infty)}$ be the Hodge--Tate point. Let $\mathcal{Y}^\dagger := \left(\mathrm{Spa}(C) \subset \mathcal{Y}_{(0,\infty)}\right)^\dagger$ be the overconvergent neighbourhood of the Hodge--Tate point,  let $I \subset \calO_{\mathcal{Y}^\dagger}$ be the corresponding invertible ideal, and write $M\{1\} = M \otimes_C I/I^2$ for $C$-modules $M$. Let $\mathrm{Spa}(B_2(C)) \subset \mathcal{Y}^\dagger$ be the first order infinitesimal neighbourhood of the Hodge--Tate point. Fix a smooth partially proper rigid space $X/C$.
\end{notation}

In the above context, to formulate our main theorem, we shall need the following geometric constructions, both of which are motivated by analogy with similarly named constructions in the formal scheme context; more such analogies are recorded in Table~\ref{table:analogy} at the end.

\begin{enumerate}
\item {\em The analytic prismatization near the Hodge--Tate point:} The paper \cite{ALBRCSRatAnalPrism} (see also \cite{ABLBCSdRFF,LeBrasHabThesis,HauckRatAnalSyn}) introduces an analytic prismatization $(X/\mathcal{Y}_{(0,\infty)})^{\Prism} \in \Fun(\NilPerfd_{\mathcal{Y}_{(0,\infty)}},\Sp)$.  Our aim is to describe its base change to $\mathcal{Y}^\dagger \subset \mathcal{Y}_{(0,\infty)}$, so it is convenient to set
\[ (X/{\mathcal{Y}^\dagger})^{\Prism} := (X/\mathcal{Y}_{(0,\infty)})^{\Prism}|_{\mathcal{Y}^\dagger} \in\Fun(\NilPerfd_{\mathcal{Y}^\dagger},\Sp).\]
(This is the object defined in \cite[Definition 4.5.1 (2)]{ALBRCSRatAnalPrism}, with different notation.) The preceding base change roughly corresponds to the $\varphi$-crystallization procedure (see the discussion following Definition~\ref{def:phicrystallineprism} and Footnote~\ref{fn:AcrysNbhdFontainePrism}), so $(X/{\mathcal{Y}^\dagger})^{\Prism}$ is an analytic counterpart of prismatization over a $\varphi$-crystalline prism. Moreover, this analogy is robust: $(X/{\mathcal{Y}^\dagger})^{\Prism}$ can be checked to be the transmutation of $X$ with respect to the top left term in the obvious nilperfectoid analog of the pullback square \eqref{eq:WmodIsqzero} of animated ring stacks, giving a self-contained definition of $(X/\mathcal{Y}^\dagger)^\Prism$ for our purposes. Write $(X/\mathcal{Y}^\dagger)^\HT$ for its Hodge--Tate reduction, i.e., its restriction along $\mathrm{Spa}(C) \subset \mathcal{Y}^\dagger$; this is the analytic Hodge--Tate stack of Remark~\ref{rmk:AnalyticHTstack}.

\item {\em The analytic $I$-crystallization:} Given a lift $\widetilde{X}/B_2(C)$, there is a stack 
\[ (\widetilde{X}/\mathcal{Y}^\dagger)^{I,\dR} \in\Fun(\NilPerfd_{\mathcal{Y}^\dagger},\Sp) \]
controlling $I$-connections on lifts: define it via transmutation from the evident nilperfectoid  analog of \eqref{eq:IconnRingStack}. The Hodge--Tate reduction of this stack is identified with $BT_X^\dagger\{1\}$.
\end{enumerate}

With the above notation, we expect that the statement and the proof of Theorem~\ref{thm:QuantizationPhiCrys} adapt\footnote{In fact, certain statements are {\em simpler} in the analytic setting: the functor $\Ga^\dR$ is valued in ordinary $p$-complete rings, and is Cartier self-dual. On the flipside, the foundations of the theory are more intricate.}, thanks also to robust Cartier duality statements by Rodr\'iguez Camargo \cite{CamargoCartierDual}, to give the following rigid analytic variant (where $\mathcal{D}(-)$ is interpreted as the derived $\infty$-category of solid quasi-coherent sheaves):

\begin{theorem}[Quantized $p$-adic Simpson equivalence for rigid spaces]
\label{thm:QuantNAHTrigid}
Let $X/C$ be a smooth partially proper rigid space with a lift $\widetilde{X}/B_2(C)$. Write $\Omega = T^*X\{-1\}$ for the twisted cotangent bundle of $X$. Then there is a natural $\mathcal{D}((\Omega/{\mathcal{Y}^\dagger})^\dR)$-module structure on $\mathcal{D}( (X/{\mathcal{Y}^\dagger})^{\Prism})$ and $\mathcal{D}( (\widetilde{X}/{\mathcal{Y}^\dagger})^{I,\dR})$, a multiplicative $\mathbf{G}_m$-gerbe $\mathcal{P}_{\widetilde{X}} \to (\Omega/{\mathcal{Y}^\dagger})^\dR$, and a natural identification of $\mathcal{D}( (\Omega/{\mathcal{Y}^\dagger})^\dR)$-modules
    \[ \mathcal{D}( (X/{\mathcal{Y}^\dagger})^\Prism) \simeq
    \mathcal{D}( (\widetilde{X}/{\mathcal{Y}^\dagger})^{I,\dR})_{\mathcal{P}_{\widetilde{X}}},\]
    where the subscript $\mathcal{P}_{\widetilde{X}}$ indicates twisting by the invertible $\mathcal{D}( (\Omega/{\mathcal{Y}^\dagger})^\dR)$-module attached to $\mathcal{P}_{\widetilde{X}}$.
\end{theorem}

\begin{remark}
\label{rmk:QuantNAHTrigidHT}
On the Hodge--Tate locus, Theorem~\ref{thm:QuantNAHTrigid}  is closely related to Theorem~\ref{thm:padicSimpson}. More precisely, restricting to dualizable objects on the Hodge--Tate locus in Theorem~\ref{thm:QuantNAHTrigid} turns the LHS into $\Perf((X/\mathcal{Y}^\dagger)^\HT)$ and the RHS into $\mathcal{D}_{X\text{-}\Perf}(T^*X\{-1\})_{\mathcal{P}_X}$. The argument in the last sentence of the proof sketch of Theorem~\ref{thm:padicSimpson} then identifies the LHS with $\Perf(X_\proet,\widehat{\calO}_X)$.
\end{remark}

Let us explain how to deduce the untwisted Theorem~\ref{thm:untwistedDeformedSimpson} from the twisted Theorem~\ref{thm:QuantNAHTrigid} above.

\begin{proof}[Proof sketch for Theorem~\ref{thm:untwistedDeformedSimpson}]
We shall obtain the resulting by completing Theorem~\ref{thm:QuantNAHTrigid} at the Hodge--Tate point and passing to dualizable objects. By Remark~\ref{rmk:QuantNAHTrigidHT} (plus deformation theory along $B_\dR^+ \to C$), the LHS of the result is $\Perf(X_\proet,\mathbb{B}_\dR^+)$. The RHS is a twisted version of the RHS of Theorem~\ref{thm:untwistedDeformedSimpson}, so it remains to explain why the twist vanishes. Following the argument over the Hodge--Tate locus in the proof of Theorem~\ref{thm:padicSimpsonUntwisted}, we must show that $\mathcal{P}_{\widetilde{X}}$ trivializes on $(Z/\mathcal{Y}^\dagger)^\dR$, naturally in $Z$ for any $Z \subset \Omega$ finite over $X$, and compatibly with its multiplicative structure  (under the assumptions in Theorem~\ref{thm:untwistedDeformedSimpson}).  Using the description of $\mathcal{P}_X$ appearing in the proof of Theorem~\ref{thm:padicSimpsonUntwisted}, one shows that the choice of $\widetilde{X}$ provides a reduction of structure group to $\mu_{p^\infty} \subset \mathbf{G}_m$ for the multiplicative gerbe $\mathcal{P}_X$. By the topological invariance of the \'etale site, this multiplicative $\mu_{p^\infty}$-gerbe lifts uniquely to one on $(\Omega/\mathcal{Y}^\dagger)^\dR$. Using $\mathbf{A}^1$-homotopy invariance for de Rham cohomology, one shows that the  $\mathbf{G}_m$-gerbe associated this last $\mu_{p^\infty}$-gerbe on $(\Omega/\mathcal{Y}^\dagger)^\dR$ identifies with $\mathcal{P}_{\widetilde{X}}$, and there is a unique such identification extending the tautological one after pullback to $\Omega$. Thanks to this compatibility, we can now simply invoke the trivializations provided by Theorem~\ref{thm:padicSimpsonUntwisted} to finish.
\end{proof}

The following table of analogies might be helpful in translating statements from formal schemes to rigid spaces and vice versa:

\begin{table}[h]
\centering
\renewcommand{\arraystretch}{2.2}
\begin{tabular}{|>{\itshape\raggedright\arraybackslash}p{3.2cm}|>{\centering\arraybackslash}p{5.2cm}|>{\centering\arraybackslash}p{6.6cm}|}
\hline
 & Prismatization over $\varphi$-crystalline prisms $(A,I)$ & Analytic prismatization over $\mathcal{Y}^\dagger$ \\
\hline
Base & $\Spf(\overline{A}) \subset \Spf(A)$ & $\mathrm{Spa}(C) \subset {\mathcal{Y}^\dagger} := \left(\mathrm{Spa}(C) \subset \mathcal{Y}_{(0,\infty)}\right)^\dagger$  \\
Test objects &  $\CAlg^\an_{\Spf(A)}$ (and variants) & $\NilPerfd_{\mathcal{Y}^\dagger}$ (and variants) \\
Small neighbourhood of $0 \in \Ga$ & $\Ga^\sharp \to \Ga$ & $\Ga^\dagger \hookrightarrow \Ga$ \\
Cartier dual of the small neighbourhood & $\widehat{\Ga} \leftarrow \widehat{\Ga^\sharp}$ & $\Ga \leftarrow \Ga^\dagger$ \\
Space of interest &  $X/\Spf(\overline{A})$ smooth scheme & $X/C$ a smooth partially proper rigid space \\
Prismatization \& \par its HT reduction & $(X/A)^\Prism$ \par $(X/A)^\HT \to X$, a $BT_X^\sharp\{1\}$-torsor & $(X/{\mathcal{Y}^\dagger})^{\Prism}$ \par $(X/\mathcal{Y}^\dagger)^{\HT} \to X$, a $BT_X^\dagger\{1\}$-torsor \\
Simpson gerbe & $\mathcal{P}_X \to \widehat{\Omega} := \widehat{T^*X}\{-1\}$ & $\mathcal{P}_X \to \Omega := T^*X\{-1\}$ \\
Simpson equivalence & $\mathcal{D}((X/A)^\HT) \simeq \mathcal{D}(\widehat{\Omega})_{\mathcal{P}_X}$ \par Theorem~\ref{thm:HTstackformal} (3) & $\Perf(X_\proet,\widehat{\calO}_X) \simeq \mathcal{D}_{X\text{-}\Perf}(\Omega)_{\mathcal{P}_X}$ \par (in fact, $\mathcal{D}((X/\mathcal{Y}^\dagger)^{\HT}) \simeq \mathcal{D}(\Omega)_{\mathcal{P}_X}$) \par Theorem~\ref{thm:padicSimpson} (and \cite{CamargoCartierDual}) \\
Crystallization & $(X/A)^\crys$ & $(X/{\mathcal{Y}^\dagger})^\dR$ \\
Lifts & $\widetilde{X}/\Spf(A/I^2)$ & $\widetilde{X}/\mathrm{Spa}(B_2(C))$ \\
$I$-crystallization \& \par its HT reduction & $(\widetilde{X}/A)^{I,\crys}$ \par $BT_X^\sharp\{1\}$ & $(\widetilde{X}/{\mathcal{Y}^\dagger})^{I,\dR}$ \par $BT_X^\dagger\{1\}$ \\
Quantized Simpson gerbe & $\mathcal{P}_{\widetilde{X}} \to (\widehat{\Omega}/A)^\crys$ & $\mathcal{P}_{\widetilde{X}} \to (\Omega/{\mathcal{Y}^\dagger})^\dR$ \\
Quantized Simpson equivalence & $\mathcal{D}((X/A)^\Prism) \simeq \mathcal{D}((\widetilde{X}/A)^{I,\crys})_{\mathcal{P}_{\widetilde{X}}}$ \par Theorem~\ref{thm:QuantizationPhiCrys} & $\mathcal{D}((X/\mathcal{Y}^\dagger)^\Prism) \simeq \mathcal{D}((\widetilde{X}/\mathcal{Y}^\dagger)^{I,\dR})_{\mathcal{P}_{\widetilde{X}}}$ \par Theorem~\ref{thm:QuantNAHTrigid} \\
\hline
\end{tabular}
\caption{Analogies between prismatization over $\varphi$-crystalline prisms and analytic prismatization over $\mathcal{Y}^\dagger$ for the Simpson equivalence and its quantization}
\label{table:analogy}
\end{table}

\clearpage
\bibliographystyle{amsalpha}
\bibliography{mybib}
\end{document}